\documentclass{amsart}

\usepackage{amssymb}
\usepackage{graphicx}

\newcounter{rmkakkocounter}
\newcommand{\setrmkakko}{\setcounter{rmkakkocounter}{1}}
\newcommand{\rmkakko}{{\rm (\thermkakkocounter)} \addtocounter{rmkakkocounter}{1}}

\newtheorem{theorem}{Theorem}[section]
\newtheorem{lemma}[theorem]{Lemma}
\newtheorem{proposition}[theorem]{Proposition}
\newtheorem{proposition_definition}[theorem]{Proposition-Definition}

\newtheorem{corollary}[theorem]{Corollary}
\newtheorem{conjecture}[theorem]{Conjecture}

\newtheorem{_propositionrmkakko}[theorem]{Proposition}

\newtheorem{_definition}[theorem]{Definition}
\newenvironment{definition}{\begin{_definition}\rm}{\end{_definition}}

\newtheorem{_remark}[theorem]{\it Remark}
\newenvironment{remark}{\begin{_remark}\rm}{\end{_remark}}

\newtheorem{_example}[theorem]{Example}
\newenvironment{example}{\begin{_example}\rm}{\end{_example}}

\numberwithin{equation}{section}
\numberwithin{table}{section}
\numberwithin{figure}{section}

\newcommand{\A}{\mathord{\mathbb A}}
\newcommand{\C}{\mathord{\mathbb C}}
\renewcommand{\P}{\mathord{\mathbb  P}}
\newcommand{\R}{\mathord{\mathbb R}}
\newcommand{\mathbbS}{\mathord{\mathbb S}}
\newcommand{\T}{\mathord{\mathbb T}}
\newcommand{\Z}{\mathord{\mathbb Z}}

\newcommand{\CCC}{\mathord{\mathcal C}}
\newcommand{\GGG}{\mathord{\mathcal G}}
\newcommand{\LLL}{\mathord{\mathcal L}}
\newcommand{\MMM}{\mathord{\mathcal M}}
\newcommand{\NNN}{\mathord{\mathcal N}}
\newcommand{\OOO}{\mathord{\mathcal O}}
\newcommand{\RRR}{\mathord{\mathcal R}}
\newcommand{\SSS}{\mathord{\mathcal S}}
\newcommand{\UUU}{\mathord{\mathcal U}}
\newcommand{\YYY}{\mathord{\mathcal Y}}

\font\mathgot=eufm10

\newcommand{\SSSS}{\mathord{\hbox{\mathgot S}}}

\newcommand{\maprightsp}[1]{\; \smash{\mathop{\; \longrightarrow \; }\limits\sp{#1}}\; }

\newcommand{\isomarrow}{\maprightsp{\sim}}

\newcommand{\mapleftspsb}[2]{\; \smash{\mathop{\; \longleftarrow \; }\limits\sp{#1}\limits\sb{#2}}\; }

\newcommand{\mapdown}{\phantom{\Big\downarrow}\hskip -8pt \downarrow}

\newcommand{\mapdownleftright}[2]{\rlap{$\vcenter{\hbox{$\scriptstyle#1$}}$}%
\,\,\mapdown\rlap{$\vcenter{\hbox{$\scriptstyle#2$}}$}}
\newcommand{\mapdownsurj}{
\hbox{$\bigm\downarrow$}
\llap{\hbox{\raise 2pt\hbox{$\bigm\downarrow$}}}%
\vstrechmapdown
}

\newcommand{\inj}{\hookrightarrow}

\newcommand{\surj}{\mathbin{\to \hskip -7pt \to}}

\newcommand{\isom}{\mathbin{\,\raise -.6pt\rlap{$\to$}\raise 3.5pt%
\hbox{\hskip .3pt$\mathord{\sim}$}\,}}

\newcommand{\set}[2]{\{\; {#1} \; \mid \; {#2} \;  \}}
\newcommand{\shortset}[2]{\{ {#1} \,|\, {#2}   \}}

\newcommand{\map}[3]{ #1 \, : \, #2 \, \to \, #3}
\newcommand{\mapisom}[3]{ #1 \, : \, #2 \; \isom \; #3}
\newcommand{\mapinj}[3]{ #1 \, : \, #2 \; \inj \; #3}
\newcommand{\mapsurj}[3]{ #1 \, : \, #2 \; \surj \; #3}

\newcommand{\shortmap}[3]{ #1  :  #2 \to #3}
\newcommand{\shortmapinj}[3]{ #1  :  #2 \inj #3}

\newcommand{\sprime}{\sp\prime}

\newcommand{\spprime}{\sp{\prime\prime}}

\newcommand{\sptimes}{\sp{\times}}

\newcommand{\dual}{\sp{\vee}}

\newcommand{\semidirectproduct}{\rtimes}

\newcommand{\inv}{\sp{-1}}

\newcommand{\Ker}{\operatorname{\rm Ker}\nolimits}

\renewcommand{\Im}{\operatorname{\rm Im}\nolimits}

\newcommand{\Aut}{\operatorname{\rm Aut}\nolimits}

\newcommand{\pr}{\operatorname{\rm pr}\nolimits}
\newcommand{\Sing}{\operatorname{\rm Sing}\nolimits}

\newcommand{\Pt}{\P^2}
\newcommand{\pione}{\pi_1}

\newcommand{\ol}[1]{\overline{#1}}

\newcommand{\rmand}{\textrm{and}}

\newcommand{\rmwhere}{\textrm{where}}

\newcommand{\quand}{\quad\rmand\quad}

\newcommand{\erase}[1]{}

\newcommand{\sphere}{\mathbbS}

\newcommand{\PStextplot}[3]{\rlap{\hskip 0pt  \raise 0pt \hbox{\hskip #1pt  \raise #2pt \hbox{#3}}}}

\newcommand{\ratmap}{\cdot\hskip -2.2pt \cdot \hskip -2.8pt \to}

\newcommand{\PN}{\P^N}
\newcommand{\PM}{\P^M}

\newcommand{\Grass}{\mathord{\rm Gr}} 
\newcommand{\U}{\mathord{U}}
\newcommand{\ngen}[1]{\langle\hskip -1.2pt\langle #1 \rangle\hskip -1.2pt\rangle}
\newcommand{\gen}[1]{\langle #1 \rangle}
\newcommand{\ZQ}{/\hskip -2.2pt/\hskip 1pt}
\newcommand{\Cinf}{\CCC^\infty}

\newcommand{\SB}{\mathord{ SB}}
\newcommand{\B}{\mathord{ B}}
\newcommand{\Div}{\mathord{\rm Div}}
\newcommand{\rDiv}{\mathord{\rm rDiv}}
\newcommand{\Pic}{\mathord{\rm Pic}}
\newcommand{\symgroup}{\mathord{\hbox{\mathgot S}}}

\newcommand{\Fb}{F_b}

\newcommand{\cond}[1]{{\rm (#1)}}

\newcommand{\Xc}{X\sp{\circ}}

\newcommand{\bdr}{\partial}
\newcommand{\IbI}{(I, \bdr I)}

\newcommand{\spc}{\sp{\circ}}
\newcommand{\spsh}{\sp{\sharp}}

\newcommand{\tlb}{\tilde{b}}

\newcommand{\lift}[1]{\tilde{#1}}
\newcommand{\Lift}[1]{\widetilde{#1}}

\newcommand{\splift}{^{\sim}}

\newcommand{\phantominv}{\sp{\phantom{-1}}}

\newcommand{\barlambda}{\bar\lambda}
\newcommand{\red}{\mathord{\rm red}}

\newcommand{\spred}{\sp{\red}}

\newcommand{\unitdisc}{\Delta}
\newcommand{\cunitdisc}{\bar\Delta}

\newcommand{\FL}{{\mathord{\rm FL}}}

\newcommand{\Conj}{{\mathord{\rm C}}}

\renewcommand{\Lift}[1]{(#1)^{\sim}}

\newcommand{\fspshinv}{(f^{\sharp})\inv}

\newcommand{\vexp}{\varepsilon}

\newcommand{\bdre}{\bdr_{\vexp}}

\newcommand{\rightya}{\raise -2.48pt \hbox{\hskip -3.6pt $>$}}
\newcommand{\leftya}{\raise -2.48pt \hbox{\hskip -3.6pt $<$}}
\newcommand{\upya}{\raise -2.3pt \hbox{\hskip -3.3pt $\wedge$}}
\newcommand{\downya}{\raise -2.3pt \hbox{\hskip -3.3pt $\vee$}}

\newcommand{\kuromaru}{{\hskip -2.3pt \lower 2.1pt \hbox{$\bullet$}}}

\newcommand{\PNdual}{(\PN)\dual}

\newcommand{\II}{\mathord{\bf I}}

\begin{document}

\title[Zariski-van Kampen theorem and Grassmannian dual varieties]%
 {Generalized Zariski-van Kampen theorem and \\ its application 
 to Grassmannian dual varieties}

\author{Ichiro Shimada}
\address{
Department of Mathematics,
Graduate School of Science,
Hiroshima University, Kagamiyama, Higashi-Hiroshima, 739-8526,  JAPAN
\\
fax: 81-(0)82-424-0710}

\dedicatory{Dedicated to the memory of Professor Nguyen Huu Duc}
\email{shimada@math.sci.hiroshima-u.ac.jp}

\thanks{Partially supported by
 JSPS Core-to-Core Program 18005:
``New Developments of Arithmetic Geometry, Motive, Galois Theory, and Their Practical Applications"}


\begin{abstract}
We formulate and prove a generalization of Zariski-van Kampen theorem 
on the topological fundamental groups of  smooth complex algebraic varieties.
As  an application, we prove 
a hyperplane section theorem of Lefschetz-Zariski-van Kampen type
for  the fundamental groups of
the complements to the Grassmannian dual varieties.
\end{abstract}

\maketitle

\section{Introduction}\label{sec:Introduction}
We work over the complex number field $\C$.
By a \emph{variety},
we mean a reduced irreducible quasi-projective scheme.
The fundamental group $\pione(V)$ of a variety $V$
is the topological fundamental group of the analytic space
underlying $V$.
The conjunction of paths is read from left to right;
that is, for paths  $\alpha: I:=[0, 1]\to V$ and $\beta: I\to V$,
we define   $\alpha\beta: I\to V$ only when $\alpha(1)=\beta (0)$.
\par
%
For a subset $S$ of a group $G$,
we  denote by $\gen{S}$ the subgroup of $G$ generated by the elements of $S$.
Let a group $\Gamma$ act on $G$ from the right.
Then the subgroup
$$
N_\Gamma:=\gen {\;\shortset{g\inv g\sp\gamma}{g\in G, \gamma\in \Gamma}\;}
$$ 
of $G$ is normal, because
$h\inv (g\inv g^\gamma) h =((gh)\inv (gh)\sp\gamma) (h\inv h\sp\gamma)\inv$.
We then put
$$
G\ZQ \Gamma :=G/N_\Gamma,
$$
and call $G\ZQ \Gamma$ the \emph{Zariski-van Kampen quotient }of $G$ by $\Gamma$.
\par
\bigskip
Let $\shortmap{f}{X}{Y}$ be a dominant morphism
from a smooth variety $X$ to a smooth variety $Y$
with a connected general fiber.
There exists a non-empty Zariski open  subset $Y\spc \subset Y$
such that $f$ is  locally trivial in the $\Cinf$-category over $Y\spc$.
We put $X\spc:=f\inv (Y\spc)$,
and denote by $\shortmap{f\spc}{X\spc}{Y\spc}$ the restriction of $f$ to $X\spc$.
We choose a base point $b\in Y\spc$, put $F_b:=f\inv (b)$,
and choose a base point $\tlb\in F_b$.

We investigate the kernel of
the homomorphism
$$
\map{\iota_*}{\pione (F_b, \tlb)}{\pione (X, \tlb)}
$$
induced by the inclusion $\iota: F_b\inj X$.
The classical Zariski-van Kampen theorem,
which started from~\cite{vanKampen},
describes $\Ker (\iota_*)$ in terms of the monodromy action of $\pione (Y\spc, b)$ on
$\pione (F_b, \tlb)$
\emph{under  the assumption that  a cross-section of $f$ passing through $\tlb$ exists}. 
(See~\cite{MR0366922} for an account of the proof.) 
The cross-section plays a double role;
one is to define the monodromy action of $\pione (Y\spc, b)$ on
$\pione (F_b, \tlb)$,
and the other is to prevent $\pi_2 (Y)$
from contributing  to $\Ker (\iota_*)$.
However,  the cross-section rarely exists in applications.
If we do not have any  cross-section,
then the monodromy of $\pione(Y\spc, b)$ on $\pione (F_b)$ is not well-defined, 
and moreover $\pi_2 (Y)$ may contribute to $\Ker (\iota_*)$.
(See Example~\ref{example:L}.)
\par
In this paper,
we give a generalization of Zariski-van Kampen theorem~(Theorem~\ref{thm:ZvK}),
which describes $\Ker (\iota_*)$
under weaker  conditions on the existence of the cross-section.
Informally,
our  theorem
states that,
if there exists a cross-section  on a subspace of $Y$ whose $\pi_2$  surjects to $\pi_2(Y)$, 
then,
under additional  assumptions  on the  singular fibers  of $f$,
 $\Ker (\iota_*)$ is generated by the monodromy relations arising from  the \emph{lifted monodromy},
which is defined as follows.
\par
Since $\shortmap{f\spc}{X\spc}{Y\spc}$ is locally trivial, 
the groups $\pione (f\inv (f(x)), x)$ form a locally constant system on $X\spc$
when $x$ moves on $X\spc$,
and hence 
  $\pione (X\spc, \tlb)$ acts on $\pione (F_b, \tlb)$ from the right 
  in a natural way.
We denote this action by
\begin{equation}\label{eq:mu}
\map{\mu}{\pione (\Xc, \tlb)}{\Aut(\pione (\Fb, \tlb))},
\end{equation}
and call $\mu$ the \emph{lifted monodromy}.
\par
\medskip
Combining  our main result  with Nori's lemma~\cite{MR732347} (see Proposition~\ref{prop:nori}),
we obtain the following:
\begin{corollary}\label{cor:RRReq}
Suppose that the following three conditions are satisfied:
\begin{itemize}
\item[\cond{C1}]
the locus $\Sing (f)$ of critical points of $f$ is of codimension $\ge 2$
in $X$,
\item[\cond{C2}]
there exists  a Zariski closed subset $\Xi_0$ of $Y$ with codimension $\ge 2$
such that 
 $F_y:=f\inv (y)$ is non-empty and irreducible
 for any $y\in Y\setminus \Xi_0$, and
\item[\rlap{\cond{Z}}\phantom{\cond{C2}}]
there exist a subspace $Z\subset Y$ containing $b$ 
and a continuous cross-section $s_Z: Z\to f\inv (Z)$ 
of $f$ over $Z$ satisfying
$s_Z(Z)\cap \Sing (f)=\emptyset$ and $s_Z(b)=\tlb$
such that
 the inclusion $Z\inj Y$ induces 
a surjection  $\pi_2 (Z, b)\surj \pi_2(Y, b)$.
\end{itemize}
Let $i_{X*}: \pione(X\spc, \tlb)\to \pione (X, \tlb)$ be the homomorphism
induced by the inclusion $i_{X}: X\spc\inj X$.
Then $\Ker (\iota_*)$ is equal to
\begin{equation}\label{eq:RRR}
\RRR:=\gen{\;\shortset{g\inv g^{\mu(\gamma)}}{g\in \pione(F_b, \tlb),\; \gamma\in \Ker (i_{X*})}\;},
\end{equation}
and   we have the exact sequence
$$
1\;\maprightsp{}\; 
\pione(F_b, \tlb)\ZQ \Ker (i_{X*})\;\maprightsp{\iota_*}\;
\pione(X, \tlb)\;\maprightsp{f_*}\;
\pione(Y, b)\;\maprightsp{}\;
1.
$$
\end{corollary}
\medskip
\begin{remark}
The condition~\cond{Z} is trivially satisfied
if $\pi_2(Y)=0$;
for example,
when $Y$ is an affine space $\A^N$,
an abelian variety, or a Riemann surface of genus $>0$.
\end{remark}
In our previous papers~\cite{MR1341806}, 
~\cite{MR1988200} 
and~\cite{MR1952329}, 
we have given three different proofs to
 a  special case of  Theorem~\ref{thm:ZvK},   
 where $Y$ is an affine space $\A^N$.
Even this special case has yielded many applications
(\cite{MR1282219, 
MR1354002, 
MR1421396, 
MR1428061, 
MR1474860, 
 MR2011641,  
 MR1952330}). 
Thus we can expect more applications of the generalized Zariski-van Kampen theorem of this paper.
\par
\medskip
As an easy application,
we obtain the following:
\begin{corollary}\label{cor:proj}
Let $\shortmap{f}{X}{Y}$ be a  morphism
from a smooth variety $X$ to a smooth variety $Y$.
Suppose that $\pi_2(Y)=0$,
that $f$ is projective
with the general fiber $F_b$ being  connected,
and that $\Sing(f)$ is of codimension $\ge 3$ in $X$.
Let $\shortmapinj{\iota}{F_b}{X}$ be the inclusion. 
Then the sequence 
$$
1\;\maprightsp{}\;\pione (F_b)\;\maprightsp{\iota_*} \; \pione (X) \;\maprightsp{f_*} \; \pione (Y) \;\maprightsp{}\; 1
$$
is exact.
\end{corollary}
As the next  application,  we investigate 
the fundamental group of the complement
of the \emph{Grassmannian dual variety},
and prove a hyperplane section theorem of
Zariski-Lefschetz-van Kampen type.
\par
A Zariski closed subset of a projective space  $\PN$ is said to be \emph{non-degenerate}
if it is not contained in any hyperplane of $\PN$.
We denote by $\Grass^c(\PN)$
the Grassmannian variety of $(N-c)$-dimensional linear subspaces  
of  $\PN$.
For a point  $t\in (\PN)\dual=\Grass^1(\PN)$ of the dual projective space,
let  $H_t\subset \PN$ denote the corresponding hyperplane.
\par
%
Let $W$ be a closed subscheme of $\PN$
such that every irreducible component is of dimension $n$.
For $c\le n$, 
the \emph{Grassmannian dual variety of $W$ in $\Grass^c(\PN)$}
is defined to be the locus of $L\in \Grass^c(\PN)$
such that 
the scheme-theoretic intersection 
of $W$ and the linear subspace $L\subset\PN $ 
\emph{fails} to be smooth of dimension $n-c$.
For a non-negative integer $k$, 
we denote by $\U_k(W,\PN)$ the complement of the Grassmannian dual variety of $W$ in $\Grass^{n-k}(\PN)$;
that is, 
 $\U_k(W,\PN)\subset \Grass^{n-k}(\PN)$ is the  Zariski open subset 
 of  all $L\in \Grass^{n-k}(\PN)$ that intersect  $W$ along  a smooth scheme  of dimension $k$.
\par
%
%
Let $X\subset \PN$ be a smooth non-degenerate projective  variety of dimension $n\ge 2$.
The fundamental group  $\pione ((\PN)\dual\setminus X\dual)=\pione (U_{n-1}(X,\PN))$
of the complement of the dual variety 
has been studied in  several papers~(for example,~\cite{MR1682991, MR644816}).
However, there seem to be few studies on its generalization to Grassmannian varieties.
We will investigate the fundamental groups 
 $\pione (U_k(X,\PN))$ for $k=0, \dots, n-2$.
 \par
We choose a \emph{general} line $\Lambda $ in  $(\PN)\dual$,
and consider  the corresponding  pencil  $\{H_{t}\}_{t\in \Lambda}$ of hyperplanes.
Let $A:=\bigcap H_t\cong \P^{N-2}$ denote 
the axis of the pencil.
We  put
$$
Y_t:=X\cap H_t \quand Z_{\Lambda}:=X\cap A.
$$
Let $k$ be an integer such that $0\le k\le n-2$.
Regarding $\Grass^{c-1} (H_t)$ as a closed subvariety of $\Grass^c(\PN)$,
and $\Grass^{c-2} (A)$ as a closed subvariety of $\Grass^{c-1}(H_t)$,
we have  canonical   inclusions 
$$
\U_k (Z_\Lambda, A)\;\;\inj\;\; \U_k(Y_t,H_t)\;\;\inj\;\; 
\U_k(X,\PN).
$$
Since $k\le n-2$,  the space $\U_k (Z_\Lambda, A)$ is non-empty.
(When $k=n-2$,  the space $\U_{n-2} (Z_\Lambda, A)$ is equal to
the one-point set  $\Grass^0(A)=\{A\}$.)
We choose a base point
$$
L_o\;\;\in\;\; \U_k (Z_\Lambda, A),
$$
which serves also as a base point of $\U_k(X,\PN)$ and of $\U_k(Y_t,H_t)$
by the natural inclusions above.
Consider the space 
$$
\UUU_k (X,\PN,\Lambda):=\set{(L, t)\in \U_k(X,\PN)\times \Lambda}{L\subset H_t}
$$
with  the projection  
$$
\map{f_{\Lambda}}{\UUU_k (X,\PN,\Lambda)}{\Lambda}.
$$
The fiber of $f_{\Lambda}$ over $t\in \Lambda$ is canonically identified with 
$\U_k(Y_t, H_t)$, and 
the point $L_o$ furnishes  us with  a holomorphic section
$$
\map{s_o}{\Lambda}{\UUU_k (X,\PN,\Lambda)}
$$
of $f_{\Lambda}$.
There exists a proper Zariski closed subset $\Sigma_\Lambda$ of $\Lambda$
such that $f_{\Lambda}$ is locally trivial over $\Lambda\setminus \Sigma_{\Lambda}$
in the $\Cinf$-category.
We choose a base point  $0\in \Lambda\setminus \Sigma_\Lambda$.
By the section $s_o$,
the fundamental group $\pione (\Lambda\setminus \Sigma_{\Lambda}, 0)$
acts on $\pione (U_k(Y_0, H_0), L_o)$
in  the classical (not lifted) monodromy.
\par
\medskip
Using the fact that $\Lambda\inj (\PN)\dual$ induces an isomorphism
$\pi_2(\Lambda)\cong \pi_2((\PN)\dual)$,
we derive from Theorem~\ref{thm:ZvK} the following:
\begin{theorem}\label{thm:ULZvK}
\setrmkakko
Consider the homomorphism
$$
\map{\iota_*}{\pione (\U_k (Y_0, H_0), L_o)}{\pione (\U_k (X, \PN), L_o)}
$$
induced by the inclusion $\iota: \U_k (Y_0, H_0)\inj \U_k (X, \PN)$.

\rmkakko
If $k\le n-2$, then $\iota_*$ is surjective and induces an isomorphism
$$
\pione (\U_k (Y_0, H_0), L_o)\ZQ \pione (\Lambda\setminus \Sigma_{\Lambda}, 0)
\;\;\isom\;\; \pione (\U_k (X, \PN), L_o).
$$

\rmkakko 
If $k< n-2$,   the monodromy action of $\pione (\Lambda\setminus \Sigma_{\Lambda}, 0)$
on $\pione (\U_k (Y_0, H_0), L_o)$ is trivial.
In particular,
the homomorphism $\iota_*$ is an isomorphism
for $k< n-2$.
\end{theorem}
Remark that this theorem resembles the classical Lefschetz hyperplane section theorem on the homotopy groups
of smooth projective varieties:
namely, the inclusion $Y_0\inj X$
induces  surjective homomorphisms
$\pi_k(Y_0)\surj \pi_k(X)$
for $k\le n-1$,
and isomorphisms
$\pi_k(Y_0)\isom  \pi_k(X)$
for $k< n-1$.
\par
\medskip
The isomorphism in the assertion (2) of Theorem~\ref{thm:ULZvK} 
seems to fail to hold for $k=n-2$, as can be seen from the argument in \S\ref{sec:ADKY} of this paper.
\par
\medskip
As the third application,  we study
$\pione (\U_k (X, \PN), L_o)$
for  $k=0$.
By Theorem~\ref{thm:ULZvK},
it is enough to investigate   the case where $\dim X=2$, and 
to  study the monodromy action of $\pione (\Lambda\setminus \Sigma_{\Lambda}, 0)$
on $\pione (\U_0 (Y_0, H_0),  L_o)$,
where $Y_0=X\cap H_0$ is a smooth compact Riemann surface.
\par
\medskip
First we define the simple braid group $\SB^d_g $
of $d$ strings on a compact Riemann surface $C$ of genus $g>0$.
We denote by $\Div^d (C)$
the variety of effective divisors of degree $d$ on $C$,
and by $\rDiv^d(C)\subset \Div^d(C)$
the Zariski open subset consisting of \emph{reduced}  divisors.
We fix a base point
$$
D_0=p_1+\dots+p_d
$$
of $\rDiv^d (C)$.
The braid group
$\B^d_g =\B(C, D_0)$ is defined to be the fundamental group
$\pione (\rDiv^d(C), D_0)$.
(See~\cite{MR0375281}.)
\begin{definition}\label{def:SB}
The \emph{simple braid group}
$\SB^d_g =\SB(C, D_0)$ is defined to be the kernel
of the homomorphism
$\B(C, D_0)\to \pione (\Div^d(C),D_0)$
induced by the inclusion $\rDiv^d(C)\inj \Div^d(C)$.
\end{definition}
Let $\MMM^d_g=\MMM(C, D_0)$
be the topological group 
of orientation-preserving diffeomorphisms
$\gamma$ of $C$  acting from the right  that satisfy  ${p_i}^\gamma =p_i$ 
for each point $p_i$ of $D_0$.
We denote by 
$$
\varGamma^d_g=\varGamma(C, D_0):=\pi_0(\MMM(C, D_0))
$$
the group of isotopy classes of diffeomorphisms in $\MMM^d_g=\MMM(C, D_0)$,
which acts on 
$\SB^d_g =\SB(C, D_0)$ from the right   in a natural way.
\par
\medskip
Let $C\subset \PM$ be a smooth non-degenerate projective curve of degree $d$ and genus $g>0$,
and let $D_0\in \rDiv^d(C)$ be  a general  hyperplane section.
We will investigate $\pione (\U_0 (C, \PM), D_0)$;
that is, the fundamental group of the complement of the \emph{dual hypersurface} of $C$.
\par
\medskip
In~\cite{KulikovMPIpreprint} and~\cite{MR1988200},
we studied 
this group 
  under conditions that
$d\ge 2g+2$ and that the invertible sheaf  $\OOO_C(D_0)$
corresponds to a \emph{general} point of the Picard variety $\Pic^d(C)$
of isomorphism classes of line bundles of degree $d$. 
\par
\medskip
Using the fact that $\pi_2(\Pic^d(C))=0$,
we derive from our main theorem (Theorem~\ref{thm:ZvK})  
the following result, which states the same result as in~\cite{KulikovMPIpreprint} and~\cite{MR1988200}
under weaker conditions.
\begin{definition}
We say that $C\subset \P^M$ is \emph{Pl\"ucker general}
if the dual curve $\rho(C)\dual\subset (\Pt)\dual$ of the image 
$\rho(C)\subset \P^2$ of the  general projection
$\rho: C\to \Pt$ has only ordinary nodes and ordinary cusps 
as its singularities.
\end{definition}
\begin{theorem}\label{thm:SB}
Suppose that $d\ge g+4$  and that $C$ is Pl\"ucker general in $\PM$.
Then $\pione (\U_0 (C, \PM), D_0)$
is isomorphic to $\SB(C, D_0)$.
\end{theorem}
Let $X\subset \PN$ be a smooth non-degenerate projective surface of degree $d$,
and let $\{Y_t\}_{t\in \Lambda}$ be a  pencil of hyperplane sections of $X$
parameterized by a general line $\Lambda\subset(\PN)\dual$
with  the base locus $Z_{\Lambda}:=X\cap A$,
where $A=\bigcap H_t$ is the axis of the pencil $\{H_t\}_{t\in \Lambda}$ of hyperplanes.
Let
$$
\map{\varphi}{\YYY:=\set{(x, t)\in X\times \Lambda}{x\in H_t}}{\Lambda}
$$
be the fibration of the pencil.
Then $\varphi$ is locally trivial over $\Lambda\setminus \Sigma\sprime_{\Lambda}$
in the $\Cinf$-category,
where $\Sigma\sprime_{\Lambda}$ is the set of critical values of $\varphi$.
Let $0$ be a general point of $\Lambda$.
The corresponding member $Y_0$ is a compact Riemann surface of genus 
$$
g:=(d+H_0\cdot K_X)/2 +1.
$$
Note that $\U_0 (Z_{\Lambda}, A)=\{A\}$, and that 
each point of $Z_{\Lambda}$ yields a holomorphic section of $\varphi:\YYY\to\Lambda$.
By the classical monodromy,
we obtain a homomorphism
\begin{equation}\label{eq:monhom}
\pione (\Lambda\setminus \Sigma\sprime_{\Lambda}, 0)\;\to\; \varGamma^d_g=\varGamma(Y_0, Z_{\Lambda}),
\end{equation}
and hence $\pione (\Lambda\setminus \Sigma\sprime_{\Lambda}, 0)$ acts on the simple braid group
$\SB^d_g =\SB(Y_0, Z_{\Lambda})$ from the right.
We denote by 
$$
\Gamma_{\Lambda}\;\subset\;  \varGamma^d_g=\varGamma(Y_0, Z_{\Lambda})
$$
the image of the monodromy homomorphism~\eqref{eq:monhom}.
Combining Theorems~\ref{thm:ULZvK}~and~\ref{thm:SB},
we obtain the following:
\begin{corollary}\label{cor:SB}
Let $X$, $\{Y_t\}_{t\in \Lambda}$, $Z_{\Lambda}=X\cap A$ and  $\Gamma_{\Lambda}$ be as above.
Suppose that $g>0$, $d\ge g+4$,
and that a general hyperplane section of $X$ is Pl\"ucker general.
Then $\pione (\U_0 (X, \PN), A)$ is isomorphic to 
the Zariski-van Kampen quotient $\SB(Y_0, Z_\Lambda)\ZQ \Gamma_{\Lambda}$.
\end{corollary}
A motivation of the study of the fundamental group $\pione (U_0(X, \PN))$
for a surface $X\subset \PN$ is the  conjecture 
 of Auroux,  Donaldson, Katzarkov and Yotov~\cite{MR2081427}
 about
  the fundamental group $\pione (\Pt\setminus B)$ of the complement of the branch curve
$B\subset \Pt$ of the general projection $X\to \Pt$,
which had been intensively  studied by Moishezon, Teicher, Robb.
The  weakening of the conditions 
from our previous works~(\cite{KulikovMPIpreprint}, ~\cite{MR1988200}) to 
the present result (Theorem~\ref{thm:SB}) is  important 
with respect to this application.
See Remark~\ref{rem:LSXm}.
\par
\bigskip
The plan of this paper is as follows.
In~\S\ref{sec:ZvKQ},
we state  some elementary facts about Zariski-van Kampen quotients.
In~\S\ref{sec:pionefib},
we prove the generalized Zariski-van Kampen theorem~(Theorem~\ref{thm:ZvK}).
We then prove its variant~~(Theorem~\ref{thm:C}),
 and  deduce Corollaries~\ref{cor:RRReq}~and~\ref{cor:proj}.
The main ingredient of the proof is the notion of
\emph{free loop pairs of monodromy relation type}
(Definitions~\ref{def:frp} and~\ref{def:frpmrt}),
and Proposition~\ref{prop:TD}.
Using these results,
we prove Theorem~\ref{thm:ULZvK} in~\S\ref{sec:proof1},
and  Theorem~\ref{thm:SB} in~\S\ref{sec:SB}.
In the last section, we explain the relation between 
$\pione (U_0(X, \PN))$
and the conjecture of Auroux,  Donaldson, Katzarkov, Yotov.%
\par
\bigskip
This paper is dedicated to the memory of Professor Nguyen Huu Duc.
\par
\bigskip
{\bf Conventions and Notation}
\begin{itemize}
\setrmkakko
\item[\rmkakko] The constant map to a point $P$ is denoted by $1_P$.
\item[\rmkakko]
 We denote by $I\subset \R$ the interval $[0,1]$,
 by  $\unitdisc\subset\C$ the  open unit disc, 
and by $\cunitdisc\subset\C$ the  closed  unit disc. 
\item[\rmkakko]
For a continuous map
$\shortmap{\delta}{\cunitdisc}{T}$
to a topological space $T$,
we denote by
$$
\map{\bdre\delta}{I}{T}
$$
the loop   given by $t\mapsto \delta(\exp(2\pi\sqrt{-1}t))$. 
\end{itemize}
\section{Zariski-van Kampen quotient}\label{sec:ZvKQ}
%
%
%
%
\begin{definition}
Let $G$ be a group,
and let $S$ be a subset of $G$.
We denote by $\gen{S}_G$ or simply by $\gen{S}$
the smallest subgroup of $G$ containing $S$, 
and by $\ngen{S}_G$ or simply by $\ngen{S}$
the smallest \emph{normal} subgroup of $G$ containing $S$.
\end{definition}
We let  a group $\Gamma$ act on a group $G$ from the right.
The following are easy:
\begin{lemma}\label{lem:normal}
For any $\gamma\in \Gamma$, the subgroup
$\gen{\shortset{g\inv g\sp{\gamma}}{g\in G}}_G$ of $G$
is normal.
Hence, 
for any subset $\Sigma\subset \Gamma$, 
the subgroup 
$\gen{\shortset{g\inv g\sp\sigma}{g\in G, \sigma\in \Sigma}}_G$ is  normal.
\end{lemma}
\begin{lemma}\label{lem:S}
Let $S$ be a subset of $G$, and let $\Sigma$ be a subset of $\Gamma$.
If $G=\gen{S}_G$ and $\Gamma=\gen{\Sigma}_\Gamma$,
then we have
$$
\ngen{\shortset{s\inv s\sp\sigma}{s\in S, \sigma\in \Sigma}}_G=
\gen{\shortset{g\inv g\sp\sigma}{g\in G, \sigma\in \Sigma}}_G=
\gen{\shortset{g\inv g\sp\gamma}{g\in G, \gamma\in \Gamma}}_G.
$$
\end{lemma}
\begin{definition}\label{def:sdp}
We define $G\semidirectproduct \Gamma$ to be the group
with the underlying set $G\times \Gamma$ and with the product
defined by
$$
(g, \gamma)(h, \delta):=(g \cdot \left(h^{(\gamma\inv)}\right), \gamma\delta).
$$
We then define homomorphisms
$i : G\to G\semidirectproduct \Gamma$, $p: G\semidirectproduct \Gamma\to\Gamma$
and $s:\Gamma\to G\semidirectproduct \Gamma$ by 
$i(g):=(g, 1)$,  $p(g, \gamma):=\gamma$ and $s (\gamma):=(1, \gamma)$.
Then we obtain  an exact sequence
\begin{equation}\label{eq:sdp}
1 \;\;\maprightsp{}\;\; 
G \;\;\maprightsp{i}\;\;
 G\semidirectproduct \Gamma \;\;\maprightsp{p}\;\;
  \Gamma\;\;\maprightsp{} \;\;1 
\end{equation}
with  the cross-section  $s$  of $p$,
and the action $g\mapsto g^\gamma$ of $\gamma\in \Gamma$ on $G$
coincides with the inner-automorphism $g\mapsto s(\gamma)\inv g s (\gamma)$
by $s(\gamma)\in  G\semidirectproduct \Gamma$
on the normal subgroup $G=i(G)$ of  $ G\semidirectproduct \Gamma$.
\end{definition}
The following two lemmas are elementary:
\begin{lemma}\label{lem:GGG}
Let $\GGG$ be a group. 
Suppose that we are given an exact sequence
\begin{equation}\label{eq:GGG}
1\;\;\maprightsp{}\;\; 
G \;\;\maprightsp{i\sprime} \;\;
\GGG \;\;\maprightsp{p\sprime}\;\; 
\Gamma\;\;\maprightsp{}\;\; 1 
\end{equation}
with a cross-section  $s\sprime: \Gamma\to \GGG$ of $p\sprime$
that is a homomorphism of groups.
Suppose also that the action of $\gamma\in \Gamma$ on $g\in G$
is equal to  
the inner-automorphism by $s\sprime(\gamma)$;
that is, we have
$i\sprime(g\sp\gamma)=s\sprime(\gamma)\inv i\sprime (g) s\sprime(\gamma)$
for any $g\in G$ and $\gamma\in \Gamma$.
Then there exists an isomorphism $\GGG\cong G\semidirectproduct \Gamma$
such that 
the exact sequences~\eqref{eq:sdp} and~\eqref{eq:GGG} coincide
and the cross-section $s$ corresponds to $s\sprime$
by this isomorphism.
\end{lemma}
\begin{lemma}\label{lem:sGamma}
The composite homomorphism 
$$
G \;\maprightsp{i}\;
 G\semidirectproduct \Gamma\;\maprightsp{}\; (G\semidirectproduct \Gamma)/ \ngen{s(\Gamma)}_{G\semidirectproduct \Gamma}
$$
is surjective, and its kernel is equal to
$\gen{\shortset{g\inv g^\gamma }{g\in G, \gamma\in \Gamma}}$;
that is, 
the Zariski-van Kampen quotient $G\ZQ\Gamma$ is 
isomorphic to $(G\semidirectproduct \Gamma) / \ngen{s(\Gamma)}$.
\end{lemma}
\section{Fundamental groups of algebraic fiber spaces}\label{sec:pionefib}
Let $X$ and $Y$ be smooth  varieties,
and let $f: X\to Y$ be  a dominant morphism.
We denote by $\Sing (f)\subset X$
the Zariski closed subset of 
the critical points of $f$.
For a point $y\in Y$, we put
$$
F_y:=f\inv (y).
$$
Let $\alpha : T\to Y$ be a continuous map
from a topological space $T$.
Then a continuous map
$\lift{\alpha}: T\to X$
is  said to be  a \emph{lift of $\alpha$} if $f\circ \lift{\alpha}=\alpha$.
\par
\medskip
We fix, once and for all,
a proper Zariski closed subset 
$$
\Sigma\subset Y
$$
such that $f\spc:X\spc\to Y\spc$  is locally trivial in the $\Cinf$-category,
where 
$$
Y\spc :=Y\setminus \Sigma,
\quad 
X\spc :=f\inv (Y\spc)
\quand
f\spc:=f|_{X\spc} : X\spc \to Y\spc.
$$
(In particular, $\Sing (f)$ is contained in $f\inv (\Sigma)$.)
%
It follows from Hironaka's resolution of singularities 
that such a proper Zariski closed subset $\Sigma\subset Y$  exists.
%
We then fix base points
$$
b\in Y\spc\quand
\tlb \in \Fb\subset X\spc,
$$
and consider the homomorphisms
$$
\map{\iota_*}{\pione (\Fb, \tlb)}{\pione(X, \tlb)}
\quand
\map{f_*}{\pione (X, \tlb)}{\pione (Y, b)}
$$
induced by the inclusion $\iota :\Fb\inj X$ and the morphism $f:X\to Y$, respectively.
The aim of Zariski-van Kampen theorem is  to describe  $\Ker (\iota_*)$.
\par
\medskip
The following result of Nori~\cite{MR732347} will be used throughout this paper:
\begin{proposition}\label{prop:nori}
Suppose that $F_b$  is connected,
and that there exists  a Zariski closed subset $\Xi\sprime\subset Y$ of codimension $\ge 2$
such that $F_y\setminus(F_y\cap \Sing(f))\ne\emptyset $ for any $y\in Y\setminus \Xi\sprime$.
Then    $f_*: \pione (X, \tlb)\to\pione(Y, b)$ is surjective,
and its  kernel  
 is equal to the image of
$\iota_*:\pione (\Fb, \tlb)\to \pione (X, \tlb)$.
\end{proposition}
\begin{proof}
See Nori~\cite[Lemma 1.5]{MR732347} and~\cite[Proposition 3.1]{MR1988200}.
\end{proof}
Let $\lift{\alpha} : I\to \Xc$ be a path,
and we put $\alpha:=f\spc \circ\lift{\alpha}$.
Then $\lift{\alpha}$ induces an isomorphism
$\pione (F_{\alpha(0)}, \lift{\alpha}(0))\isom \pione (F_{\alpha(1)}, \lift{\alpha}(1))$,
which depends only on the homotopy class (relative to $\bdr I$)
of the path $\lift{\alpha}$.
Hence we can  write this isomorphism 
as
$$
\mapisom{[\lift{\alpha}]_*}{\pione (F_{\alpha(0)}, \lift{\alpha}(0))}%
{\pione (F_{\alpha(1)}, \lift{\alpha}(1))}.
$$
The \emph{lifted monodromy}
$$
\map{\mu}{\pione (\Xc, \tlb)}{\Aut(\pione (\Fb, \tlb))}
$$
introduced  in \S\ref{sec:Introduction} (see~\eqref{eq:mu})
is obtained by applying this construction to the loops in $\Xc$ with the base point $\tlb$.
By the definition, we have the following:
\begin{proposition}\label{prop:lift}
For any $[\lift{\alpha}]\in \pione(X\spc, \tlb)$ and $g\in \pione (F_b, \tlb)$,
we have 
$$
\iota\spc_*(g^{\mu([\lift{\alpha}])})=[\lift{\alpha}]\inv \cdot \iota\spc_*(g) \cdot [\lift{\alpha}]
$$
in $\pione(X\spc, \tlb)$, where 
$\iota\spc_*:\pione (F_b, \tlb)\to \pione (X\spc, \tlb)$
is the homomorphism induced by  the inclusion $\iota\spc: F_b\inj X\spc$.
\end{proposition}
First we prove the following:
\begin{proposition}\label{prop:relisinKer}
Suppose that a loop $\lift{\alpha}:(I, \bdr I)\to (X\spc, \tlb)$ is null-homotopic in $(X, \tlb)$.
Then $g\inv g^{\mu([\lift{\alpha}])}\in \Ker (\iota_*)$
for any $g\in \pione (\Fb, \tlb)$.
\end{proposition}
\begin{proof}
We put $\alpha:=f\spc\circ \lift{\alpha}$, and 
$\sqcup := (I\times \{0\})\cup (\bdr I\times I)$.
Let $g\in \pione (\Fb, \tlb)$ be represented by a loop
$\gamma :(I, \bdr I)\to (F_b, \tlb)$.
We define $\phi_{\sqcup} : \sqcup \to X\spc$ by
$$
\phi_{\sqcup}(s, 0):=\gamma(s),
\quad
\phi_{\sqcup}(0, t):=\lift{\alpha}(t),
\quand
\phi_{\sqcup}(1, t):=\lift{\alpha}(t).
$$
Then we have
$f\spc \circ \phi_\sqcup=({\alpha}\circ \pr_2)|_{\sqcup}$,
where $\pr_2: I\times I \to I$ is the second projection.
Since $\sqcup$ is a strong deformation retract of $I\times I$ and $f\spc$ is locally trivial,
the extension of $({\alpha}\circ \pr_2)|_{\sqcup}: \sqcup\to Y\spc $ to
${\alpha}\circ \pr_2: I\times I\to Y\spc$ lifts to an extension from $\phi_{\sqcup}:\sqcup \to X\spc$ 
to a continuous map 
$\shortmap{\phi}{I\times I}{X\spc}$
that satisfies  $\phi|_\sqcup=\phi_\sqcup$ and  $f\spc \circ \phi={\alpha}\circ \pr_2$.
(See Figure~\ref{figphi}.)
Then the loop
$$
\map{\gamma\sprime:=\phi|_{I\times\{1\}}}{(I, \bdr I)}{(F_b, \tlb)}
$$
represents $g^{\mu([\lift{\alpha}])}$.
Since $\phi|_{\{0\}\times I}=\lift{\alpha}$ and $\phi|_{\{1\}\times I}=\lift{\alpha}$,
we have
$$
[\gamma]\inv [\lift{\alpha}] [\gamma\sprime][\lift{\alpha}]\inv =1
$$
in $\pione(X\spc, \tlb)$.
Since $[\lift{\alpha}] =1$ in $\pione(X, \tlb)$ by the assumption,
we have $[\gamma]\inv  [\gamma\sprime] =1 $ in $\pione(X, \tlb)$.
\end{proof}
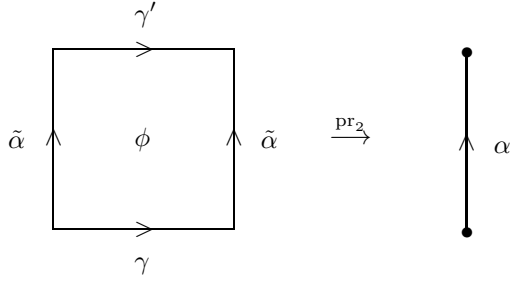
\begin{figure}
\begin{center}
$
\setlength{\unitlength}{1.2mm}
\vcenter{\hbox{
\begin{picture}(30,25)(-5,-3)
\put (0,0) {\line (1,0){20}}\put (10,0){\rightya}
\put (9, -4.4) {$\gamma$}
\put (20,0) {\line (0,1){20}}  \put (20,10){\upya}
\put (23, 9) {${\lift{\alpha}}$}
\put (20,20) {\line (-1,0){20}}  \put (10,20){\rightya}
\put (9, 23.0) {${\gamma\sprime}$}
\put (0,20) {\line (0,-1){20}} \put (0,10){\upya}
\put (-5, 9) {${\lift{\alpha}}$}
\put (9, 9.2) {$\phi$}
\end{picture}
}}
\quad
\vcenter{\hbox{$\maprightsp{\pr_2}$}}
\quad
\vcenter{\hbox{
\begin{picture}(10,30)(-5,-5)
\put (0,0) {\line (0, 1){20}} \put (0,10){\upya}
\put (0,0){$\kuromaru$}
\put (0,20){$\kuromaru$}
\put (3,9) {${{\alpha}}$}
\end{picture}
}}
$
\end{center}
\caption{The extension $\phi$}\label{figphi}
\end{figure}
By  Proposition~\ref{prop:relisinKer},
the normal subgroup  $\RRR$  
defined by~\eqref{eq:RRR} is contained in $\Ker (\iota_*)$.
However $\RRR$ is  not equal to $\Ker (\iota_*)$ in general.
We give two examples.
\begin{example}\label{example:L}
Let $L\to \P^1$ be a line bundle of degree $d>0$,
and let $L\sptimes\subset L$ be the complement of the zero-section.
Since  the projection  $\shortmap{f}{X=L\sptimes}{Y=\P^1}$
is  locally trivial,
we can put $\Sigma=\emptyset$,
and hence   $\RRR=\{1\}$.
However,  the kernel of 
$$
\map{\iota_*}{\pione (F_b)=\pione(\C\sptimes)\cong\Z}{\pione(L\sptimes)\cong \Z/d\Z}
$$
is non-trivial.
Indeed, $\Ker(\iota_*)$  is equal to the image of the boundary homomorphism
$\pi_2(\P^1)\to \pione(\C\sptimes)$
in the homotopy exact sequence.
\end{example}
\begin{example}\label{example:nonirred}
Consider the morphism
$$
\map{f}{X=\C^2}{Y=\C}
$$
given by $f(x, y):=xy$.
We can put
$\Sigma=\{0\}$, 
 and hence the fundamental group of 
 $X\spc=\C^2\setminus\{xy=0\}$ 
is  isomorphic to $\Z^2$.
 The general fiber $F_b$ is isomorphic to $\P^1$ minus two points, and 
 the lifted monodromy action of $\pione( X\spc)$ on $\pione (F_b)\cong \Z$ is
 trivial.
 Therefore we have $\RRR=\{1\}$,
 while we have $\Ker (\iota_*)=\pione (F_b)\cong\Z$.
\end{example}
Our ultimate goal is to show that the three conditions in Corollary~\ref{cor:RRReq}
is sufficient for $\RRR=\Ker(\iota_*)$ to hold.
\par
\medskip
From now on,  we suppose that $f:X\to Y$ satisfies
the first two of the three conditions in Corollary~\ref{cor:RRReq}; namely, we assume the following:
\begin{itemize}
\item[\cond{C1}]
$\Sing (f)$ is of codimension $\ge 2$ in $X$, and 
\item[\cond{C2}]
there exists a Zariski closed subset
$\Xi_0\subset Y$  of codimension $\ge 2$ such that 
 $F_y$ is non-empty and irreducible for any $y\in Y\setminus \Xi_0$.
\end{itemize}
\begin{remark}\label{rem:C0C3}
By the conditions~\cond{C1} and~\cond{C2}, the following hold:
\begin{itemize}
\item[\cond{C0}]
for $y\in Y\spc$, the  fiber  $F_y$ is  connected, and 
\item[\cond{C3}]
there exists  a Zariski closed subset $\Xi_1\subset Y$ of codimension $\ge 2$
such that 
$F_y\setminus (F_y\cap \Sing(f))$ is non-empty and connected for every $y\in Y\setminus \Xi_1$.
\end{itemize}
In particular,
 we see that $ f_*$  is surjective 
and   $\Im (\iota_*)=\Ker (f_*)$ holds by 
Nori's lemma~(Proposition~\ref{prop:nori}).
\end{remark}
Let  $\Sigma_1$, \dots, $\Sigma_N$ be 
the irreducible components of $\Sigma$ with codimension $1$ in $Y$.
There exists a proper Zariski closed subset $\Xi\subset \Sigma$ 
with the following properties.
We put 
$$
Y\spsh:=Y\setminus \Xi,
\quad
\Sigma_i\spsh :=\Sigma_i \setminus (\Sigma_i\cap \Xi)=\Sigma_i\cap Y\spsh,
\quad
\Sigma\spsh :=\Sigma\setminus \Xi =\Sigma\cap Y\spsh.
$$
\begin{itemize}
\item[($\Xi 0$)]  
The codimension of $\Xi$ in $Y$ is $\ge 2$.
\item[($\Xi 1$)] 
The  Zariski closed subsets
$\Xi_0\subset Y$ in the condition~\cond{C2} 
and
$\Xi_1\subset Y$ in the condition~\cond{C3}
are  contained in $\Xi$.
\item[($\Xi 2$)] 
Each  $\Sigma_i\spsh$ is a smooth hypersurface of $Y\spsh$,
and $\Sigma\spsh$ is a disjoint union of
$\Sigma_1\spsh, \dots, \Sigma_N\spsh$;
that is,
$\Xi$ contains all the irreducible components of $\Sigma$ with codimension $\ge 2$ in $Y$
and  the singular locus  of $\Sigma$. 
\item[($\Xi 3$)]
For each $y\in \Sigma_i\spsh$,
there exist an open neighborhood 
$U\subset Y\spsh$ of $y$ in $Y\spsh$ and an analytic isomorphism
$$
\phi : (U, U\cap \Sigma) \isomarrow \unitdisc^{m-1}\times (\unitdisc , 0),
\qquad \textrm{where $m=\dim Y$,}
$$
with the following properties.
Let $\psi: U\to\unitdisc^{m-1}$ be the composite of
$\phi : U\cong \unitdisc^{m-1}\times \unitdisc$ 
and the projection $\unitdisc^{m-1}\times \unitdisc\to \unitdisc^{m-1}$.
Then 
$$
\map{\Psi := \psi\circ f}{ f\inv (U)}{\unitdisc^{m-1}}
$$
is smooth, and the commutative diagram
$$
\renewcommand{\arraystretch}{1.2}
\begin{array}{ccccc}
f\inv (U) &&\maprightsp{f} && U \\
\hskip 15pt {}_\Psi\hskip -15pt  &\searrow && \swarrow &\hskip -7pt {}_\psi \hskip 7pt \\
&& \unitdisc^{m-1}&& 
\end{array}
$$
is a  trivial family of $\Cinf$-maps over $\unitdisc^{m-1}$ in the $\Cinf$-category.

\end{itemize}
Because of  the choice of $\Xi$,
for \emph{any} point $y\in \Sigma_i\spsh$,
there exists an open disc $\unitdisc\subset Y\spsh$
with the following properties:
\begin{itemize}
\item[\cond{$\Delta\spsh$1}]  $\unitdisc\cap \Sigma=\{y\}$, and 
 $\unitdisc$ intersects $\Sigma_i\spsh$ transversely at $y$,
 \item[\cond{$\Delta\spsh$2}]  $f\inv (\unitdisc)$ is a complex manifold, 
 \item[\cond{$\Delta\spsh$3}] $\shortmap{f|_{f\inv (\unitdisc)}}{f\inv (\unitdisc)}{ \unitdisc}$
is a one-dimensional family of complex analytic spaces 
that is locally trivial in the $\Cinf$-category over $\unitdisc\setminus\{y\}$, and 
 \item[\cond{$\Delta\spsh$4}] the central fiber $F_y:=f\inv (y)$ is an irreducible  hypersurface of $f\inv (\unitdisc)$,
 and $F_y\setminus (F_y\cap \Sing(f))$ is non-empty and connected.
\end{itemize}
Moreover
the diffeomorphism type of $\shortmap{f|_{f\inv (\unitdisc)}}{f\inv (\unitdisc)}{ \unitdisc}$
depends only on the index $i$ of  $\Sigma_i$.
\par
\medskip
We put
$$
X\spsh:=f\inv (Y\spsh),\;\; f\spsh:=f|_{X\spsh}:X\spsh\to Y\spsh,
\;\; \Theta\spsh_i:=\fspshinv (\Sigma_i\spsh) \;\;\textrm{and}\;\; \Theta\spsh:=\fspshinv (\Sigma\spsh).
$$
Then each $\Theta\spsh_i$ is an irreducible   hypersurface of $X\spsh$, 
and 
$\Theta\spsh$ is a disjoint union of $\Theta\spsh_1, \dots, \Theta\spsh_N$.
Note that  we have $X\spc=X\spsh\setminus \Theta\spsh$.
\begin{remark}\label{rem:C1}
By  the condition  \cond{C1}, 
the Zariski closed subset $f\inv (\Xi)$ of $X$ is also of codimension $\ge 2$,
and hence the inclusions induce isomorphisms
$\pione (X\spsh, \tlb)\cong \pione (X, \tlb)$ and $\pione (Y\spsh, b)\cong \pione (Y, b)$.
\end{remark}
We introduce notions of \emph{transversal discs},  \emph{leashed discs} and  \emph{lassos}.
\begin{definition}\label{def:defs1}
\setrmkakko
Let $H\subset M$ be a reduced hypersurface of a complex manifold $M$
of dimension $m$,
and let $H_1, \dots, H_l$
be the irreducible components of $H$.
We fix a base point $b_M\in M\setminus H$ .

\rmkakko
Let $N$ be a real $k$-dimensional $\Cinf$-manifold with $2\le k\le 2m$
(possibly with boundaries and corners),
and let $\phi: N\to M$ be a continuous map.
Let $p$ be a point of $N$ that is not in the corner of $N$.
If $k=2$, we further assume that $p\notin \bdr N$.
We say that 
\emph{ $\phi: N \to M$  intersects $H$  at  $p$ transversely}
 if  the following hold:
 \begin{itemize}
 \item[\cond{$\phi 1$}] $\phi(p)\in H\setminus  \Sing (H)$, and 
\item[\cond{$\phi 2$}]
 there exist local coordinates $(u_1, \dots, u_k)$
  of $N$ at $p$
 and local coordinates $(v_1, \dots, v_{2m})$
 of the $\Cinf$-manifold underlying $M$ at $\phi(p)$ such that
 \begin{itemize}
 \item[$\bullet$] $p=(0, \dots, 0)$, $\phi(p)=(0, \dots, 0)$, 
 \item[$\bullet$]  if $p\in \bdr N$, then 
 $N$ is given by $u_k\ge 0$ locally at $p$, 
 \item[$\bullet$]  $H$ is locally defined by $v_1=v_2=0$ in $M$, and 
 \item[$\bullet$] $\phi$ is given by
 $(u_1, \dots, u_k)\mapsto (v_1, \dots, v_{2m})=(u_1, \dots, u_k, 0, \dots, 0)$.
 \end{itemize}
\end{itemize}
We say that 
\emph{$\phi: N \to M$  intersects $H$ transversely}
if $\phi\inv (H)$ is disjoint from the corner of $N$
(when $k=2$, we  assume that $\phi\inv (H)\cap \bdr N=\emptyset$)
and  $\phi$ intersects $H$ transversely at every point of $\phi\inv (H)$.

If $\phi$ intersects $H$ transversely,
then $\phi\inv (H)$ is a real $(k-2)$-dimensional sub-manifold of $N$.
If $k>2$, then  the  boundary of $\phi\inv (H)$ is equal to $\phi\inv (H)\cap \bdr N$,
while if $k=2$, then $\phi\inv (H)$ is a finite set of points in the interior of $N$.

\rmkakko
A continuous map $\delta: \cunitdisc\to M$ is called a \emph{transversal disc  around $H_i$}
if $\delta\inv (H)=\{0\}$, $\delta(0)\in H_i$ and
$\delta$ intersects $H$ transversely at $0$.
In this case, the \emph{sign} of  $\delta$ 
is the  local intersection number ($+1$ or $-1$) of
$\delta$ with $H_i$ at $\delta (0)$.

\rmkakko
An \emph{isotopy} between  transversal discs $\delta$ and $\delta\sprime$  around  $H_i$
is a continuous map
$$
\map{h}{\cunitdisc\times I}{M}
$$
such that, for each $t\in I$,
the restriction $\delta_t:=\shortmap{h|_{\cunitdisc\times\{t\}}}{\cunitdisc}{M}$
of $h$ to $\cunitdisc\times\{t\}$ is a transversal disc  around $H_i$,
and such that $\delta_0=\delta$ and $\delta_1=\delta\sprime$ hold.

\rmkakko
A \emph{leashed disc} around $H_i$ with the base point $b_M$ is a pair $\rho=(\delta, \eta)$
of a  transversal  disc $\delta: \cunitdisc\to M$   around $H_i$
and a path $\eta: I\to M\setminus H$ from  $\delta (1)=\bdre\delta(0)=\bdre\delta(1)$
to $b_M$.
(Recall that $\bdre\delta$ is the loop given by $t\mapsto \delta(\exp(2\pi\sqrt{-1}t))$.
See Convention (3).)
The \emph{sign} of a leashed disc $\rho=(\delta, \eta)$  
is the sign of $\delta$.

\rmkakko
The  \emph{lasso} $\lambda(\rho)$ associated with a leashed disc $\rho=(\delta, \eta)$
is the loop $\eta\inv  \cdot (\bdre \delta)\cdot \eta$ in $M\setminus H$ with the base point $b_M$. 

\rmkakko
An \emph{isotopy} 
of  leashed discs
around $H_i$ with the base point $b_M$ 
is the pair of continuous maps
$$
\map{(h_{\cunitdisc}, h_I)}{ (\cunitdisc, I)\times I}{ (M, M\setminus H)}
$$
such that,
for each $t\in I$, 
the restriction of $(h_{\cunitdisc}, h_I)$ 
to $(\cunitdisc, I)\times\{t\}$
is a leashed disc around $H_i$ with the base point $b_M$.
\end{definition}
\begin{remark}\label{rem:lambda}
The isotopy class of a leashed disc $\rho$ is denoted by $[\rho]$.
If $[\rho]=[\rho\sprime]$,
then  $[\lambda(\rho)]=[\lambda(\rho\sprime)]$  holds in $\pione (M\setminus H, b_M)$.
\end{remark}
The following  is obvious:
\begin{proposition}\label{prop:rho}
\setrmkakko

\rmkakko
Any  two transversal discs  around   $H_i$
with the same sign  are isotopic.

\rmkakko
The homotopy classes of lassos associated with 
all the leashed discs around $H_i$
with a fixed  sign
form a conjugacy class
in $\pione (M\setminus H, b_M)$.

\rmkakko
The kernel of the homomorphism $\pione (M\setminus H, b_M)\to \pione (M,b_M)$
induced by the inclusion 
is generated by the homotopy classes of all lassos around $H_1, \dots, H_l$.
\end{proposition}
We apply these notions to the hypersurfaces
$$
\Sigma\spsh =\Sigma\spsh_1\cup\dots\cup \Sigma\spsh_N
\;\;\textrm{of $Y\spsh$},
\quand
\Theta\spsh =\Theta\spsh_1\cup\dots\cup \Theta\spsh_N
\;\;\textrm{of $X\spsh$}.
$$
\begin{definition}
\setrmkakko 
\rmkakko
A \emph{transversal lift} of a transversal  disc $\shortmap{\delta}{\cunitdisc}{Y\spsh}$
around $\Sigma_i\spsh$ is a lift
$\shortmap{\lift{\delta}}{\cunitdisc}{X\spsh}$
of $\delta$ with  
 $\lift{\delta} (0)\notin \Sing(f)$ 
 such that $\lift{\delta}$ intersects
 the irreducible hypersurface  $\Theta\spsh_i$ transversely at $0$.

\rmkakko
Let $\rho=(\delta, \eta)$ be a leashed disc around $\Sigma\spsh_i$ with the base point $b$.
A \emph{transversal  lift} of $\rho$ is a pair  $\lift{\rho}=(\lift{\delta},\lift{\eta})$ 
such that $\shortmap{\lift{\delta}}{\cunitdisc}{X\spsh}$
is a transversal  lift  of $\shortmap{\delta}{\cunitdisc}{Y\spsh}$
and $\shortmap{\lift{\eta}}{I}{X\spc}$ is a lift of 
$\shortmap{\eta}{I}{Y\spc}$ such that   
$\lift{\eta}(0)=\lift{\delta} (1)$
and $\lift{\eta}(1)=\tlb$.
\end{definition}
\begin{remark}
Any transversal  lift of a transversal disc (resp.~a leashed disc)
around $\Sigma_i\spsh$ 
is a transversal disc  (resp.~a leashed disc) around   $\Theta\spsh_i$.
Moreover the lifting does not change the sign.
\end{remark}
\begin{definition}\label{def:homotopylift}
\setrmkakko 
\rmkakko
Let $\delta_0$ and $\delta_1$ be two transversal discs on $Y\spsh$
around  $\Sigma_i\spsh$, and 
let $\shortmap{h}{\cunitdisc\times I }{Y\spsh}$
be an isotopy  of transversal discs from $\delta_0$ to $\delta_1$.
A \emph{lift} of the isotopy $h$ is a continuous map
$$
\map{\lift{h}}{\cunitdisc\times I }{X\spsh}
$$
such that, for each $t\in I$, the restriction $\lift{\delta}_t:=\lift{h}|_{\cunitdisc\times \{t\}}$ is
a  transversal lift of the transversal disc $\delta_t:=h|_{\cunitdisc\times \{t\}}$ on $Y\spsh$.
In particular,
we have $f\circ\lift{h}=h$ and 
$\lift{h}(\cunitdisc \times I)\cap \Sing (f)=\emptyset$.
Moreover $\lift{h}$ is an isotopy of transversal discs
around $\Theta_i\spsh$ from $\lift{\delta}_0$ to $\lift{\delta}_1$.
By abuse of notation,
we sometimes say that the isotopy $\lift{\delta}_t$ is the transversal lift
of the isotopy ${\delta}_t$,
understanding that $t$ is the homotopy parameter.

\rmkakko
Let $\rho_0$ and $\rho_1$ be two leashed discs on $Y\spsh$
around to $\Sigma_i\spsh$, and 
let $\shortmap{(h_{\cunitdisc}, h_I)}{(\cunitdisc, I)\times I }{(Y\spsh, Y\spc)}$
be an isotopy  of leashed discs from $\rho_0$ to $\rho_1$.
A \emph{lift} of the isotopy $(h_{\cunitdisc}, h_I)$ is a pair of  continuous maps
$$
\map{(\lift{h}_{\cunitdisc}, \lift{h_I})}{(\cunitdisc, I)\times I  }{(X\spsh, X\spc)}
$$
such that, for each $t\in I$, the restriction 
$\lift{\rho}_t:=(\lift{h}_{\cunitdisc}, \lift{h}_I)|_{(\cunitdisc, I)\times \{t\}}$ is
a transversal lift of the leashed  disc $\rho_t:=(h_{\cunitdisc}, h_I)|_{(\cunitdisc, I)\times \{t\}}$ on $Y\spsh$.
\end{definition}
The following are obvious
from the condition \cond{$\Delta\spsh$4}:
\begin{proposition}\label{prop:homotopyliftexistence}
Every transversal disc 
around  $\Sigma_i\spsh$ 
has a transversal lift 
on $X\spsh$.
Moreover, every isotopy $\delta_t$ of transversal discs  around $\Sigma_i\spsh$
from $\delta_0$ to $\delta_1$ 
lifts to an isotopy   $\lift{\delta}_t$
from  a given transversal lift $\lift{\delta}_0$ 
of $\delta_0$ to a given  transversal lift $\lift{\delta}_1$
of $\delta_1$.
\end{proposition}
\begin{remark}\label{rem:homotopyliftexistence2}
Every leashed disc 
on $Y\spsh$ around $\Sigma_i\spsh$ has a transversal lift 
on $X\spsh$.
Moreover, 
every isotopy $\rho_t$ of leashed  discs  on $Y\spsh$
has a lift  $\lift{\rho}_t$
on $X\spsh$ from  a given transversal lift $\lift{\rho}_0$ of $\rho_0$,
but the ending lift $\lift{\rho}_1$ cannot be arbitrarily given.
\end{remark}
%
%
\begin{definition}
Let  $\rho$ be a leashed disc on $Y\spsh$ around $\Sigma_i\spsh$,
and let $\lift{\rho}$ be a transversal lift of $\rho$.
Then we have the lasso $\lambda(\lift{\rho})$,
which is a loop in $X\spc$
with the base point $\tlb$.
Recall that $\mu$ is the lifted monodromy.
We put
$$
N(\lift{\rho}):=
\gen{\;\shortset{g\inv g^{\mu([\lambda(\lift{\rho})])}}{g\in \pione (\Fb, \tlb)}\;}_{\pione (\Fb, \tlb)}.
$$
\end{definition}
\begin{proposition_definition}\label{prop:N}
Let $\rho\sprime$ be a leashed disc on $Y\spsh$ isotopic  to $\rho$,
and let $\lift{\rho}\sprime$ be a transversal lift of $\rho\sprime$.
Then we have 
$$
N(\lift{\rho})\;=\;N(\lift{\rho}\sprime).
$$
Therefore,
for an  isotopy class $[\rho]$ of  leashed discs on $Y\spsh$,
 we can  define a normal subgroup
$N^{[\rho]}$ of $\pione (\Fb, \tlb)$ by choosing a transversal lift
$\lift{\rho}$ of a representative $\rho$ of $[\rho]$, and putting
$$
N^{[\rho]}:=N(\lift{\rho}).
$$
\end{proposition_definition}
\begin{proof}
By Remarks~\ref{rem:lambda} and~\ref{rem:homotopyliftexistence2},
the isotopy from $\rho$ to $\rho\sprime$ lifts to
an isotopy from $\lift{\rho}$ to some lift $\lift{\rho}\sprime_1$ of $\rho\sprime$,
and we have $[\lambda(\lift{\rho})]=[\lambda(\lift{\rho}_1\sprime)]$ in $\pione (X\spc, \tlb)$.
(However $[\lambda(\lift{\rho}\sprime_1)]$ and $[\lambda(\lift{\rho}\sprime)]$ 
may be distinct in general.)
Therefore it is enough to show that
$N(\lift{\rho}^{(1)})=N(\lift{\rho}^{(2)})$ holds
for any two transversal
lifts $\lift{\rho}^{(1)}=(\lift{\delta}^{(1)}, \lift{\eta}^{(1)})$ and 
$\lift{\rho}^{(2)}=(\lift{\delta}^{(2)}, \lift{\eta}^{(2)})$ 
of a single leashed disc $\rho=(\delta, \eta)$ on $Y\spsh$.
We can assume that the transversal disc $\delta: \cunitdisc\to Y\spsh$
around  $\Sigma\spsh_i$
is an embedding of a complex manifold.
We denote by $\cunitdisc_\rho$ the image of $\delta$,
and by $\unitdisc_\rho$ the interior of $\cunitdisc_\rho$.
We can further assume that 
$\cunitdisc_\rho$ is sufficiently small,
and that 
$$
E_\rho:=f\inv(\unitdisc_\rho)
$$
is a smooth complex manifold
by the condition~\cond{$\Delta\spsh$2}.
We then put
$$
\ol{E}_\rho\:=f\inv(\cunitdisc_\rho),
\quad
\ol{E}\sptimes_\rho\:=f\inv(\cunitdisc_\rho\sptimes),
$$
where $\cunitdisc_\rho\sptimes:=\cunitdisc_\rho\setminus \{\delta(0)\}=\cunitdisc_\rho\cap Y\spc$.
We also put $q:=\delta(1)=\eta(0)\in \bdr \cunitdisc_\rho$ and 
$$
\lift{q}^{(1)}:=\lift{\delta}^{(1)}(1)=\lift{\eta}^{(1)}(0)\in F_q,
\quad
\lift{q}^{(2)}:=\lift{\delta}^{(2)}(1)=\lift{\eta}^{(2)}(0)\in F_q.
$$
Since $f$ is locally trivial over $\eta(I)\subset Y\spc$
and $\sqcap=(\bdr I \times I)\cup (I\times\{1\})$ is a strong deformation retract of $I\times I$,
there exists a continuous map
$\shortmap{\Omega}{I\times I}{\Xc}$
such that the following hold for any $s, t\in I$:
$$
f(\Omega (s, t))=\eta (t),
\quad
\Omega(s, 1)=\tlb,
\quad 
\Omega(0, t)=\lift{\eta}^{(1)} (t),
\quad
\Omega(1, t)=\lift{\eta}^{(2)} (t).
$$
(See Figure~\ref{figOmega}.)
Then,
for each $t\in I$,  the map $s\mapsto \Omega(s, t)$ is a path in $F_{\eta (t)}$ from 
$\lift{\eta}^{(1)} (t)$ to $\lift{\eta}^{(2)} (t)$.
We denote by $\omega: I\to F_q$ the path
in $F_q$ from $\lift{q}^{(1)}$ to $\lift{q}^{(2)}$
defined by  $\omega(s):=\Omega(s, 0)$.
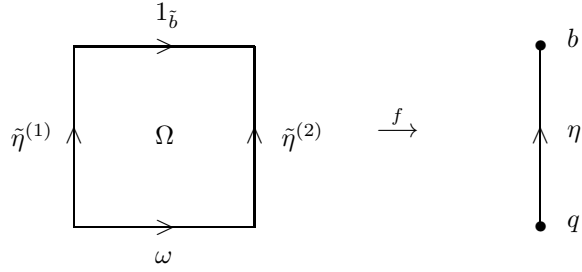
\begin{figure}
\begin{center}
$
\setlength{\unitlength}{1.2mm}
\vcenter{\hbox{
\begin{picture}(30,30)(-5,-5)
\put (0,0) {\line (1,0){20}}\put (10,0){\rightya}
\put (9, -4) {$\omega$}
\put (20,0) {\line (0,1){20}}  \put (20,10){\upya}
\put (23, 9) {${\lift{\eta}}^{(2)}$}
\put (20,20) {\line (-1,0){20}}  \put (10,20){\rightya}
\put (9, 23) {$1_{\tlb}$}
\put (0,20) {\line (0,-1){20}} \put (0,10){\upya}
\put (-7, 9) {${\lift{\eta}}^{(1)}$}
\put (9, 9.2) {$\Omega$}
\end{picture}
}}
\qquad
\vcenter{\hbox{$\maprightsp{f}$}}
\qquad
\vcenter{\hbox{
\begin{picture}(10,30)(-5,-5)
\put (0,0) {\line (0, 1){20}} \put (0,10){\upya}
\put (3,10) {$\eta$}
\put (3,0) {$q$}
\put (0,0) {\kuromaru}
\put (3,20) {$b$}
\put (0,20) {\kuromaru}
\end{picture}
}}
$
\end{center}
\caption{The map $\Omega$}\label{figOmega}
\end{figure}
Then we have the following commutative diagram:
$$
\begin{array}{ccccc}
\pione (\Fb, \tlb) & \mapleftspsb{\sim}{[\lift{\eta}^{(1)}]_*} 
& \pione (F_q, \lift{q}^{(1)}) & \maprightsp{i_{q *}} & \pione(\ol{E}_\rho, \lift{q}^{(1)}) \\
\parallel& & \mapdownleftright{\hskip -10pt [\omega]_*}{\wr}  & &  \mapdownleftright{\hskip -10pt [\omega]_*}{\wr} \\
\pione (\Fb, \tlb) & \mapleftspsb{\sim}{[\lift{\eta}^{(2)}]_*} 
& \pione (F_q, \lift{q}^{(2)}) & \maprightsp{i_{q *}} 
& \pione(\ol{E}_\rho, \lift{q}^{(2)}),\\
&&
\end{array}
$$
where $i_{q}: F_q\inj \ol{E}_\rho$ is the inclusion.
Hence,
in order to prove $N(\lift{\rho}^{(1)})=N(\lift{\rho}^{(2)})$,  it is enough to show the following equality:
$$
[\lift{\eta}^{(1)}]_*\inv (N(\lift{\rho}^{(1)}))=\Ker(i_{q *} : \pione (F_q, \lift{q}^{(1)})  \to \pione(\ol{E}_\rho, \lift{q}^{(1)})).
$$
Since
$\shortmap{f|_{\ol{E}_\rho}}{\ol{E}_\rho}{\cunitdisc_\rho}$
is locally trivial over $\cunitdisc_\rho\sptimes$
with the general fiber being connected by \cond{C0},
and since there exists a cross-section
$$
\map{{}^s\lift{\delta}^{(1)}}{\cunitdisc_\rho}{\ol{E}_\rho}
$$
of $f|_{\ol{E}_\rho}$
given by the transversal lift  $\lift{\delta}^{(1)}$ of $\delta$, 
we have an exact sequence 
$$
1\;\maprightsp{}\; 
\pione (F_q, \lift{q}^{(1)})\;\maprightsp{i_{q*}}\;
\pione (\ol{E}_\rho\sptimes, \lift{q}^{(1)} )\;\maprightsp{(f|_{\ol{E}\sptimes_\rho})_*}\;
\pione(\cunitdisc_{\rho}\sptimes, q)\;\maprightsp{}\;
1
$$
with the cross-section 
$$
\map{s}{\pione(\cunitdisc_{\rho}\sptimes, q)}{\pione (\ol{E}_\rho\sptimes, \lift{q}^{(1)} )}
$$ 
of $(f|_{\ol{E}\sptimes_\rho})_*$ that maps the positive generator $[\bdre \delta]$ of $\pione(\cunitdisc_{\rho}\sptimes, q) \cong \Z$
to $[\bdre \lift{\delta}^{(1)}]\in \pione (\ol{E}_\rho\sptimes, \lift{q}^{(1)} )$.
By the cross-section ${}^s\lift{\delta}^{(1)}$ of $f|_{\ol{E}_\rho}$
 over $\cunitdisc_\rho$,
we have the classical monodromy action of $\pione(\cunitdisc_{\rho}\sptimes, q)$ on $\pione (F_q, \lift{q}^{(1)})$.
By the definition,
the  action of $[\bdre \delta]\in \pione(\cunitdisc_{\rho}\sptimes, q)$ is equal to 
$$
g\;\;\mapsto\;\; g^{\mu([\bdre \lift{\delta}^{(1)}])}=
[\bdre \lift{\delta}^{(1)}]\inv \cdot   g  \cdot[\bdre \lift{\delta}^{(1)}]
\quad\textrm{for}\quad g\in \pione (F_q, \lift{q}),
$$
where the product is taken in $\pione (\ol{E}_\rho\sptimes, \lift{q}^{(1)} )$
and $\pione (F_q, \lift{q}^{(1)})$ is regarded as a normal subgroup of 
$\pione (\ol{E}_\rho\sptimes, \lift{q}^{(1)} )$ by $i_{q*}$.
Hence, by Lemma~\ref{lem:GGG},
$\pione (\ol{E}_\rho\sptimes, \lift{q}^{(1)} )$ is isomorphic to the semi-direct product
$\pione (F_q, \lift{q}^{(1)})\semidirectproduct \pione(\cunitdisc_{\rho}\sptimes, q)$
constructed by the monodromy  action.
On the other hand,
by the condition~\cond{$\Delta\spsh$4},
the central fiber  $F_{\delta(0)}$ of $\ol{E}_\rho\to \cunitdisc_\rho$ is 
an irreducible hypersurface of $\ol{E}_\rho$, and hence 
the kernel of
$$
\map{j_*}{\pione (\ol{E}_\rho\sptimes, \lift{q}^{(1)} )}{\pione (\ol{E}_\rho, \lift{q}^{(1)} )}
$$
induced by the inclusion 
$j: \ol{E}_\rho\sptimes\inj \ol{E}_\rho$ is generated by the conjugacy class of lassos
around  $F_{\delta(0)}$.
(See Proposition~\ref{prop:rho}.)
Since $\bdre \lift{\delta}^{(1)}=\lambda(\lift{\delta}^{(1)})$ is a lasso around $F_{\delta(0)}$, 
 the kernel of $j_*$
is equal to the normal subgroup
 $\ngen{\{[\bdre \lift{\delta}^{(1)}]\}}=\ngen{\Im (s)}$.
By Lemmas~\ref{lem:S} and~\ref{lem:sGamma},
the kernel of 
the composite
$$
\pione (F_q, \lift{q}^{(1)})
\;\maprightsp{i_{q*}}\;
\pione (\ol{E}_\rho\sptimes, \lift{q}^{(1)} )
\;\maprightsp{j_*}\;
\pione (\ol{E}_\rho, \lift{q}^{(1)} )=\pione (\ol{E}_\rho\sptimes, \lift{q}^{(1)} )/\ngen{\Im (s)}
$$
is equal to
$$
N\sprime:=\gen{\set{g\inv g^{\mu([\bdre \lift{\delta}^{(1)}])}}{g\in \pione (F_q, \lift{q}^{(1)})}}.
$$
Since 
$[\lift{\eta}^{(1)}]_* (g^{\mu([\bdre \lift{\delta}^{(1)}])})=([\lift{\eta}^{(1)}]_* (g))^{\mu([\lambda(\lift{\rho}^{(1)})])}$
for any $g\in \pione(F_q, \lift{q}^{(1)})$,
we see that
$[\lift{\eta}^{(1)}]_*$ induces an isomorphism $N\sprime \isom N(\lift{\rho}^{(1)})$.
\end{proof}
\begin{proposition}\label{prop:Nalpha}
Let $\shortmap{\lift{\gamma}}{(I, \bdr I)}{(X\spc, \tlb)}$
be a loop, and we put $\gamma:=f\circ \lift{\gamma}$.
Then,
for any leashed disc 
$\rho=(\delta, \eta)$ on $Y\spsh$ around $\Sigma_i\spsh$, we have
$$
({N^{[\rho]}})\sp{\mu([\lift{\gamma}])} =N^{[(\delta, \eta\gamma)]}.
$$
\end{proposition}
\begin{proof}
Let $g$ be an element of $\pione(\Fb, \tlb)$,
and let $h$ denote $g^{\mu([\lift{\gamma}])}$. 
Then, for a  transversal lift $\lift{\rho}=(\lift{\delta}, \lift{\eta})$ of $\rho$, 
we have
$$
(g\inv g^{\mu([\lambda(\lift{\rho})])})^{\mu([\lift{\gamma}])}
=h\inv h ^{\mu([\lift{\gamma}]\inv [\lambda(\lift{\rho})] [\lift{\gamma}])}.
$$
Since  
$\lift{\gamma}\inv  \lambda(\lift{\rho}) \lift{\gamma}=
\lift{\gamma}\inv  \lift{\eta}\inv  \cdot \bdre{\lift{\delta}} \cdot \lift{\eta}  \lift{\gamma}$
is a lasso
associated with the transversal lift $(\lift{\delta},\lift{\eta}\lift{\gamma})$
of the leashed disc $(\delta, \eta\gamma)$,
we obtain the proof.
\end{proof}
\begin{corollary}\label{cor:foranyrho}
If $N^{[\rho]}=1$ holds
for one leashed disc $\rho$ around $\Sigma_i\spsh$,
then we have 
$N^{[\rho]}=1$
for any leashed disc $\rho$ around $\Sigma_i\spsh$.
\end{corollary}
We can now state  the main result of this section.
\begin{theorem}\label{thm:ZvK}
Suppose that the conditions \cond{C1},~\cond{C2} and the following condition~\cond{Z}  are satisfied:
\begin{itemize}
\item[\cond{Z}]
There exists a continuous cross-section $s_Z: Z\to f\inv (Z)$ 
of $f$ over a subspace $Z\subset Y$ satisfying
$b\in Z$, $s_Z(b)=\tlb$, $s_Z(Z)\cap \Sing (f)=\emptyset$ and 
such that
 the inclusion $Z\inj Y$ induces 
a surjection  $\pi_2 (Z, b)\surj \pi_2(Y, b)$.
\end{itemize}
Let $\LLL$ be the set of 
isotopy classes of all leashed discs on $Y\spsh$ 
around  $\Sigma\spsh_1, \dots, \Sigma\spsh_N$.
Then $\Ker (\iota_*)$
is equal to
$$
\NNN:=\gen{\;\;\textstyle{\bigcup}_{[\rho]\in \LLL} N^{[\rho]}\;\;}_{\pione (\Fb, \tlb)}.
$$
\end{theorem}
\begin{remark}\label{rem:pitwo}
If $\pi_2(Y)=0$, then the condition~\cond{Z} is always satisfied,
because we can put $Z=\{b\}$ and $s_Z(b)=\tlb$.
\end{remark}
For the proof,
we define the notion of \emph{free loop pairs of monodromy relation type}.
Let $\sphere^1$ denote the oriented circle.
\begin{definition}
Let $T$ be a topological space.
A \emph{free loop}
on $T$ is a continuous map
$\shortmap{\varphi}{\sphere^1}{T}$.
A \emph{homotopy} from a free loop $\varphi$ to a free loop  $\varphi\sprime$
is a  continuous map
$\shortmap{\Phi}{\sphere^1\times I}{T}$
such that $\Phi|_{{\sphere^1}\times\{0\}} =\varphi$
and $\Phi|_{{\sphere^1}\times\{1\}} =\varphi\sprime$.
The homotopy class of a free loop $\varphi$ is denoted by $[\varphi]_\FL$.
\end{definition}
Suppose that $T$ is path-connected, and let $b_T$ be  a  base point of $T$.
Then the natural map $[\alpha]\mapsto [\alpha]_\FL$ induces a bijection from the set of
conjugacy classes of $\pione (T, b_T)$ to the set of homotopy classes 
of free loops on $T$.
\par
\medskip
Let $D$ be a topological space homeomorphic to $\cunitdisc$,
let $b_D$ be a point  of $D$, 
and let $\bdr D$ be the boundary of $D$
with an orientation.
\begin{definition}\label{def:frp}
A \emph{free loop pair} is a pair 
$$
\map{(\psi, \Lift{\psi|_{\bdr D}})}{(D, \bdr D)}{(Y\spc, X\spc)}
$$
 of 
a continuous map $\shortmap{\psi}{ D}{ Y\spc}$
and a  lift $\shortmap{\Lift{\psi|_{\bdr D}}}{\bdr D}{X\spc}$ of the restriction $\shortmap{\psi|_{\bdr D}}{\bdr D }{Y\spc}$
of $\psi$ to $\bdr D$.
\end{definition}
\begin{remark}
The notation $\shortmap{(\psi, \Lift{\psi|_{\bdr D}})}{(D, \bdr D)}{(Y\spc, X\spc)}$ 
for a pair of maps is different from 
the usual notation in topology,  
because $X\spc$ is \emph{not} a subspace of $Y\spc$.
The same warning is also applied to Definition~\ref{def:hfrp}.
\end{remark}
Let $\shortmap{(\psi, \Lift{\psi|_{\bdr D}})}{(D, \bdr D)}{(Y\spc, X\spc)}$
be a free loop pair.
Consider the pull-back 
$$
\map{\psi^*(f\spc)}{\psi^*(X\spc):=X\spc\times_{Y\spc} D}{ D }
$$
of the locally trivial map $\shortmap{f\spc}{X\spc}{Y\spc}$ by $\psi$.
Since $ D $ is contractible,
 we have a contraction $c: \psi^*(X\spc)\to F_{\psi(b_D)}$,
 which is the homotopy inverse of the inclusion $F_{\psi(b_D)}\inj \psi^*(X\spc)$.
Then the cross-section 
$$
\map{{}\sp{s}\Lift{\psi|_{\bdr D}}}{\bdr D }{\psi^*(X\spc)}
$$
of $\psi^*(f\spc) $ over $\bdr D $
obtained from $\shortmap{\Lift{\psi|_{\bdr D}}}{\bdr D}{X\spc}$ defines a homotopy class
$[\Lift{\psi|_{\bdr D}}]_\FL$ of free loops on $F_{\psi(b_D)}$ via the contraction $c$,
and hence a conjugacy class $\Conj (\psi, \Lift{\psi|_{\bdr D}})$ of 
$\pione (F_{\psi(b_D)}, \tlb\sprime)$,
where $\tlb\sprime\in F_{\psi(b_D)}$ is an arbitrary base point.
Remark that $\Conj (\psi, \Lift{\psi|_{\bdr D}})$ 
does not depend on the choice of the contraction $c$.
\begin{definition}\label{def:frpmrt}
We choose a path $\lift{\alpha}$ in $X\spc$ from $\tlb\in F_b$ to $\tlb\sprime\in F_{\psi(b_D)}$.
We say that the free loop pair 
$$
\map{(\psi, \Lift{\psi|_{\bdr D}})}{(D, \bdr D)}{(Y\spc, X\spc)}
$$
 is  \emph{of monodromy relation type around $\Sigma_i\spsh$}
if the pull-back of the conjugacy class 
$\Conj (\psi, \Lift{\psi|_{\bdr D}})\subset \pione (F_{\psi(b_D)}, \tlb\sprime)$ by the isomorphism
$[\lift{\alpha}]_* : \pione(F_b, \tlb)\isom \pione (F_{\psi(b_D)}, \tlb\sprime)$
is contained in $N^{[\rho]}$ for some leashed disc $\rho$ on $Y\spsh$ around $\Sigma_i\spsh$.
\end{definition}
\begin{remark}
It is obvious that 
this definition does not depend on the choice of the orientation of $\bdr D$.
It also follows from 
Proposition~\ref{prop:Nalpha} that 
this definition does not depend on the choice of the path $\lift{\alpha}$
connecting $\tlb\in F_b$ and $\tlb\sprime\in F_{\psi(b_D)}$.
\end{remark}
\begin{definition}\label{def:hfrp}
A \emph{homotopy of free loop pairs}
is a pair of continuous maps
$$
\map{(h, \Lift{h|_{\bdr D}})}{(D, \bdr D)\times I}{(Y\spc, X\spc)}
$$
such that, for each $u\in I$,
the restriction of $(h, \Lift{h|_{\bdr D}})$ to $(D, \bdr D)\times \{u\}$ is a free loop pair.
\end{definition}
\begin{remark}\label{rem:homotopymonrel}
Suppose that two free loop pairs are homotopic.
If one is of monodromy relation type around $\Sigma_i\spsh$,
then so is the other.
\end{remark}
\begin{remark}\label{rem:homotopyD}
Let 
$\shortmap{\psi_u}{D}{Y\spc}$ 
 be a homotopy of continuous maps from $\psi_0$ to $\psi_1$ parametrized by $u\in I$.
Since $f\spc$ is locally trivial,
the homotopy 
$\shortmap{\psi_u|_{\bdr D}}{\bdr D}{Y\spc}$
lifts to a homotopy
$\shortmap{\Lift{\psi_u|_{\bdr D}}}{\bdr D}{X\spc}$
that starts from any  given lift $\Lift{\psi_0|_{\bdr D}}$
of  ${\psi_0|_{\bdr D}}$ 
and hence we obtain a homotopy $(\psi_u, \Lift{\psi_u|_{\bdr D}})$  of free loop pairs
starting from a given $(\psi_0, \Lift{\psi_0|_{\bdr D}})$.
(The ending lift $\Lift{\psi_1|_{\bdr D}}$
cannot be arbitrarily given.)
\end{remark}
\begin{proposition}\label{prop:TD}
Let $\delta_0$ and $\delta_1$ be  two transversal discs on $Y\spsh$ around $\Sigma_i\spsh$,
and let 
$\shortmap{h}{\cunitdisc\times I}{Y\spsh}$
be an isotopy  of transversal discs from $\delta_0=h|_{\cunitdisc\times\{0\}}$ to 
$\delta_1=h|_{\cunitdisc\times\{1\}}$.
Let $D$ be a closed subset of $\bdr\cunitdisc\times (I\setminus\bdr I)$
homeomorphic to $\cunitdisc$, and  put
$$
T:=\bdr(\cunitdisc \times I)\setminus (D\setminus \bdr D),
$$
so that $\bdr T=\bdr D$.
Suppose that we are given a lift
$$
\map{\Lift{h|_T}}{T}{X\spsh}
$$
of $\shortmap{h|_T}{T}{Y\spsh}$
such that the restrictions 
$$
\shortmap{\lift{\delta}_0:=\Lift{h|_T}|_{\cunitdisc\times\{0\}}}{\cunitdisc}{X\spsh}
\quand
\shortmap{\lift{\delta}_1:=\Lift{h|_T}|_{\cunitdisc\times\{1\}}}{\cunitdisc}{X\spsh}
$$
are transversal lifts of $\delta_0$ and $\delta_1$,  respectively.
Then the free loop pair
$$
\map{(h|_D, \Lift{h|_T}|_{\bdr D})}{(D, \bdr D)}{(Y\spc, X\spc)}
$$
is  of monodromy relation type around $\Sigma_i\spsh$.
\end{proposition}
\begin{remark}
In Figure~\ref{figLiftH},
the closed subset  $D$ is the region surrounded by the dashed curve on the right tube $\cunitdisc\times I$.
\end{remark}
%
%
\renewcommand{\PStextplot}[3]{\rlap{\hskip -300.000000 pt \hbox{\hskip #1pt  \raise #2pt \hbox{#3}}}}%
\begin{figure}
 \begin{center}
 \includegraphics{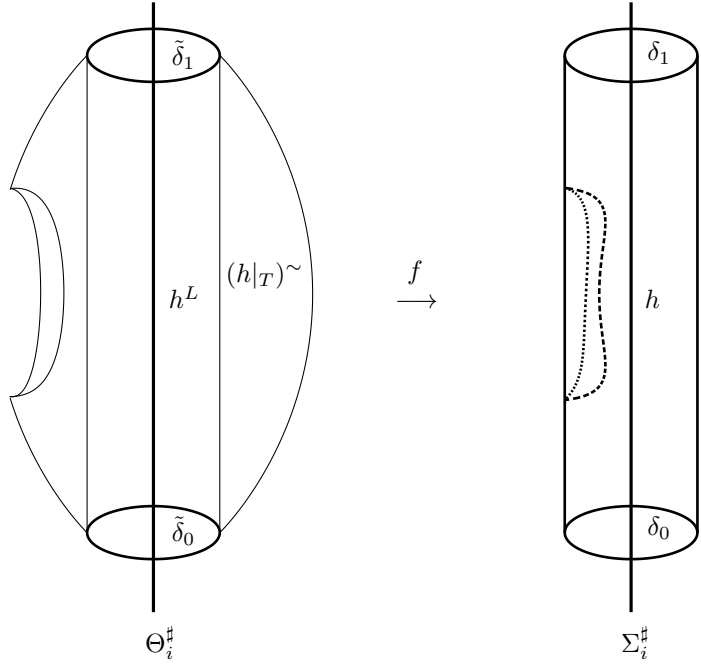}%
\PStextplot{251.000000}{209.000000}{${\delta}_1$}%
\PStextplot{251.000000}{29.000000}{${\delta}_0$}%
\PStextplot{71.000000}{207.000000}{$\lift{\delta}_1$}%
\PStextplot{71.000000}{27.000000}{$\lift{\delta}_0$}%
\PStextplot{241.000000}{-15.000000}{$\Sigma_i\spsh$}%
\PStextplot{61.000000}{-15.000000}{$\Theta_i\spsh$}%
\PStextplot{70.000000}{115.000000}{$h\sp{L}$}%
\PStextplot{91.500000}{125.000000}{$\Lift{h|_T}$}%
\PStextplot{150.000000}{115.000000}{$\maprightsp{}$}%
\PStextplot{160.000000}{126.000000}{$f$}%
\PStextplot{250.000000}{115.000000}{$h$}%
\end{center}
\caption{$(h|_T)^\sim$ and $h^L$}%
\label{figLiftH}%
\end{figure}%
\begin{proof}[Proof of Proposition~\ref{prop:TD}]
First note that,
since $h$ is an isotopy of transversal discs,
the image of $\bdr\cunitdisc \times I$ by $h$ is contained in $Y\spc$,
and hence we have $h|_D(D)\subset Y\spc$.
\par
\medskip
By Remarks~\ref{rem:homotopymonrel} and~\ref{rem:homotopyD}, we can assume that
$D\cap (\{1\} \times I)=\emptyset$
by moving $D$  by a homeomorphism of  $\bdr\cunitdisc\times I$
homotopic to the identity.
We consider the continuous map 
$$
\map{\tau }{I^2 }{ \bdr\cunitdisc \times I}
$$
given by $\tau(s, t):=(\exp(2\pi\sqrt{-1}s), t)$.
Then we have $D\subset \tau (I^2\setminus \bdr I^2)$ and 
 $\tau (\bdr I^2)\subset T$.
Under a suitable homeomorphism between $D$ and $I^2$,
the inclusion $D\inj \bdr\cunitdisc \times I$ is homotopic to
$\tau$.
We put
$$
H_0:=h\circ \tau\;:\; I^2\to Y\spc
$$
and define a lift $\Lift{H_0 |_{\bdr I^2}}$ of $H_0 |_{\bdr I^2}$ by
$$
\Lift{H_0 |_{\bdr I^2}}:=\Lift{h|_T} \circ(\tau|_{\bdr I^2}) \;: \; \bdr I^2\to X\spc.
$$
By Remarks~\ref{rem:homotopymonrel} and~\ref{rem:homotopyD} again, it is enough to prove that 
the free loop pair
$$
\map{(H_0, \Lift{H_0 |_{\bdr I^2}})}{(I^2, \bdr I^2)}{(Y\spc, X\spc)}
$$
is  of monodromy relation type around $\Sigma_i\spsh$.
For simplicity, we put
\begin{eqnarray*}
&&q:=\delta_0 (1)=h(1,0)=H_0 (0,0)=H_0(1,0), \quand \\
&&\lift{q}:=\lift{\delta}_0(1)=\Lift{h|_T}(1,0)=\Lift{H_0|_{\bdr I^2}} (0,0)=\Lift{H_0|_{\bdr I^2}} (1,0)\in F_q.
\end{eqnarray*}
By Proposition~\ref{prop:homotopyliftexistence}, 
we have an isotopy
$$
\map{h\sp{L}}{\cunitdisc\times I}{X\spsh}
$$
of transversal discs around $\Theta\spsh_i$ 
from $\lift{\delta}_0=\Lift{h|_T}|_{\cunitdisc\times\{0\}}$
to $\lift{\delta}_1=\Lift{h|_T}|_{\cunitdisc\times\{1\}}$
that is a lift of 
the isotopy  $\shortmap{h}{\cunitdisc\times I}{Y\spsh}$;
$$
f\circ h\sp{L}=h.
$$
In Figure~\ref{figLiftH}, 
the left tube is $h^L$,
while the barrel with a hole is $\Lift{h|_T}$.
We  put
$$
\map{\delta_t:=h|_{\cunitdisc\times\{t\}}}{\cunitdisc}{Y\spsh}
\quand
\map{\lift{\delta}_t:=h\sp{L}|_{\cunitdisc\times\{t\}}}{\cunitdisc}{X\spsh}.
$$
Then $\lift{\delta}_t$ 
is a transversal lift of 
$\delta_t$.
Next we  put
$$
\map{k_0:=h|_{\{1\}\times I}}{I}{Y\spc},
$$
which is a path on $Y\spc$ from $q=\delta_0 (1)$ to $\delta_1 (1)$,
and 
$$
\lift{k}_0:=\Lift{h|_T}|_{\{1\}\times I}=\Lift{H_0 |_{\bdr I^2}}|_{\{0\}\times I}=\Lift{H_0 |_{\bdr I^2}}|_{\{1\}\times I},
$$
which is 
 a lift of $k_0$ from $\lift{q}=\lift{\delta}_0 (1)$ to $\lift{\delta}_1 (1)$. 
Note that,
with the base point $(0,0)$ and the
orientation of $\bdr I^2$ given in Figure~\ref{figorientation},
the map  $\shortmap{\Lift{H_0 |_{\bdr I^2}}}{\bdr I^2}{X\spc}$
is equal to 
$$
\lift{k}_0\cdot \bdre\lift{\delta}_1\cdot \lift{k}_0\inv \cdot \bdre\lift{\delta}_0 \inv
$$
as a loop with the base point $\lift{q}=\Lift{H_0 |_{\bdr I^2}}(0,0)\in F_{q}$.
\begin{figure}
\begin{center}
\setlength{\unitlength}{1.2mm}
\begin{picture}(30,30)(-5,-2)
\put (0,0) {\line (1,0){20}}\put (10,0){\leftya}
\put  (20, 0) {\;\;$(1,0)$}
\put (20,0) {\line (0,1){20}}  \put (20,10){\downya}
\put  (20, 20) {\;\;$(1,1)$}
\put (20,20) {\line (-1,0){20}}  \put (10,20){\rightya}
\put  (0, 20) {\llap{$(0,1)$\;\;}}
\put (0,20) {\line (0,-1){20}} \put (0,10){\upya}
\put  (0, 0) {\llap{the base point $(0,0)$\;\;}}
\put  (0, 0) {\kuromaru}
\put (9, 9.2) {$I^2$}
\end{picture}
\end{center}
\caption{An  orientation of  $\bdr I^2$}\label{figorientation}
\end{figure}
We define a homotopy 
$$
\map{H_u}{I^2}{Y\spc}\qquad (u\in I)
$$
with $u$ being the homotopy parameter 
by $H_u(s, t):=H_0(s, (1-u)t)$,
and will construct a homotopy
$\shortmap{\Lift{H_u|_{\bdr I^2}}}{\bdr I^2}{X\spc}$
that covers the homotopy ${H_u|_{\bdr I^2}}$ and starts from $\Lift{H_0|_{\bdr I^2}}$ above.
We define
$$
\map{K}{I\times I}{Y\spc}
$$
by $K(t, u):=k_0 ((1-u)t)$,
and put
$k_u:=K|_{I\times\{u\}}$
for $u\in I$.
Then  $k_u$ gives  a homotopy  with parameter $u\in I$
from $k_0$ to the constant map $k_1=1_{q}$. 
We then define
a lift
$\shortmap{\Lift{K|_\sqcup}}{\sqcup}{X\spc}$
of $\shortmap{K|_\sqcup}{\sqcup}{Y\spc}$,
where $\sqcup:=(\bdr I\times I)\cup (I\times\{0\})$, 
by the following:
$$
\Lift{K|_\sqcup}(t, u):=
\begin{cases}
\lift{q} & \textrm{if $t=0$,}\\
\lift{k}_0 (t) & \textrm{if $u=0$,}\\
\lift{\delta}_{1-u} (1)=h^L(1, 1-u)  & \textrm{if $t=1$.}\\
\end{cases}
$$
Since $f\spc$ is locally trivial,
the lift $\Lift{K|_\sqcup}$ extends to a lift
$\shortmap{\lift{K}}{I\times I}{X\spc}$
of $K$.
(See Figure~\ref{figK}.)
\begin{figure}
\begin{center}
$
\setlength{\unitlength}{1.2mm}
\hbox{
\begin{picture}(30,30)(-5,-5)
\put (0,0) {\line (1,0){20}}\put (10,0){\rightya}
\put (9, -5) {$\lift{k}_0$}
\put (20,0) {\line (0,1){20}}  \put (20,10){\upya}
\put (23, 9) {$u\mapsto \lift{\delta}_{1-u} (1) $}
\put (20,20) {\line (-1,0){20}}  \put (10,20){\rightya}
\put (9, 22.5) {$\lift{k}_1$}
\put (0,20) {\line (0,-1){20}} \put (0,10){\upya}
\put (-7, 9) {$1_{\lift{q}}$}
\put (9, 9.2) {$\lift{K}$}
\end{picture}
}
$
\end{center}
\caption{The map $\lift{K}$}\label{figK}
\end{figure}
Then we obtain a lift 
$$
\lift{k}_u:=\lift{K}|_{I\times\{u\}},
$$
of $k_u$, which is a path 
from   $\lift{q}\in F_{q}$ to the point
$\lift{\delta}_{1-u} (1)=h^L(1, 1-u)$ of $ F_{\delta_{1-u} (1)}$.
(See Figure~\ref{figLoopH}.)
%
\renewcommand{\PStextplot}[3]{\rlap{\hskip -240.000000 pt \hbox{\hskip #1pt  \raise #2pt \hbox{#3}}}}%
\begin{figure}
 \begin{center}
 \includegraphics{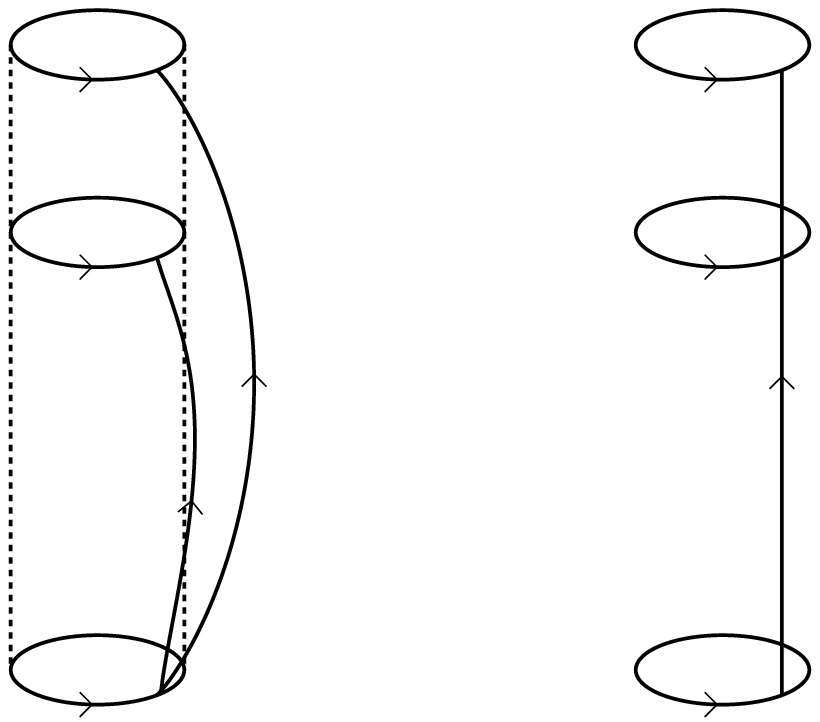}%
\PStextplot{247.000000}{194.000000}{${\bdre\delta}_1$}%
\PStextplot{247.000000}{14.000000}{${\bdre \delta}_0$}%
\PStextplot{247.000000}{140.000000}{${\bdre \delta}_{1-u}$}%
\PStextplot{-25.000000}{194.000000}{$\bdre \lift{\delta}_1$}%
\PStextplot{-25.000000}{14.000000}{$\bdre \lift{\delta}_0$}%
\PStextplot{-25.000000}{140.000000}{$\bdre \lift{\delta}_{1-u}$}%
\PStextplot{240.000000}{95.000000}{$k_0$}%
\PStextplot{88.000000}{95.000000}{$\lift{k}_0$}%
\PStextplot{47.000000}{75.000000}{$\lift{k}_u$}%
\PStextplot{120.000000}{100.000000}{$\maprightsp{}$}%
\PStextplot{130.000000}{111.000000}{$f$}%
\PStextplot{232.113678}{-2.289686}{$q$}%
\PStextplot{52.113678}{-2.289686}{$\lift{q}$}%
\end{center}
\caption{The loop $(H_u|_{\bdr I^2})^\sim$}%
\label{figLoopH}%
\end{figure}%
We then define a lift
$$
\map{\Lift{ H_u|_{\bdr I^2}}}{\bdr I^2}{X\spc}\qquad (u\in I)
$$
of $H_u|_{\bdr I^2}$ as a loop by 
$$
\lift{k}_u\cdot \bdre\lift{\delta}_{1-u}\cdot \lift{k}_u\inv \cdot \bdre\lift{\delta}_0 \inv,
$$
where $\bdr I^2$ is oriented and segmented as Figure~\ref{figorientation} above. 
%
%
%
%
 Then 
$(H_u, \Lift{ H_u|_{\bdr I^2}})$
is a homotopy of free loop pairs parametrized by $u\in I$.
By Remarks~\ref{rem:homotopymonrel} and~\ref{rem:homotopyD} again,
it is enough to prove that the free  loop pair
$$
\map{(H_1, \Lift{ H_1|_{\bdr I^2}})}{(I^2, \bdr I^2)}{(Y\spc, X\spc)}
$$
is of monodromy relation type around $\Sigma_i\spsh$.
Note that
$$
\Lift{ H_1|_{\bdr I^2}}=\lift{k}_1\cdot \bdre\lift{\delta}_{0}\cdot \lift{k}_1\inv \cdot \bdre\lift{\delta}_0 \inv,
$$
(see Figure~\ref{figLK1}),
and that the lift $\lift{k}_1$ of the constant map $k_1=1_q$ is a loop in $F_q$ with the base point $\lift{q}$.
%
\renewcommand{\PStextplot}[3]{\rlap{\hskip -126.000000 pt \hbox{\hskip #1pt  \raise #2pt \hbox{#3}}}}%
\begin{figure}
 \begin{center}
 \includegraphics{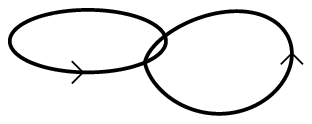}%
\PStextplot{0.000000}{47.700000}{$\bdre \lift{\delta}_0$}%
\PStextplot{116.100000}{42.300000}{$\lift{k}_1$}%
\qquad\quad
\hbox{
\setlength{\unitlength}{1.2mm}
\begin{picture}(30,30)(-5,-5)
\put (0,0) {\line (1,0){20}}\put (10,0){\rightya}
\put (8.5, -5) {$\bdre \lift{\delta}_0$}
\put (20,0) {\line (0,1){20}}  \put (20,10){\upya}
\put (23, 9) {$\lift{k}_1$}
\put (20,20) {\line (-1,0){20}}  \put (10,20){\rightya}
\put (8.5, 22.5) {$\bdre \lift{\delta}_0$}
\put (0,20) {\line (0,-1){20}} \put (0,10){\upya}
\put (-7, 9) {$\lift{k}_1$}
\end{picture}
}
\end{center}
\caption{Two figures for $\Lift{H_1|_{\bdr I^2}}=\lift{k}_1\cdot \bdre\lift{\delta}_{0}\cdot \lift{k}_1\inv \cdot \bdre\lift{\delta}_0 \inv$}%
\label{figLK1}%
\end{figure}%
Since $H_1(s, t)=H_0(s, 0)=\bdre\delta_0(s)$ for any $t$,
the pull-back
$$
\map{H_1\sp* (f\spc)}{H_1\sp* (X\spc)}{I^2}
$$
of $\shortmap{f\spc}{X\spc}{Y\spc}$ by $H_1$ is the product  of the pull-back 
$$
\map{(\bdre \delta_0) \sp* (f\spc)}{(\bdre \delta_0) \sp* (X\spc)}{I}
$$
of $f\spc$ by $\shortmap{\bdre \delta_0}{I}{Y\spc}$
and 
the identity map of the second factor $I$.
Let 
$$
\map{{}\sp{s}\Lift{ H_1|_{\bdr I^2}}}{\bdr I^2}{H_1\sp* (X\spc)=(\bdre \delta_0) \sp* (X\spc)\times I}
$$
be the cross-section of 
$H_1 \sp* (f\spc)$
over $\bdr I^2$
obtained from 
$\Lift{ H_1|_{\bdr I^2}}$.
We will describe the image of the free loop
${}\sp{s}\Lift{ H_1|_{\bdr I^2}}$
by a contraction 
$$
\map{c\sprime}{H_1\sp* (X\spc)}{F_q}.
$$
We construct the contraction $c\sprime$ as the  composite of 
the projection 
$$
\map{\pr_1}{\Lift{ H_1|_{\bdr I^2}}}{ (\bdre \delta_0) \sp* (X\spc)}
$$
onto the first factor
and a contraction 
$\shortmap{c}{(\bdre \delta_0) \sp* (X\spc)}{F_q}$. 
Let 
$$
\map{\sigma}{\bdr I^2}{(\bdre \delta_0) \sp* (X\spc)}
$$
be the composite of ${}\sp{s}\Lift{ H_1|_{\bdr I^2}}$ with the projection $\pr_1$.
The fibers $F_q^{(0)}$ and  $F_q^{(1)}$ of $\shortmap{(\bdre \delta_0 ) \sp* (f\spc)}{(\bdre \delta_0 ) \sp* (X\spc)}{I}$
over $0\in I$ and  $1\in I$ are canonically identified with $F_q$.
Let $\lift{q}^{(0)}\in F_q^{(0)}$ and $\lift{q}^{(1)}\in F_q^{(1)}$
be the points corresponding to $\lift{q}\in F_q$.
Then
$\Lift{ H_1|_{\bdr I^2}}|_{\{0\}\times I}$
(resp.~$\Lift{ H_1|_{\bdr I^2}}|_{\{1\}\times I}$)
gives rise 
 to  a loop $\lift{k}_1^{(0)}$ in $F_q^{(0)}$ 
 with the base point $\lift{q}^{(0)}$
 (resp.~a loop~$\lift{k}_1^{(1)}$ in $F_q^{(1)}$
with the base point $\lift{q}^{(1)}$).
Each of them
corresponds to the loop $\lift{k}_1$
by the obvious identifications
$(F_q, \lift{q})=(F_q^{(0)}, \lift{q}^{(0)})=(F_q^{(1)}, \lift{q}^{(1)})$.
On the other hand,  the loop $\bdre\lift{\delta}_0$ gives rise to a cross-section 
$$
\map{{}^{s} \bdre\lift{\delta}_0}{I}{(\bdr \delta_0 ) \sp* (X\spc)}
$$
of  
$(\bdre \delta_0 ) \sp* (f\spc)$
that connects $\lift{q}^{(0)}$ and $\lift{q}^{(1)}$.
The loop $\sigma$ on $(\bdre \delta_0 ) \sp* (X\spc)$ 
is then equal to the conjunction
$$
(\lift{k}_1^{(0)})\cdot ({}^{s} \bdre\lift{\delta}_0) \cdot  (\lift{k}_1^{(1)})\inv \cdot ({}^{s} \bdre\lift{\delta}_0)\inv.
$$
(See  Figure~\ref{figsigma}.)
\begin{figure}
\begin{center}
\setlength{\unitlength}{.9mm}
\begin{picture}(70,63)(-15, -4)
    \qbezier(0,50)(10,20)(0,20)
     \qbezier(0,50)(-10,20)(0,20)
     \put (0, 20) {\rightya}
     \put (-14,30) {$\lift{k}_1^{(0)}$}
     \put (0, 49.8) {\kuromaru}
      \put (-7, 49.8) {$\lift{q}^{(0)}$}
     \qbezier(40,50)(50,20)(40,20)
     \qbezier(40,50)(30,20)(40,20)
      \put (40, 20) {\rightya}
      \put (50,30) {$\lift{k}_1^{(1)}$}
     \put (40, 49.8) {\kuromaru}
     \put (43, 49.8) {$\lift{q}^{(1)}$}
     \put (0, 50) {\line (1,0){40}}
     \put (0, 3) {\line (1,0){40}}
     \put (0, 3) {\kuromaru}
     \put (40, 3) {\kuromaru}
         \put (20, 50) {\rightya}
     \put (16, 55) {$^s\bdre \lift{\delta}_0$}
     \put (19, 18) {\vector (0,-1){6}}
     \put (21, 14) {$(\bdre \delta_0)^*(f\spc)$}
      \put (20, 3) {\rightya}
     \put (17, -4) {$\bdre \delta_0$}
  \end{picture}
  \end{center}
  \caption{The loop $\sigma=(\lift{k}_1^{(0)})\cdot ({}^{s} \bdre\lift{\delta}_0) \cdot  (\lift{k}_1^{(1)})\inv \cdot ({}^{s} \bdre\lift{\delta}_0)\inv$}\label{figsigma}
  \end{figure}
We denote by $S\subset (\bdre \delta_0 ) \sp* (X\spc)$ the image of the section ${}^{s} \bdre\lift{\delta}_0$,
and choose a contraction 
$$
\map{c}{((\bdre \delta_0 ) \sp* (X\spc),  S)}{ (F_q^{(0)}, \lift{q}^{(0)})=(F_q, \lift{q})}
$$
 to the fiber over $0\in I$
that contracts  the section $S$ to the point $\lift{q}$.
We put
$$
\gamma:=\mu([\bdre \lift{\delta}_0]) \in \Aut(\pione (F_q, \lift{q})).
$$
By the definition of the lifted monodromy,
 the loop 
$$
({}^s \bdre\lift{\delta}_{0})\cdot (\lift{k}_1^{(1)})  \cdot ({}^s \bdre\lift{\delta}_0 )\inv
$$
on $\bdre \delta_0 \sp* (X\spc)$ is contracted  by $c$ to a loop in $F_q$ that represents
$$
[\lift{k}_1]^{(\gamma\inv)}\in \pione (F_q, \lift{q}),
$$
while the loop $\lift{k}_1^{(0)}$ on $F_q^{(0)}$
obviously represents 
$[\lift{k}_1]\in\;\pione (F_q, \lift{q})$.
Therefore, by the contraction $c$, the loop $\sigma$  on 
$(\bdre \delta_0) \sp* (X\spc)$ is mapped to a loop that represents
$$
[\lift{k}_1] ([\lift{k}_1]^{(\gamma\inv)})\inv =(\kappa\inv \kappa\sp{\gamma})\inv,
$$
where $\kappa:=([\lift{k}_1]^{(\gamma\inv)})\inv$.
Hence  the conjugacy class of $\pione (F_q, \lift{q})$
corresponding  to the free loop pair $(H_1, \Lift{ H_1|_{\bdr I^2}})$
is contained in the normal subgroup $N(\bdre \lift{\delta}_0)=N^{[\bdre {\delta}_0]}$
generated by the monodromy relations along $[\bdre {\delta}_0]$.
\end{proof}
\begin{corollary}\label{cor:A}
We put
\begin{eqnarray*}
\T\phantom{_\zeta}&:=&\set{(x,y,z)\in \R^3}{x^2+y^2\le 1, z\in I},\\
A_\zeta&:=&\set{(x,y,z)\in \T}{z=\zeta}, \quand\\
\Upsilon\,&:=&\set{(x,y,z)\in \T}{x^2+y^2=1}\cup A_1\;\;=\;\;\bdr \,\T\setminus A_0\spc,
\end{eqnarray*}
where $A_0\spc$ is the interior of the closed disc $A_0$.
Let $\shortmap{\varphi}{\T}{Y\spsh}$ be a continuous map
such that
$\varphi (\T)\cap \Sigma\spsh\subset \Sigma\spsh_i$ and 
$$
\varphi\inv (\Sigma\spsh_i)=\set{(x, 0, z)\in \T}{x^2+(z-1)^2=1/2}
$$
hold, and  such that
$\shortmap{\varphi|_{A_1}}{A_1}{Y\spsh}$
intersects $\Sigma\spsh$ transversely at  $(\pm 1/\sqrt{2}, 0,1)$.
Suppose that we have a lift
$\shortmap{\Lift{\varphi|_{\Upsilon}}}{\Upsilon}{X\spsh}$
of $\shortmap{\varphi|_{\Upsilon}}{\Upsilon}{Y\spsh}$
that intersects $\Theta\spsh_i$ transversely at the  two points $(\pm 1/\sqrt{2}, 0,1)$.
Let $\shortmap{\Lift{\varphi|_{\Upsilon}}|_{\bdr A_0}}{\bdr A_0}{X\spc}$ be 
the restriction of $\Lift{\varphi|_{\Upsilon}}$ to $\bdr \Upsilon=\bdr A_0$.
Then the free loop pair 
$$
\map{(\varphi|_{A_0}, \Lift{\varphi|_{\Upsilon}}|_{\bdr A_0})}{(A_0, \bdr A_0)}{(Y\spc, X\spc)}
$$
is of monodromy relation type around $\Sigma_i\spsh$.
\end{corollary}
\begin{corollary}\label{cor:Gamma}
Let $\shortmap{\delta}{\cunitdisc}{Y\spsh}$ be a transversal disc around $\Sigma_i\spsh$,
and let $\lift{\delta}$ and $\lift{\delta}\sprime$ be two transversal lifts of $\delta$.
We put $q:=\delta(1)$ and  
$\lift{q}:=\lift{\delta}(1)\in F_q$, 
$\lift{q}\sprime:=\lift{\delta}\sprime(1)\in F_q$.
Suppose that we are given a path
$\shortmap{\gamma_0}{I}{F_q}$
from $\lift{q}$ to $\lift{q}\sprime$.
Then we can deform $\gamma_0$ to a path $\gamma_t$ on $F_{\bdre\delta(t)}$ 
from $\bdre\lift\delta (t)$ to  $\bdre\lift\delta\sprime(t)$; that is, 
we have a continuous map
$\shortmap{\Gamma}{I\times I}{X\spsh}$
such that
$$
f(\Gamma(s, t))=\bdre\delta(t),
\quad
\Gamma(s, 0)=\gamma_{0}(s),
\quad
\Gamma(0, t)=\bdre\lift\delta (t),
\quad
\Gamma(1, t)=\bdre\lift\delta\sprime (t),
$$
and $\gamma_t:=\Gamma|_{I\times\{t\}}$.
Consider the path  $\gamma_1$
on ${F_q}$
from $\lift{q}$ to $\lift{q}\sprime$.
The conjunction  $\gamma_0^{\phantom{1}}\gamma_1\inv$ is a loop on $F_q$,
which we  write 
$\shortmap{\gamma_0^{\phantom{1}}\gamma_1\inv}{\bdr D}{F_q}$,
where $D$ is homeomorphic to $\cunitdisc$.
Then the free loop pair 
$$
\map{(1_q, \gamma_0^{\phantom{1}}\gamma_1\inv)}{(D, \bdr D)}{(Y\spc, X\spc)}
$$ 
is of monodromy relation type around $\Sigma_i\spsh$.
\end{corollary}
Now we start the proof of Theorem~\ref{thm:ZvK}.
\begin{proof}[Proof of Theorem~\ref{thm:ZvK}]
By Proposition~\ref{prop:relisinKer}, 
we have $N^{[\rho]}\subset \Ker (\iota_*)$  for any $[\rho]\in \LLL$,
because the lasso $\lambda(\lift{\rho})$ is null-homotopic in $X$
for any transversal lift $\lift{\rho}$ of $\rho$.
Therefore $\NNN\subset \Ker (\iota_*)$ follows.
\par
\medskip
Let a loop $\gamma: \IbI\to (\Fb,\tlb)$ represent 
an element $[\gamma]$ of $\Ker (\iota_*)$.
We will show that $[\gamma]\in \NNN$.
There exists a homotopy
$$
\map{h}{({I^2}, \sqcap)}{(X,\tlb)}
$$
from $\gamma$ to $1_{\tlb}$ in $X$
stationary on $\bdr I$;
that is,
$h|_{I\times\{0\}}=\gamma$ and $h|_{\sqcap}=1_{\tlb}$, 
where
$ \sqcap:=(\bdr I\times I)\cup (I\times \{1\})\subset I^2$.
By the condition \cond{C1}, we can perturb $h$ so that
\begin{equation}\label{eq:hempty}
h({I^2})\cap \Sing (f)=\emptyset
\end{equation}
 holds.
 Since $(f\circ h)|_{\bdr {I^2}}=1_b$,
 the map $f\circ h: I^2\to Y$ represents an element of $\pi_2 (Y, b)$.
 By the condition \cond{Z},
 we have a continuous map 
 $$
 \map{l}{({I^2},\bdr{I^2})}{(Z, b)}
 $$
 such that $[f\circ h]+[i_Z\circ l]=0$ holds in $\pi_2 (Y, b)$,
 where $i_Z: Z\inj Y$ is the inclusion.
 We then consider the continuous map
 $\shortmap{s_Z\circ i_Z\circ l}{({I^2}, \bdr{I^2})}{(X,\tlb)}$.
Replacing $h$ with 
$\shortmap{h\sprime}{({I^2}, \sqcap)}{(X,\tlb)}$
defined by
$$
h\sprime (x, y):=\begin{cases}
h(x, 2y) & \textrm{if $2y\le 1$,}\\
s_Z\circ  i_Z\circ l(x, 2y-1) & \textrm{if $2y\ge 1$,}
\end{cases}
$$
we have
\begin{equation}\label{eq:hpitwo}
[f\circ h]=0\quad\textrm{in}\quad \pi_2 (Y, b). 
\end{equation}
(See Figure~\ref{fighsprime}.) 
\begin{figure}
\begin{center}
$
\setlength{\unitlength}{1.5mm}
\vcenter{\hbox{
\begin{picture}(30,30)(-5,-5)
\put (0,0) {\line (1,0){20}}\put (10,0){\rightya}
\put (9, -3.3) {$\gamma$}
\put (20,0) {\line (0,1){20}}  \put (20,15){\upya} \put (20,5){\upya}
\put (22, 9) {$1_{\tlb}$}
\put (20,20) {\line (-1,0){20}}  \put (10,20){\rightya}
\put (9, 22.5) {$1_{\tlb}$}
\put (0,20) {\line (0,-1){20}}  \put (0,15){\upya} \put (0,5){\upya}
\put (-4, 9) {$1_{\tlb}$}
\put (8, 4) {$h$}
\put (5, 14.6) {$s_Z\circ  i_Z\circ l$}
\put (0,10) {\line (1,0){20}}\put (10,10){\rightya}
\put (12, 7.5) {$1_{\tlb}$}
\end{picture}
}}
$
\end{center}
\caption{The map $h\sprime$}\label{fighsprime}
\end{figure}
Moreover, since $s_Z(Z)\cap\Sing (f)=\emptyset$ by the condition~\cond{Z}, 
we still have~\eqref{eq:hempty}.
Then any small perturbation of $f\circ h$ can be lifted 
to a small perturbation of $h$.
Since $\Xi$ is of codimension 
$\ge 2$ in $Y$,
we can assume that $(f\circ h) ({I^2}) \cap \Sigma \subset  \Sigma\spsh$,
and that $f\circ h$ intersects $\Sigma\spsh$ transversely
(see Definition~\ref{def:defs1}).
We put
$$
(f\circ h)\inv (\Sigma\spsh)=\{P_1, \dots, P_n\}\;\subset\;  {I^2}\setminus \bdr I^2.
$$
%
We will construct a continuous map
$$
\map{j}{V:=I^2\setminus (D_1\spc\cup  \dots \cup D_m\spc)}{X\spsh}
$$
with the following properties:
\begin{itemize}
\item[\cond{j1}] $D_1, \dots, D_m$ are mutually disjoint closed discs in $I^2\setminus(\bdr I^2 \cup \{P_1, \dots, P_n\})$,
and $D_\mu\spc$ is the interior of $D_\mu$;
in particular,
$V$ contains $P_1$, \dots, $P_n$ in its interior, 
\item[\cond{j2}] $j(\bdr I^2)=\{\tlb\}$,
\item[\cond{j3}] $f\circ j=f\circ h|_V$ holds,  and hence we have 
$j\inv (\Theta\spsh)=\{P_1, \dots, P_n\}$,
\item[\cond{j4}] $j$ intersects $\Theta\spsh$ transversely at the points $P_\nu$ for $\nu=1, \dots, n$,  and
\item[\cond{j5}] for each $D_\mu$, the free loop pair
$$
\map{((f\circ h)|_{D_\mu}, j|_{\bdr D_\mu})}{(D_\mu, \bdr D_\mu)}{(Y\spc, X\spc)}
$$
is of monodromy relation type.
\end{itemize} 
By~\eqref{eq:hpitwo},
there exists a homotopy
$$
\map{H}{({I^2}\times I, B)}{(Y, b)}
$$
from $f\circ h$ to $1_b$ that is stationary on $\bdr{I^2}$;
that is, 
$H|_{I^2\times\{0\}}=f\circ h$ and $H|_{B}=1_b$,
where
$$
B:=(\bdr {I^2}\times I)\cup ({I^2}\times \{1\})\;\subset\; I^2\times I.
$$
Since $\Xi$ is of real codimension $\ge 4$ in $Y$, 
we can perturb $H$ and assume 
the following:
\begin{itemize}
\item[\cond{H1}] $H({I^2}\times I)\cap \Sigma$ is contained in $\Sigma\spsh$,
\item[\cond{H2}] $H$ intersects $\Sigma\spsh$ transversely
(in the sense of Definition~\ref{def:defs1}), 
so that 
$$
L:=H\inv (\Sigma\spsh)
$$
is a disjoint union of smooth real curves,
and
\item[\cond{H3}] the projection  $\pr_L: L\to I$  to the second factor  of $I^2\times I$
has only ordinary critical points in $L$;
that is, $\pr_L$ is a Morse function on $L$.
\end{itemize}
We have
$$
\bdr L=L\cap (I^2\times\{0\})=(f\circ h)\inv (\Sigma\spsh)=\{P_1, \dots, P_n\}.
$$
Let $L_1, \dots, L_k$ be the connected components of $L$.
Then  each  $L_{\kappa}$  is a curve connecting  two points 
of $\{P_1, \dots, P_n\}$,
or a curve without boundary.
 In particular, the cardinality $n$ of the points $(f\circ h)\inv (\Sigma\spsh)$ is even.
\par
\medskip
We denote by  $p_1^+, \dots, p_l^+$ (resp.~$p_1^-, \dots, p_m^-$)  the critical points  in $L\setminus \bdr L$
of  the projection $\pr_L: L\to I$ 
at which the Morse function $\pr_L$ attains a local maximum (resp.~a local minimum),
and call them the \emph{positive (resp.~negative) critical points of $\pr_L$}.
(See Figure~\ref{figLT},
in which $L$ is drawn in thick curve.)
\par
\medskip
Let $\T$ and $A_\zeta$ be as in Corollary~\ref{cor:A}.
For each negative critical point $p_\mu^-$,
we can choose a continuous map
$$
\map{\tau_\mu}{\T}{I^2\times I}
$$
with the following properties:
\begin{itemize}
\item[\cond{$\tau$1}]
each $\tau_\mu$ is a homeomorphism onto its image $T_\mu:=\tau_\mu(\T)$,
and $T_1, \dots, T_m$ are mutually disjoint, 
\item[\cond{$\tau$2}]
 there exists a strictly increasing function $t_\mu: I\to I$
 with $t_\mu (0)=0$ 
 that makes the following diagram  commutative;
$$
\begin{array}{ccc}
\renewcommand{\arraystretch}{1.2}
 \T & \maprightsp{\tau_\mu} & I^2\times I \\
 \mapdown && \mapdown \\
 I & \maprightsp{t_\mu} &\phantom{,}I,
\end{array}
$$ 
  where the vertical arrows  are the projections onto the last factors, 
\item[\cond{$\tau$3}]
$\tau_\mu\inv (\bdr (I^2\times I))=A_0$
and $\tau_\mu(A_0)\subset (I^2\setminus \bdr I^2)\times\{0\}$, 
\item[\cond{$\tau$4}]
$\tau_\mu\inv (L)=\shortset{(x,0,z)\in T}{x^2+(z-1)^2=1/2}$
and $\tau_\mu(1/2, 0, 1/2)=p_\mu^-$,
so that $p_\mu^-$ is the only critical point of $\pr_L$ in $T_\mu\cap L$,
 and
\item[\cond{$\tau$5}]
$\shortmap{H\circ(\tau_\mu|_{A_1})}{A_1}{Y\spsh}$
intersects $\Sigma\spsh$ transversely at $(\pm 1/\sqrt{2}, 0, 1)\in A_1$.
\end{itemize}
We put 
$$
T :=T_1\cup\dots\cup T_m.
$$
(In Figure~\ref{figLT}, each $T_\mu$ is depicted by dashed curves.)
We also put 
$$
\T\spc :=\shortset{(x,y,z)\in \T}{x^2+y^2<1, z<1}
$$
(the union of the interior of $\T$ and the bottom open disc),
and 
$$
T_\mu\spc:=\tau_\mu(\T\spc), \quad T\spc :=T_1\spc\cup\dots\cup T_m\spc
\quand
J:=(I^2\times I)\setminus T\spc.
$$
Note that $J$ is the closure of $(I^2\times I)\setminus T$.
Then  
$$
L\sprime:=L\cap J
$$
 is a disjoint union of 
smooth real curves $L\sprime_1, \dots, L\sprime_l$, 
and each connected component $L_\lambda\sprime$ of $L\sprime$
contains exactly one positive critical point $p_\lambda\sp +$
in $L_\lambda\sprime \setminus \bdr L_\lambda\sprime$.
Moreover,
each $L_\lambda\sprime$ has two boundary points $Q_\lambda$ and $Q_\lambda\sprime$,
each of which is either one point among $\{P_1, \dots, P_n\}$
or one of $\tau_\mu (\pm 1/\sqrt{2}, 0, 1)$ for some $\mu$.
%
%
\renewcommand{\PStextplot}[3]{\rlap{\hskip -290.250000 pt \hbox{\hskip #1pt  \raise #2pt \hbox{#3}}}}%
\begin{figure}
 \begin{center}
 \includegraphics{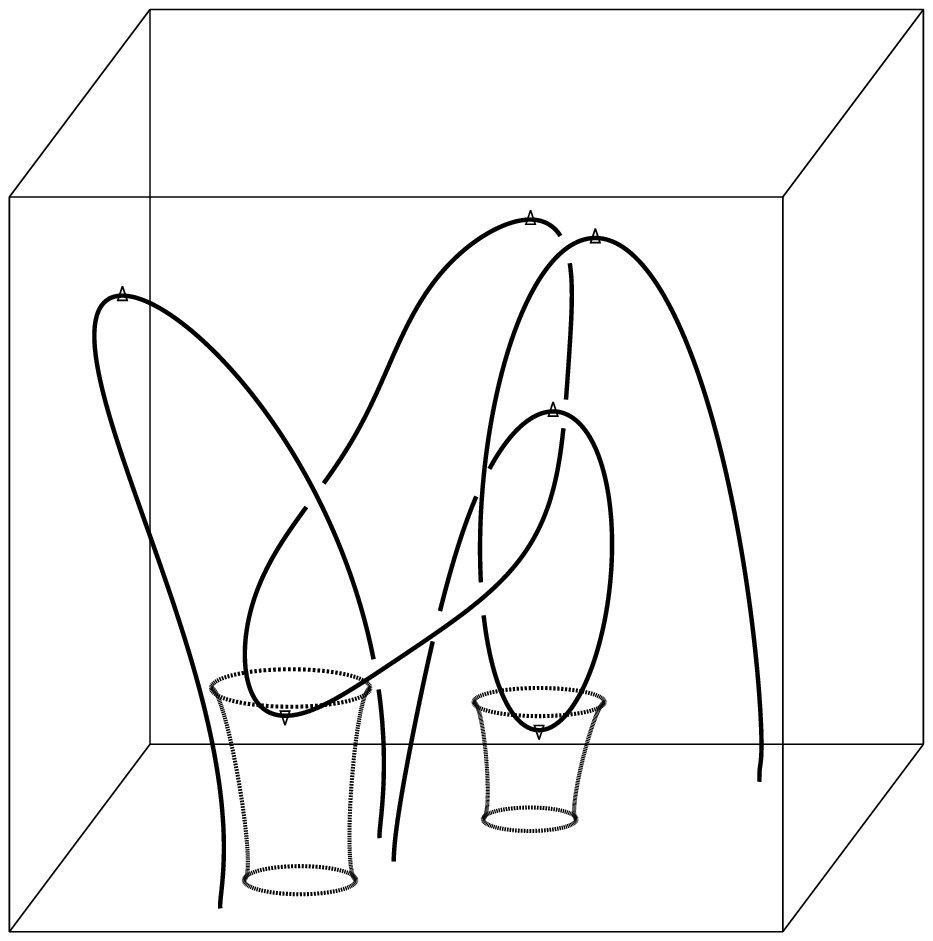}%
\PStextplot{76.950000}{-13.500000}{$\vartriangle$: the points $p_\lambda^+$, $\triangledown$: the points $p_\mu^-$.}%
\end{center}
\caption{$L$ and $T$}%
\label{figLT}%
\end{figure}%
If $Q_\lambda$ is one of $P_1, \dots, P_n$,
let $D(Q_\lambda)$ be a sufficiently small closed disc on $I^2\times\{0\}$
with the center $Q_\lambda$.
If $Q_\lambda$ is one of $\tau_\mu (\pm 1/\sqrt{2}, 0, 1))$,
let $D(Q_\lambda)$ be a sufficiently small closed disc on $\tau_\mu(A_1)$
with the center $Q_\lambda$.
We choose a closed disc $D(Q_\lambda\sprime)$
with the center $Q_\lambda\sprime$ in the same way.
Note that $\shortmap{H|_{D(Q_\lambda)}}{D(Q_\lambda)}{Y\spsh}$ and 
$\shortmap{H|_{D(Q_\lambda\sprime)}}{D(Q_\lambda)}{Y\spsh}$ are the transversal discs 
around the irreducible component 
$\Sigma\spsh_{i(\lambda)}$
of $\Sigma\spsh$
that contains $H(p_\lambda\sp+)$.
Then,
for each $\lambda=1, \dots, l$,  we have a tubular neighborhood
$$
\map{m_\lambda}{\cunitdisc\times I}{J}
$$
of $L_\lambda\sprime$ in $J$
with the following properties:
\begin{itemize}
\item[\cond{m1}]
each $m_\lambda$ is a homeomorphism onto its image $M_\lambda$,
and $M_1, \dots, M_l$ are mutually disjoint,
\item[\cond{m2}]
$m_\lambda\inv (L\sprime)=\{0\}\times I$
and 
$m_\lambda(\{0\}\times I)=L\sprime_\lambda$,
\item[\cond{m3}]
$m_\lambda$ is differentiable and locally a submersion at each point of $\{0\}\times I$, and 
\item[\cond{m4}]
$m_\lambda\inv (\bdr J)=\cunitdisc\times \bdr I$
and 
$m_\lambda(\cunitdisc \times \{0\})=D(Q_\lambda)$,
$m_\lambda(\cunitdisc \times \{1\})=D(Q_\lambda\sprime)$.
\end{itemize}
Then the composite
$\shortmap{H\circ m_\lambda}{\cunitdisc\times I}{Y\spsh}$
is an isotopy between the transversal discs
$H|_{D(Q_\lambda)}$ and $H|_{D(Q_\lambda\sprime)}$.
We put
$$
M:=M_1\cup \dots \cup M_l.
$$
Let $c_\lambda\in I$ be the real  number such that $m_\lambda(0, c_\lambda)=p_\lambda^+$.
We choose a point $p_\lambda^{+\prime}$ on $m_\lambda(\bdr\cunitdisc\times  \{c_\lambda\})\subset  \bdr {M_\lambda}$
and  a  path
$$
\map{w_\lambda}{I}{J}
$$
from  $p_\lambda^{+\prime}$ to  a point $p_\lambda^{+\prime\prime}$ of $ I^2\times\{1\}$
with the following properties:
\begin{itemize}
\item[\cond{w1}]
each $w_\lambda$ is a homeomorphism onto its image $W_\lambda$,
and $W_1, \dots, W_l$ are mutually disjoint,
\item[\cond{w2}]
$w_\lambda\inv (M)=\{0\}$, $w_\lambda\inv (\bdr J)=\{1\}$,
and 
\item[\cond{w3}]
the composite 
$\shortmap{\pr_2\circ w_\lambda}{I}{I}$ of $w_\lambda$ with the second
projection $I^2\times I\to I$ is strictly increasing.
\end{itemize}
We put 
$$
 W:=W_1\cup\dots\cup W_l.
$$
In Figure~\ref{figMW},
two of $M_{\lambda}\cup W_{\lambda}$ are illustrated.
The ceiling is $I^2\times \{1\}$,
from which $W_{\lambda}$  are dangling,
and the tubes are $M_\lambda$.
%
\renewcommand{\PStextplot}[3]{\rlap{\hskip -154.000000 pt \hbox{\hskip #1pt  \raise #2pt \hbox{#3}}}}%
\begin{figure}
 \begin{center}
 \includegraphics{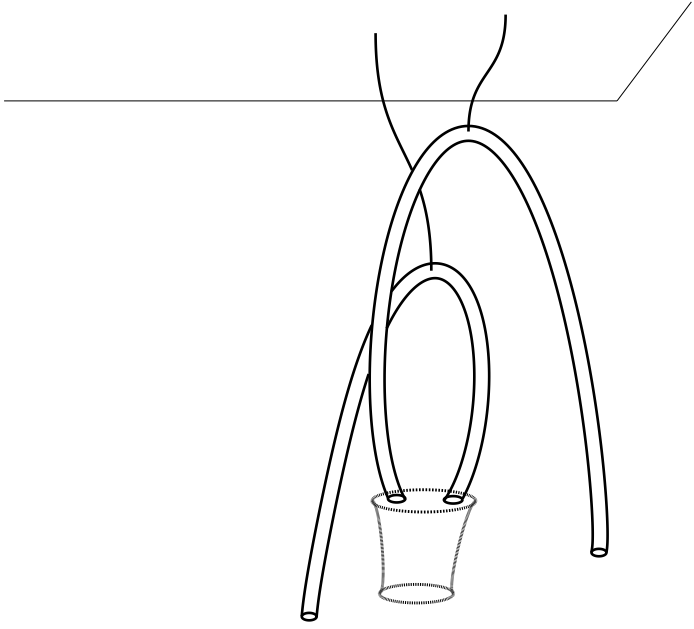}%
\end{center}
\caption{Two of $M_{\lambda}\cup W_{\lambda}$}%
\label{figMW}%
\end{figure}%
\par
\smallskip
The following fact is the crucial point in the construction of $j: V\to X\spsh$:
\begin{equation}\label{eq:sdr}
\textrm{\emph{$B\cup  M \cup W$ is a strong deformation retract of $J$}.}
\end{equation}
We choose  transversal lifts 
$\Lift{H|_{D(Q_\lambda)}}$ and $\Lift{H|_{D(Q_\lambda\sprime)}}$ of 
the transversal discs $H|_{D(Q_\lambda)}$ and $H|_{D(Q_\lambda\sprime)}$
around $\Sigma\spsh_{i(\lambda)}$, respectively.
Then 
the isotopy $\shortmap{H\circ m_\lambda}{\cunitdisc}{Y\spsh}$
 between 
$H|_{D(Q_\lambda)}$ and $H|_{D(Q_\lambda\sprime)}$
lifts to an isotopy 
between $\Lift{H|_{D(Q_\lambda)}}$ and $\Lift{H|_{D(Q_\lambda\sprime)}}$,
which yields 
 a lift $\Lift{H|_{M_\lambda}}$ of $H|_{M_\lambda}$.
Hence we obtain  a lift 
$$
\map{\Lift{H|_{M}}}{M}{X\spsh}
$$ 
of $H|_{M}$.
We define a lift $\Lift{H|_B}$ of  $H|_B$ to be the constant map $1_{\tlb}$.
Then we can lift the path $H\circ w_\lambda$ to
a path from $\Lift{H|_{M}} (p_\lambda^{+\prime})$ to $\Lift{H|_B}(p_\lambda^{+\prime\prime})=\tlb$,
and thus we obtain a lift 
$$
\map{\Lift{H|_W}}{W}{X\spsh}
$$ 
of $H|_W$.
Joining these three lifts together, we obtain a lift
$$
\map{\Lift{H|_{B\cup  M \cup W}}}{B\cup  M \cup W}{X\spsh}
$$
of $H|_{B\cup  M \cup W}$.
By the fact~\eqref{eq:sdr}, we can extend the lift $\Lift{H|_{B\cup  M \cup W}}$
to a lift 
$$
\map{\Lift{H|_J}}{J}{X\spsh}
$$ 
of $H|_J$,
because the pull-back
$(H|_J)^*(f\spsh)$ of $\shortmap{f\spsh}{ X\spsh}{Y\spsh}$ 
by $\shortmap{H|_J}{J}{Y\spsh}$ is locally trivial over the complement of the interior  of $M$ in $J$.
\par
\smallskip
Recall that the floor $I^2\times\{0\}$ of the source space $I^2\times I$ of $H$
is the source space $I^2$ of $f\circ h$.
For $\mu=1, \dots, m$, 
we put
$$
D_\mu:=\tau_{\mu}(A_0).
$$
These $D_1, \dots, D_m$ satisfy the condition~\cond{j1}.
Then  
$$
V:=I^2\setminus (D_1\spc\cup  \dots \cup D_m\spc)
$$
is identified with $J\cap(I^2\times \{0\})$.
We put
$$
j:=\Lift{H|_J}|_{V},
$$
which is a lift  of $f\circ h|_V=H|_{V}$.
Hence $j$ satisfies~\cond{j3}.
It is obvious that $j$ satisfies~\cond{j1} and~\cond{j2}.
Since $\Lift{H|_{M}}$ is constructed as a union of isotopies 
of transversal discs around $\Theta\spsh$,
the continuous map 
$$
\map{j|_{M\cap V}=\Lift{H|_{M}}|_{M\cap V}}{M\cap V}{X\spsh}
$$
intersects $\Theta\spsh$ transversely at each $P_{\nu}$.
Therefore  $j$ satisfies~\cond{j4}.
By the properties \cond{$\tau$4} and~\cond{$\tau$5} of $\tau_\mu$ and 
Corollary~\ref{cor:A}, 
we  see that  $j$ satisfies~\cond{j5}.
Thus we have constructed a continuous map $j: V\to X\spsh$
 which satisfies \cond{j1}\,-\,\cond{j5}, as is expected.
\par
\medskip
For $\nu=1, \dots, n$,
we choose a sufficiently small closed disc $D_{m+\nu}$ with the center $P_\nu$
in $I^2\setminus \bdr I^2$
in such a way that
the $m+n$ closed discs $D_1, \dots, D_{m+n}$ are mutually disjoint.

For each $\mu=1, \dots, m+n$, we choose a path
$$
\map{\alpha_{\mu}}{I}{{I^2}}
$$
from a point $R_\mu=(\rho_\mu, 0)\in I\times \{0\}$
 to a point $S_\mu\in \bdr D_\mu$
with the following properties:
\begin{itemize}
\item[\cond{$\alpha$1}] $0<\rho_{1}<\dots < \rho_{m+n}<1$, 
\item[\cond{$\alpha$2}] each $\alpha_\mu$ is injective and the images
$\alpha_\mu (I)$ $(\mu=1, \dots, m+n)$ are mutually  disjoint, and 
\item[\cond{$\alpha$3}]
$\alpha_\mu\inv (\bdr{I^2})=\{0\}$, $\alpha_\mu\inv ( D_\mu)=\{1\}$,
and $\alpha_\mu\inv (D_{\mu\sprime})=\emptyset$ if $\mu\ne \mu\sprime$.
\end{itemize}
In Figure~\ref{figAlphas},
the paths $\alpha_\mu$ are illustrated by thick curves.
%
\renewcommand{\PStextplot}[3]{\rlap{\hskip -253.000000 pt \hbox{\hskip #1pt  \raise #2pt \hbox{#3}}}}%
\begin{figure}
 \begin{center}
 \includegraphics{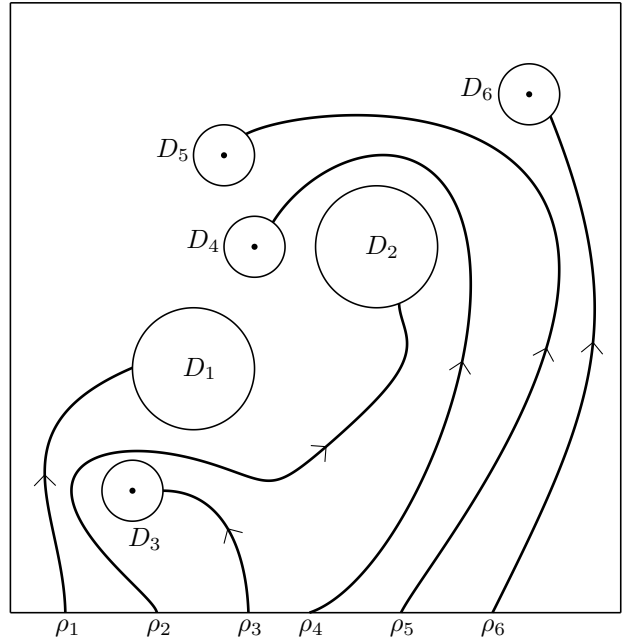}%
\PStextplot{87.400000}{89.700000}{$D_1$}%
\PStextplot{39.100000}{-6.900000}{$\rho_1$}%
\PStextplot{156.400000}{135.700000}{$D_2$}%
\PStextplot{73.600000}{-6.900000}{$\rho_2$}%
\PStextplot{66.700000}{25.300000}{$D_3$}%
\PStextplot{108.100000}{-6.900000}{$\rho_3$}%
\PStextplot{88.780000}{138.000000}{$D_4$}%
\PStextplot{131.100000}{-6.900000}{$\rho_4$}%
\PStextplot{77.280000}{172.500000}{$D_5$}%
\PStextplot{165.600000}{-6.900000}{$\rho_5$}%
\PStextplot{192.280000}{195.500000}{$D_6$}%
\PStextplot{200.100000}{-6.900000}{$\rho_6$}%
\end{center}
\caption{The paths $\alpha_\mu$}%
\label{figAlphas}%
\end{figure}%
Then there exists a continuous map
$$
\map{\ell}{{{\II}^2}}{{I^2}}
$$
with the following properties,
where $\II:=I=[0, 1]\subset \R$.
(We use the boldface $\II$ to distinguish the source plane  $\II^2$ and the target  plane $I^2$ of $\ell$.)
\begin{itemize}
\item[\cond{$\ell$1}]
$\ell$ induces a homeomorphism from ${\II}^2\setminus \bdr {\II}^2$
to  
$$
I^2\setminus \left(\bdr I^2 \cup \bigcup_{\mu=1}^{m+n} ( D_\mu\cup \alpha_{\mu}(I))\right),
$$
\item[\cond{$\ell$2}]
if $(x, y)\in \sqcap:=(\bdr {\II}\times {\II})\cup ({\II}\times \{1\})$, then
$\ell (x, y)=(x, y)$, and 
\item[\cond{$\ell$3}]
there exist real numbers $c_{\mu}, d_{\mu}, d_{\mu}\sprime, c_{\mu}\sprime\in \II$ 
for $\mu=1, \dots, m+n$ with
$$
\begin{array}{cccccccccccc}
0&<& c_{1}&<& d_{1}&<& d_{1}\sprime&<& c_{1}\sprime&<&\\
&<& c_{2}&<& d_{2}&<& d_{2}\sprime&<& c_{2}\sprime&<&\\
&&&&& \dots&&&&&\\
&<&c_{m+n}&<& d_{m+n}&<& d_{m+n}\sprime&<& c_{m+n}\sprime &<&1
\end{array}
$$
such that the following hold:
\begin{itemize}
\item 
$\ell (c_\mu, 0)=\ell (c_\mu\sprime, 0)=R_{\mu}\in I\times \{0\}$,
$\ell (d_\mu\sprime, 0)=\ell (d_\mu, 0)=S_{\mu}\in \bdr D_\mu$,
\item
$\ell|_{[c_\mu, d_\mu]\times\{0\}}$ is equal to
$\alpha_\mu$ via  a parameter change
$[c_\mu, d_\mu]\cong I$, and 
$\ell|_{[d_\mu\sprime, c_\mu\sprime]\times\{0\}}$ is equal to 
$\alpha_\mu\inv$ via  a parameter change
$[d_\mu\sprime, c_\mu\sprime]\cong I$,
\item
$\ell|_{[d_\mu, d_\mu\sprime]\times\{0\}}$ is  the loop 
that goes from $S_{\mu}$ to $S_{\mu}$ along $\bdr D_\mu$
clockwise, and
\item
$\ell|_{[c_{\mu-1}\sprime, c_\mu]\times\{0\}}$
is equal to the path
$[\rho_{\mu-1}, \rho_{\mu}]\to {I\times\{0\}}$
given by $t\mapsto (t, 0)$
via  a  parameter change $[c_{\mu-1}\sprime, c_\mu]\cong [\rho_{\mu-1}, \rho_{\mu}]$,
where we put $\rho_{0}:=0, c_{0}\sprime:=0$
and  $\rho_{m+n+1}:=1, c_{m+n+1}:=1$.
\end{itemize}
\end{itemize}
%
\renewcommand{\PStextplot}[3]{\rlap{\hskip -270.000000 pt \hbox{\hskip #1pt  \raise #2pt \hbox{#3}}}}%
\begin{figure}
 \begin{center}
 \includegraphics{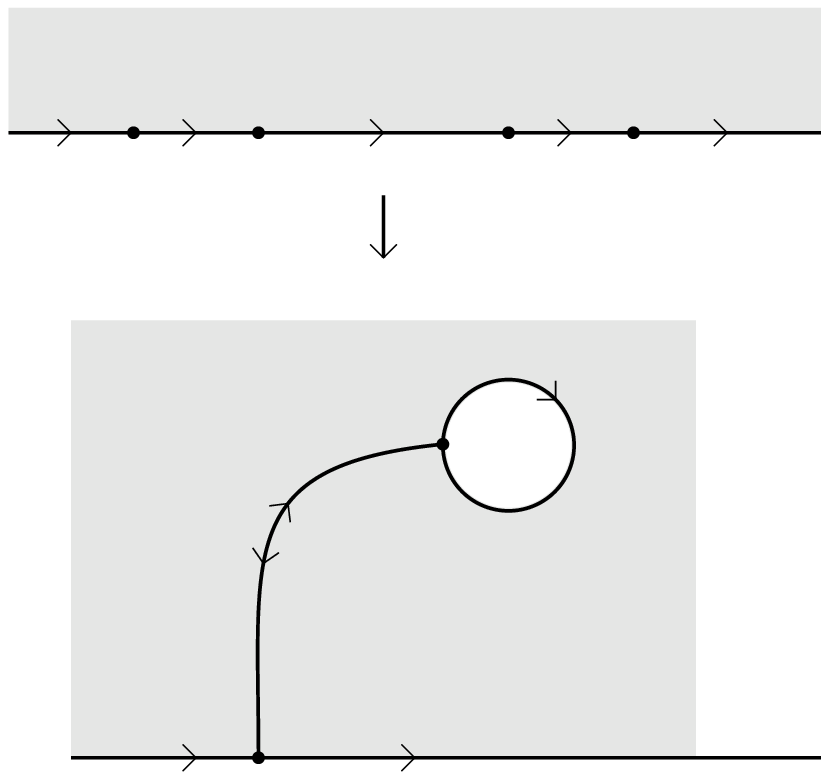}%
\PStextplot{50.400000}{196.200000}{$c_\mu$}%
\PStextplot{86.400000}{196.200000}{$d_\mu$}%
\PStextplot{194.400000}{196.200000}{$c_\mu\sprime$}%
\PStextplot{158.400000}{196.200000}{$d_\mu\sprime$}%
\PStextplot{82.800000}{-5.400000}{$R_\mu$}%
\PStextplot{130.502332}{108.296868}{$S_\mu$}%
\PStextplot{158.400000}{97.200000}{$D_\mu$}%
\PStextplot{131.400000}{162.000000}{$\ell$}%
\PStextplot{270.000000}{185.400000}{$\II\times\{0\}$}%
\PStextplot{270.000000}{5.400000}{$I\times\{0\}$}%
\end{center}
\caption{The map $\ell$}%
\label{figEll}%
\end{figure}%
(See Figure~\ref{figEll}.)
Since the image of $\ell$ is contained in $V$ and is disjoint from $\{P_1, \dots, P_n\}$,
we have  continuous maps
$$
\shortmap{j\circ \ell}{{\II}^2}{X\spc}
\quand
\shortmap{h\circ \ell}{{\II}^2}{X\spc}
$$
to $X\spc$.
They satisfy
$$
f\spc \circ j \circ \ell = f\spc \circ h \circ \ell
$$
by the property~\cond{j3}.
By  the properties~\cond{j2} and~\cond{$\ell$2},
they also satisfy
$$
j \circ \ell|_\sqcap =1_{\tlb}
\quand
h \circ \ell|_\sqcap =1_{\tlb}.
$$
We then define $G: {{\II}^2}\times {\II} \to Y\spc$ by the composition 
$$
G\;:\;
{{\II}^2}\times {\II}\;\maprightsp{\pr_1}\;
{{\II}^2} \;\maprightsp{f\spc \circ j \circ \ell = f\spc \circ h \circ \ell}\;
Y\spc,
$$
where $\pr_1$ is the first projection.
We put
$$
C:=
({{\II}^2}\times \bdr {\II})\cup (\sqcap \times {\II}) 
\;\subset\; 
{{\II}^2}\times {\II},
$$
and
define a lift
$$
\map{(G|_C)\splift}{C}{X\spc}
$$
of $G|_C: C\to Y\spc$ by the following:
$$
(G|_C)\splift (x, y, z)
:=
\begin{cases}
h(\ell(x, y)) & \textrm{if $z=0$, }\\
j(\ell(x, y)) & \textrm{if $z=1$, }\\
\tlb  & \textrm{if $(x, y, z)\in \sqcap \times {\II}$. }
\end{cases}
$$
Since $\shortmap{f\spc}{X\spc}{Y\spc}$ is locally trivial and 
$C$ is a strong deformation retract of ${{\II}^2}\times {\II}$,
the map $(G|_C)\splift$ extends to a lift
$$
\map{\lift{G}}{{{\II}^2}\times {\II}}{X\spc}
$$
of $\shortmap{G}{ {{\II}^2}\times {\II}}{Y\spc}$.
By construction,
for $(x, y)\in {{\II}^2}$, the restriction of $\lift{G}$
to $\{(x, y)\}\times {\II}$ is a path
in the fiber 
$$
F_{f\circ h\circ \ell (x, y)}=F_{f\circ j\circ \ell (x, y)}
$$
from the point $h\circ \ell (x, y)$ to the point $j\circ \ell (x, y)$.
For $x\in {\II}$, we put
$$
F_{[x]} :=F_{f\circ h\circ \ell (x, 0)}=F_{f\circ j\circ \ell (x, 0)},
\quand
\map{\xi_{[x]}:=\lift{G}|_{\{(x, 0)\}\times {\II}} }{{\II}}{F_{[x]}}.
$$
Suppose that $x\notin \textstyle{\bigcup}_{\mu=1}^{m+n}[c_{\mu}, c\sprime_{\mu}]$,
so that 
$$
(x\sprime, 0):=\ell (x, 0)\in I\times \{0\}.
$$
By~\cond{j2}, we see that
$F_{[x]}$ is equal to $F_b$ and
$\xi_{[x]}$ is a path in $F_b$ from $h(x\sprime, 0)=\gamma(x\sprime)$ to $j(x\sprime, 0)=\tlb$.
Moreover, we have
$\xi_{[0]}=\xi_{[1]}=1_{\tlb}$ because $\lift{G} |_{\sqcap\times\II}=1_{\tlb}$.
Therefore, for $\mu=0, 1, \dots, m+n$, the path
$$
\map{\gamma_{\mu}:=\gamma|_{[\rho_{\mu}, \rho_{\mu+1}]}=h|_{[\rho_{\mu}, \rho_{\mu+1}]\times\{0\}}}{[\rho_{\mu}, \rho_{\mu+1}]}{F_b}
$$ 
is homotopic to the path
$\xi_{[c_{\mu}\sprime]}\phantominv \xi_{[c_{\mu+1}]}\inv $ in $F_b$,
because the boundary of $\tilde{G}|_{[c\sprime_\mu, c_{\mu+1}]\times\{0\}\times\II}$ 
 is the loop
$\xi_{[c_{\mu}\sprime]}\phantominv \cdot 1_{\tlb} \cdot  \xi_{[c_{\mu+1}]}\inv \cdot  \gamma_{\mu}\inv$
in $F_b$,
where ${[c\sprime_\mu, c_{\mu+1}]\times\{0\}\times\II}\cong I^2$ is oriented and segmented as in Figure~\ref{figorientation}.
Since $\gamma$ is the conjunction 
$\gamma_{0}\gamma_{1}\dots\gamma_{m+n}$, 
the homotopy class $[\gamma]\in \pione (\Fb,\tlb)$ is equal to
$$
[\;\xi_{[c_{0}\sprime]}\phantominv \xi_{[c_{1}]}\inv 
\xi_{[c_{1}\sprime]}\phantominv \xi_{[c_{2}]}\inv
\dots \xi_{[c_{m+n}\sprime]}\phantominv \xi_{[c_{m+n+1}]}\inv\;]
\;\;=\;\;
[\xi_{[c_{1}]}\inv \xi_{[c_{1}\sprime]}\phantominv ]
\cdot
[\xi_{[c_{2}]}\inv \xi_{[c_{2}\sprime]}\phantominv ]
\cdot\cdots \cdot
[\xi_{[c_{m+n}]}\inv \xi_{[c_{m+n}\sprime]}\phantominv].
$$
(See Figure~\ref{figGammaXi}.)
Note that $\xi_{[c_{\mu}]}\inv \xi_{[c_{\mu}\sprime]}\phantominv$
is a loop in $F_b$ with the base point $\tlb$.
It is enough to show that
each 
$[\xi_{[c_{\mu}]}\inv \xi_{[c_{\mu}\sprime]}\phantominv]\in \pione (F_b, \tlb)$
is contained in $N^{[\rho]}$
for some transversal disc $\rho$ around an irreducible component 
of $\Sigma\spsh$.
%
\renewcommand{\PStextplot}[3]{\rlap{\hskip -220.500000 pt \hbox{\hskip #1pt  \raise #2pt \hbox{#3}}}}%
\begin{figure}
 \begin{center}
 \includegraphics{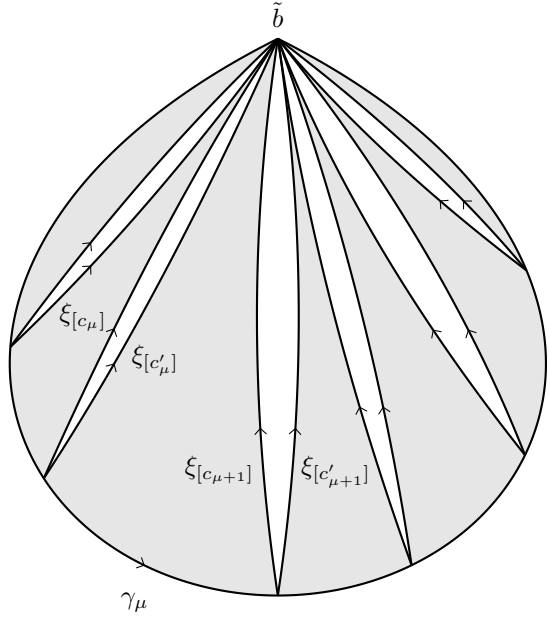}%
\PStextplot{55.650000}{7.350000}{$\gamma_{\mu}$}%
\PStextplot{32.550000}{114.450000}{$\xi_{[c_{\mu}]}$}%
\PStextplot{59.850000}{98.700000}{$\xi_{[c_{\mu}\sprime]}$}%
\PStextplot{79.800000}{56.700000}{$\xi_{[c_{\mu+1}]}$}%
\PStextplot{123.900000}{56.700000}{$\xi_{[c_{\mu+1}\sprime]}$}%
\PStextplot{112.980000}{225.750000}{$\tlb$}%
\end{center}
\caption{The paths $\gamma_{\mu}$ and $\xi_{[c_\mu]},\xi_{[c_\mu\sprime]}$}%
\label{figGammaXi}%
\end{figure}%
\par
\medskip
Consider the path 
$$
\map{\lift{\alpha}_\mu:=j\circ \alpha_\mu}{I}{X\spc}
$$
from $\tlb$ to $\lift{q}_\mu:=j(S_\mu)\in F_{q_\mu}$,
where $q_\mu:=f(j(S_\mu))=f(h(S_\mu))$,
and the induced isomorphism
$$
\mapisom{[\lift{\alpha}_\mu]_*}{\pione(F_b, \tlb)}{\pione(F_{q_\mu}, \lift{q}_\mu)}.
$$
This isomorphism maps
$[\xi_{[c_{\mu}]}\inv \xi_{[c_{\mu}\sprime]}\phantominv]\in \pione (F_b, \tlb)$
to
$$
[\xi_{[d_{\mu}]}\inv \xi_{[d_{\mu}\sprime]}\phantominv]\in \pione (F_{q_\mu}, \lift{q}_\mu).
$$
(See Figure~\ref{figLiftAlphas}.)
We  consider $\xi_{[d_{\mu}]}\inv \xi_{[d_{\mu}\sprime]}\phantominv$ as a free loop 
$\bdr\cunitdisc\to F_{q_\mu}$ in $F_{q_\mu}$.
It is enough to show that the free loop pair
$$
\map{(1_{q_\mu}, \xi_{[d_{\mu}]}\inv \xi_{[d_{\mu}\sprime]}\phantominv)}{(\cunitdisc, \bdr \cunitdisc)}{(Y\spc, X\spc)}
$$
is of monodromy relation type.
%
\renewcommand{\PStextplot}[3]{\rlap{\hskip -238.000000 pt \hbox{\hskip #1pt  \raise #2pt \hbox{#3}}}}%
\begin{figure}
 \begin{center}
 \includegraphics{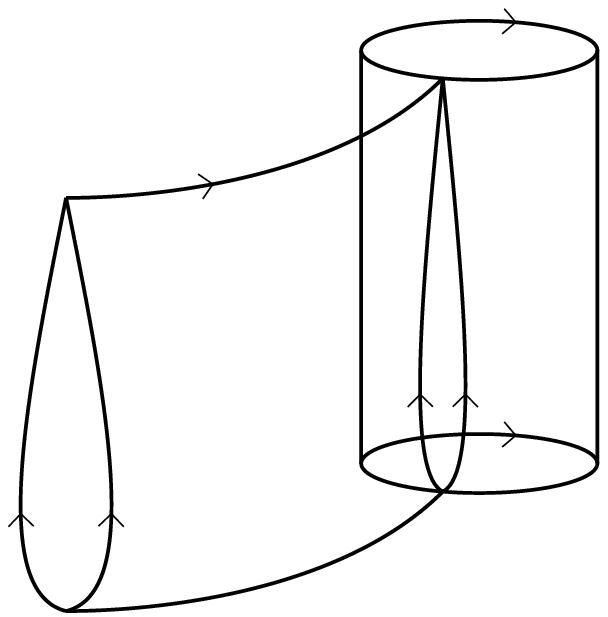}%
\PStextplot{49.300000}{125.800000}{$\tlb$}%
\PStextplot{16.150000}{42.075000}{$\xi_{c_[{\mu}]}$}%
\PStextplot{65.535000}{42.075000}{$\xi_{[c_{\mu}\sprime]}$}%
\PStextplot{136.679422}{123.224020}{$\xi_{[d_{\mu}]}$}%
\PStextplot{165.307422}{123.224020}{$\xi_{[d_{\mu}\sprime]}$}%
\PStextplot{93.269685}{134.854439}{$\lift{\alpha}_{\mu}$}%
\PStextplot{157.793422}{160.216020}{$\lift{q}_{\mu}$}%
\PStextplot{207.400000}{163.200000}{$j|_{\bdr D_{\mu}}$}%
\PStextplot{207.400000}{44.200000}{$h|_{\bdr D_{\mu}}$}%
\end{center}
\caption{Deformation of the loop along $\lift{\alpha}_\mu$}%
\label{figLiftAlphas}%
\end{figure}%
\par
\medskip
Suppose that $\mu>m$, so that $D_{\mu}$ is a  disc with the center $P_{\mu-m}
\in (f\circ h)\inv (\Sigma\spsh)$.
Then $(1_{q_\mu}, \xi_{[d_{\mu}]}\inv \xi_{[d_{\mu}\sprime]}\phantominv)$ is of monodromy relation type
by Corollary~\ref{cor:Gamma}.
Suppose that $\mu\le m$.
By~\cond{j5},
it is enough to show that the free loop pair 
$(1_{q_\mu}, \xi_{[d_{\mu}]}\inv \xi_{[d_{\mu}\sprime]}\phantominv)$ is homotopic to 
the free loop pair
$$
\map{((f\circ h)|_{D_\mu}, j|_{\bdr D_\mu})}{(D_\mu, \bdr D_\mu)}{(Y\spc, X\spc)}
$$
under a suitable homeomorphism  $\cunitdisc\cong D_\mu$.
We put
$$
\map{l_\mu:=\ell_{[d_\mu, d_\mu\sprime]\times\{0\}}}{[d_\mu, d_\mu\sprime]}{\bdr D_\mu}.
$$
Consider the continuous map 
$$
\map{\zeta_\mu}{[d_\mu, d_\mu\sprime]\times I}{X\spc}
$$
given by $\zeta_\mu(x, t):=\xi_{[x]}(t)$.
With the base point and the orientation on the boundary of $[d_\mu, d_\mu\sprime]\times I$
given in Figure~\ref{figdd},
the boundary of $\zeta_\mu$ is equal to the loop 
$$
\xi_{[d_\mu]}\inv \cdot (h\circ l_\mu)\cdot \xi_{[d\sprime _\mu]}\phantominv \cdot (j\circ l_\mu)\inv
$$
with the base point $\lift{q}_\mu$.
Since the free loop $h\circ l_\mu$ 
is the boundary of $h|_{D_\mu}$,
it is  null-homotopic in $X\spc$.
Hence the free loop $\xi_{[d_\mu]}\inv \cdot \xi_{[d\sprime _\mu]}\phantominv$
is homotopic to the free loop $j\circ l_\mu$ in $X\spc$.
It can be easily seen that we can construct a  homotopy 
of free loops from 
$j|_{\bdr D_\mu}=j\circ l_\mu$ to $\xi_{[d_\mu]}\inv \cdot \xi_{[d\sprime _\mu]}\phantominv$
in $X\spc$ as a lift of the restriction to $\bdr D_\mu$
of  a contraction from $f(h(D_\mu))$ to $q_\mu$,
because $f(h(D_\mu)) \subset Y\spc$ holds for $\mu\le m$.
Hence $(1_{q_\mu}, \xi_{[d_{\mu}]}\inv \xi_{[d_{\mu}\sprime]}\phantominv)$ is homotopic to 
$((f\circ h)|_{D_\mu}, j|_{\bdr D_\mu})$.
\begin{figure}
\begin{center}
\setlength{\unitlength}{1.2mm}
\begin{picture}(34,38)(-12,-2)
\put (0,0) {\line (1,0){17}}\put (8.5,0){\rightya}
\put  (17, 0) {\;\;$(d_\mu\sprime,0)$}
\put (17,0) {\line (0,1){30}}  \put (17,15){\upya}
\put  (17, 30) {\;\;$(d_\mu\sprime,1)$}
\put (17,30) {\line (-1,0){17}}  \put (8.5,30){\leftya}
\put  (0, 30) {\llap{the base point  $(d_\mu,1)$\;\;}}
\put (0,30) {\line (0,-1){30}} \put (0,15){\downya}
\put  (0, 0) {\llap{$(d_\mu,0)$\;\;}}
\put  (0, 29.9) {\kuromaru}
\put (2.0, 14.2) {$[d_\mu, d_\mu\sprime]\times I$}
\put (-7, 15){$\xi_{[d_\mu]}\inv$}
\put (20, 15){$\xi_{[d_\mu\sprime]}\phantominv$}
\put (6, -4){$h\circ l_\mu$}
\put (5, 33){$(j\circ l_\mu)\inv$}
\end{picture}
\end{center}
\caption{The orientation of $\bdr ([d_\mu, d_\mu\sprime]\times I)$}\label{figdd}
\end{figure}
\end{proof}
The following is a semi-classical version of Theorem~\ref{thm:ZvK}.
\begin{theorem}\label{thm:C}
Suppose that the conditions~\cond{C1} and~\cond{C2} are satisfied.
Suppose also that there exist a reduced connected curve $C$ 
(possibly singular and/or  reducible and not necessarily closed) on $Y$ 
and a continuous cross-section 
$$
\map{s_C}{C}{f\inv (C)}
$$
of $f$ over $C$
with the following properties:
\begin{itemize}
\item 
$C\spc:=C\cap Y\spc$ is non-empty and connected, and the inclusion
$C\spc\inj Y\spc$ induces a surjection
$\pione(C\spc, b)\surj \pione(Y\spc, b)$,
where $b\in C\spc$ is a base point,
\item 
the inclusion $C\inj Y$ induces a surjection
$\pi_2(C, b)\surj \pi_2(Y, b)$, 
\item $s_C(C)\cap \Sing (f)=\emptyset$, and 
\item
for each irreducible component $\Sigma_i$ of $\Sigma$ with codimension $1$ in $Y$,
there exists  a  point $p_i\in C\cap \Sigma_i$
satisfying  the following:
\begin{itemize}
\item  $C$ and $\Sigma$ are smooth  at $p_i$,
and $C$ intersects $\Sigma_i$  transversely at $p_i$, 
\item  the cross-section $s_C$ is holomorphic at $p_i$.
\end{itemize}
\end{itemize}
By the cross-section $s_C$, we have the classical monodromy action
$$
\pione (C\spc, b)\;\to\;\Aut (\pione (F_b, \tlb)),
\quad\rmwhere\quad \tlb:=s_C(b)\in F_b:=f\inv (b),
$$
which we denote by $g\mapsto g^u$ for $u\in \pione (C\spc, b)$.
Then $\Ker (\iota_*)$ is  equal to
$$
N_K:=\gen{\;\shortset{g\inv g^u}{g\in \pione (F_b, \tlb), u\in K}\;}, 
$$
where $K\subset \pione (C\spc, b)$
is the kernel of $\pione(C\spc, b)\to \pione(C, b)$
induced by the inclusion. 
\end{theorem}
\begin{proof}
First of all, remark that 
the condition \cond{Z} is
satisfied 
with $C$ and $s_C$ being $Z$ and $s_Z$ in the condition \cond{Z},
and hence $\Ker (\iota_*)$ is equal to $\NNN$.
\par
\medskip
Let $\shortmap{\gamma}{(I, \bdr I)}{(C\spc, b)}$ be a loop that
represents an element $u$ of $K$.
We have a homotopy (stationary on $\bdr I$) $h$ on $C$ from $\gamma$ to $1_b$.
Then $s_C\circ h$ 
is a homotopy on $X$ from $s_C\circ \gamma$ to $1_{\tlb}$.
By definition,
the classical monodromy action by $u$
is equal to the lifted monodromy action by
 $[s_C\circ \gamma]\in \pione(X\spc, \tlb)$.
Since  $s_C\circ \gamma$ is  null-homotopic in $X$, we see  that 
$g\inv g^u=g\inv g^{\mu([s_C\circ \gamma])}$
is contained in $\Ker (\iota_*)$ by Proposition~\ref{prop:relisinKer}.
Thus $N_K\subset \Ker(\iota_*)$ is proved.
\par
\medskip
In order to prove $\NNN=\Ker(\iota_*)\subset N_K$,
it is enough to show that,
for any leashed disc $\rho=(\delta, \eta)$
around an irreducible component $\Sigma\spsh_i$ of $\Sigma\spsh$ in
$Y\spsh$,
the normal subgroup
$N^{[\rho]}$ is contained in $N_K$.
We have a point $p_i$ of $C\cap \Sigma_i$
at which $C$ and $\Sigma$ are smooth and intersect transversely.
Let
$$
\mapinj{\delta_{i, C}}{\cunitdisc}{C}
$$
be a sufficiently small closed disc on $C$ 
such that $\delta_{i, C} (0)=p_i$.
Since $s_C$ is holomorphic at $p_i$
and $s_C(p_i)\notin \Sing (f)$ by the assumption,
$\Theta:=f\inv (\Sigma)$ is  smooth  at $s_C(p_i)$,
and $s_C\circ \delta_{i, C}$ intersects $\Theta$ at $s_C(p_i)$ transversely.
If $p_i\in \Xi$,
then we perturb $\delta_{i, C}$ to 
a $\Cinf$-map $\shortmap{\delta_{i, C}\sprime}{\cunitdisc}{Y\spsh}$
such that
$\delta_{i, C}|_{\bdr\cunitdisc}=\delta_{i, C}\sprime|_{\bdr\cunitdisc}$.
If $p_i\notin \Xi$, then we put $\delta_{i, C}\sprime:=\delta_{i, C}$.
Then $\delta_{i, C}\sprime$ is a transversal disc around $\Sigma_i\spsh$
such that $\delta_{i, C}\sprime (\bdr\cunitdisc)\subset C\spc$.
Since $s_C(p_i)\notin \Sing (f)$,
we can lift the perturbation from $\delta_{i, C}$ to $\delta_{i, C}\sprime$
to a perturbation
from $s_C\circ \delta_{i, C}$ to 
$$
\mapinj{\lift{\delta}_{i, C}\sprime}{\cunitdisc}{X\spsh}
$$
in such a way that
$$
\lift{\delta}_{i, C}\sprime|_{\bdr\cunitdisc}
=s_C\circ {\delta}_{i, C}\sprime|_{\bdr\cunitdisc}
=s_C\circ {\delta}_{i, C}|_{\bdr\cunitdisc},
$$
and that 
$\lift{\delta}_{i, C}\sprime$ is a transversal lift of $\delta_{i, C}\sprime$
around $\Theta\spsh_i$.
The transversal disc $\delta$ of the given leashed disc $\rho=(\delta, \eta)$
is isotopic  to $\delta_{i, C}\sprime$ (Proposition~\ref{prop:rho}).
Hence $\rho$ is isotopic to a leashed disc
$$
\rho\sprime=(\delta_{i, C}\sprime, \eta\sprime)
$$
for some path $\eta\sprime$ on $Y\spc$ from $\delta_{i, C}(1)=\delta_{i, C}\sprime(1)\in C\spc$ to $b$.
Since $C\spc$ is connected,
there exists a path $\zeta$ on $C\spc$ from $b$ to $\eta\sprime(0)=\delta_{i, C}(1)$.
Then $\zeta\eta\sprime$ is a loop on $Y\spc$
with the base point $b$.
Since the inclusion $C\spc\inj Y\spc$
induces a surjection $\pione (C\spc, b)\surj \pione (Y\spc, b)$,
there exists a loop $\xi$ on $C\spc$ with the base point $b$
that is homotopic to $\zeta\eta\sprime$  in $Y\spc$.
Then $\rho=(\delta, \eta)$ is isotopic to the leashed disc
$$
\rho_C:=(\delta_{i, C}\sprime, \zeta\inv\xi).
$$
Note that $\zeta\inv\xi$ is a path on $C\spc$.
Since $\lift{\delta}_{i, C}\sprime(1)=s_C(\delta_{i, C}\sprime(1))$, the pair  
$$
\lift{\rho}_C:=(\lift{\delta}_{i, C}\sprime, s_C\circ(\zeta\inv\xi))
$$
is a leashed disc,
which is a transversal lift of $\rho_C$.
Hence $N^{[\rho]}$ is generated by the monodromy relations 
$g\inv g^{\mu([\lambda(\lift{\rho}_C)])}$ along $[\lambda(\lift{\rho}_C)]$.
Note that the lasso $\lambda (\rho_C)$ is a loop on $C\spc$
that is null-homotopic in $C$,
so that we have
$[\lambda (\rho_C)]\in K$.
Because 
 $s_C\circ \lambda (\rho_C)=\lambda(\lift{\rho}_C)$,
 the generators $g\inv g^{\mu([\lambda(\lift{\rho}_C)])}$ 
of $N^{[\rho]}$ are contained in $N_K$.
\end{proof}
We give a sufficient condition under which 
 $N^{[\rho]}=1$ holds
for one (and hence any) leashed disc $\rho$ around $\Sigma_i\spsh$.
(See Corollary~\ref{cor:foranyrho}.)
%
\par
\medskip
Suppose that 
$X$ is the complement to a reduced  hypersurface
$W$ in a smooth variety $\ol{X}$,
and that $f$ is the restriction to $X$
of a  \emph{projective} morphism $\bar{f}: \ol{X}\to Y$.
For $y\in Y$,
we put $\ol{F}_y:=\bar{f}\inv (y)$,
and denote by $W_y$ the \emph{scheme-theoretic} intersection of $\ol{F}_y$ with $W$.
Let $\Sing (\bar{f})\subset \ol{X}$ be the Zariski closed subset of
critical points of $\bar{f}$.
\begin{proposition}\label{prop:proj}
We assume the conditions~\cond{C1} and~\cond{C2}.
Suppose that,
for a general point $y$ of $\Sigma_i$,
the intersection $\ol{F}_{y}\cap \Sing (\bar{f})$ 
is of codimension $\ge 2$ in $\ol{F}_{y}$
and $W_{y}\setminus (W_{y}\cap \Sing (\bar{f}))$
is a reduced hypersurface of $\ol{F}_{y}\setminus (\ol{F}_y\cap \Sing (\bar{f}))$.
Then $N^{[\rho]}=1$ holds
for a leashed disc $\rho$ around $\Sigma_i\spsh$.
\end{proposition}
\begin{proof}
Let $y_0$ be a general point $y_0$ of $\Sigma_i$,
and 
let $U\subset Y$ be a sufficiently small  
contractible neighborhood of $y_0$.
Since $\bar f$ is  projective, 
there exists an embedding over $U$
of $\bar{f}\inv (U)$ into $\P^N\times U$;
$$
\begin{array}{ccccc}
\bar{f}\inv (U)  && \inj && \P^N\times U \\
&\searrow &&\swarrow & \\
&&U.&&
\end{array}
$$
By this embedding,
we consider each $\ol{F}_y$ for $y\in U$ as a closed subscheme of $\P^N$
of dimension $\dim X-\dim Y$.
We choose a general linear subspace $P\subset \P^N$
of codimension $\dim \ol{F}_y-1$.
By the assumption
$\dim (\ol{F}_y\cap \Sing (\bar{f}))\le \dim \ol{F}_y -2$ for any $y\in U\cap\Sigma_i$,
we have $(P\times U)\cap \Sing (\bar f) =\emptyset$ and we can assume that 
$P\cap\ol{F}_y$
is a smooth projective curve for any $y\in U$.
By the  assumption on $W_y$,
we see that $P\cap W_y$ is a reduced divisor of $P\cap\ol{F}_y$
whose degree  is independent of $y\in U$.
Hence the family
$$
P\cap F_y=P\cap (\ol{F}_y\setminus W_y)\qquad (y\in U)
$$
of punctured Riemann surfaces is trivial (in the $\Cinf$-category) over $U$.
Let $\shortmap{\delta}{\cunitdisc}{Y\spsh}$
be  a transversal disc around $\Sigma_i\spsh$
such that $\delta(\cunitdisc)\subset U$.
Then we have a transversal lift
$\shortmap{\lift{\delta}}{\cunitdisc}{X\spsh}$
of $\delta$
such that $\lift{\delta}(z)\in P\cap F_{\delta(z)}$
holds for any $z\in \cunitdisc$.
We put
$$
q:=\delta (1),
\qquad \lift{q}:=\lift{\delta}(1)\in P\cap F_{q}.
$$
The lifted monodromy of $[\bdre \lift\delta]$ on $\pione(P\cap F_{q}, \lift{q})$
is trivial.
On the other hand,
the inclusion $P\cap F_{q}\inj  F_{q}$ induces a surjective homomorphism
$$
\pione (P\cap F_{q}, \lift{q}) \surj \pione (  F_{q}, \lift{q})
$$
by the Lefschetz-Zariski hyperplane section theorem.
(See, for example, ~\cite{MR932724} or~\cite{MR820315}).
Hence the lifted monodromy of $[\bdre \lift\delta]$ on $\pione(F_{q}, \lift{q})$
is also trivial.
\end{proof}
We prove the two corollaries stated in Introduction.
\begin{proof}[Proof of Corollary~\ref{cor:RRReq}]
Since the lasso of any transversal lift of a leashed disc on $Y\spsh$
around   $\Sigma_i\spsh$ is null-homotopic in $X$,
we have $\NNN\subset \RRR$.
Hence 
Corollary~\ref{cor:RRReq} follows
from Theorem~\ref{thm:ZvK}, Proposition~\ref{prop:relisinKer} and Nori's lemma~(Proposition~\ref{prop:nori} and Remark~\ref{rem:C0C3}).
\end{proof}
\begin{proof}[Proof of Corollary~\ref{cor:proj}]
It is enough to show that $f$ satisfies the condition~\cond{C2},
and that,
for each $\Sigma_i$, $N^{[\rho]}=1$ holds
for a leashed disc $\rho$ around $\Sigma_i\spsh$.

Since $f$ is projective and the general fiber is connected,
every fiber of $f$ is non-empty and connected.
Suppose that $F_y$ is reducible for a general  point $y$
of some irreducible hypersurface $\Sigma\sprime$ of $Y$.
Let  $\unitdisc\subset Y$ be a small open disc 
intersecting $\Sigma\sprime$ transversely at $y$
such that $f\inv (\unitdisc)$ is smooth.
Then $F_y$ is a reducible hypersurface of $f\inv (\unitdisc)$.
Since $F_y$ is connected and projective,
there exist distinct irreducible components $F_y\sprime$ and $F_y\spprime$ 
of $F_y$ that intersect.
Since $F_y\sprime\cap F_y\spprime$ is of codimension $2$ in  $f\inv (\unitdisc)$,
we obtain a contradiction to the assumption that $\Sing (f)$ is of codimension $\ge 3$ in $X$.
Thus the condition~\cond{C2} is satisfied.

Let $y$ be a general point of $\Sigma_i$.
By the assumption that  $\Sing (f)$ is of codimension $\ge 3$ in $X$,
we see that $F_y\cap \Sing (f)$ is of codimension $\ge 2$ in $F_y$.
Applying Proposition~\ref{prop:proj}
to the case  where $W=\emptyset$ and $X=\ol{X}$, we obtain $N^{[\rho]}=1$ 
for a leashed disc $\rho$ around $\Sigma_i\spsh$.
\end{proof}
\section{Proof of Theorem~\ref{thm:ULZvK}}\label{sec:proof1}
\begin{proof}[Proof of Theorem~\ref{thm:ULZvK}]
We assume $k\le n-2$, where $n$ is the dimension of the 
smooth non-degenerate projective variety $X\subset \PN$.
We
put
$$
\UUU_k (X,\PN,\PNdual):=\set{(L, t)\in \U_k(X,\PN)\times \PNdual}{L\subset H_t},
$$
and consider the projection  
$$
\map{f_{\PNdual}}{\UUU_k (X,\PN,\PNdual)}{\PNdual}.
$$
Then the fiber of $f_{\PNdual}$ over $t\in \PNdual$ is canonically identified with 
$\U_k(Y_t, H_t)$, where $Y_t=X\cap H_t$.
The morphism
$$
\map{f_{\Lambda}}{\UUU_k (X,\PN,\Lambda)}{\Lambda}
$$
defined in Introduction is the pull-back of
$f_{\PNdual}$ by the inclusion $\Lambda\inj \PNdual$.
Consider the following diagram:
$$
\renewcommand{\arraystretch}{1.4}
\begin{array}{ccccc}
\UUU_k (X,\PN,\Lambda) &\inj  & \UUU_k (X,\PN,\PNdual)   &\maprightsp{\pr_1} & \U_k(X,\PN) \\
\lower 3pt \llap{${}^{f_{\Lambda}}$} \mapdown & \square & \mapdown \lower 3pt \rlap{${}^{f_{\PNdual}}$}  & \\
\Lambda&\inj  & \PNdual,   &
\end{array}
$$
where $\pr_1$ is the projection onto the first factor.
The fiber of $\pr_1$ over $L\in \U_k(X,\PN)$
is isomorphic to a linear subspace
$\shortset{t\in\PNdual}{L\subset H_t}$ of $\PNdual$,
and hence $\pr_1$ is smooth and proper (and thus locally trivial) with simply-connected fibers.
Therefore $\UUU_k (X,\PN,\PNdual)$ is smooth and  irreducible, and 
$\pr_1$ induces an isomorphism
\begin{equation}\label{eq:isom1}
\pione(\UUU_k (X,\PN,\PNdual), s_o(0))\cong \pione (\U_k(X,\PN), L_o).
\end{equation}
The fiber of $f_{\PNdual}$ over $t\in \PNdual$ 
is a Zariski open subset of $\Grass^{n-1-k} (H_t)$.
Hence $f_{\PNdual}$ is smooth.
There exists a Zariski closed  subset $\Xi\spprime$
of $\PNdual$ of codimension $\ge 2$ such that, 
if $t\in \PNdual\setminus \Xi\spprime$,
then $Y_t$ has only isolated singular points.
(See~\cite{MR592569}, for example.)
Then $\U_k(Y_t, H_t)$ is non-empty and irreducible for $t\in \PNdual\setminus \Xi\spprime$.
Therefore
$f_{\PNdual}$ satisfies the conditions~\cond{C1} and~\cond{C2}.
In particular,
by Nori's lemma~(Proposition~\ref{prop:nori}),
we see that the inclusion of the general fiber induces a surjective  homomorphism
\begin{equation}\label{eq:iotasurj}
\mapsurj{\iota_*}{ \pione(\U_k (Y_0, H_0), L_o)}{ \pione (\UUU_k (X,\PN,\PNdual), s_o(0))}.
\end{equation}
On the other hand,
in virtue of the \emph{general} line $\Lambda\subset \PNdual$ and the holomorphic
section $s_o$ over $\Lambda$,
we see that $f_{\PNdual}$ satisfies the conditions of  Theorem~\ref{thm:C},
and hence $\iota_*$ induces an injective homomorphism  
\begin{equation}\label{eq:iotainj}
\pione (\U_k (Y_0, H_0), L_o)\ZQ \pione (\Lambda\setminus \Sigma_{\Lambda}, 0)
\;\;\inj\;\; \pione (\UUU_k (X,\PN,\PNdual), s_o(0)).
\end{equation}
Combining~\eqref{eq:isom1},~\eqref{eq:iotasurj} and~\eqref{eq:iotainj},
we complete the proof of   Theorem~\ref{thm:ULZvK}(1).
\par
In particular, the inclusion $\U_k (Y_0, H_0)\inj \U_k (X, \PN)$
induces a surjective homomorphism 
on the fundamental groups.
%
%
%
%
If $k<n-2$, then we can apply this result to the inclusion 
$\U_k(Z_\Lambda, A)\inj \U_k (Y_0, H_0)$,
and obtain a surjection
\begin{equation*}\label{eq:surj2}
\pione(\U_k(Z_\Lambda, A), L_o) \surj\pione( \U_k (Y_0, H_0), L_o).
\end{equation*}
By construction,
this homomorphism is equivariant under the classical monodromy action of 
$\pione (\Lambda\setminus \Sigma_{\Lambda}, 0)$ given by the cross-section $s_o$.
Since $\pione (\Lambda\setminus \Sigma_{\Lambda}, 0)$ acts on $\pione(\U_k(Z_\Lambda, A), L_o)$
trivially,
we obtain   the proof of Theorem~\ref{thm:ULZvK}(2).
\end{proof}
\section{The simple braid group}\label{sec:SB}
Let $C$ be a compact Riemann surface of genus $g>0$,
and let $D_0=p_1+\dots+p_d$ be a reduced effective divisor on $C$ of degree $d$,
which we use as a base point of the space $\rDiv^d(C)$
of reduced divisors of degree $d$ on $C$.
Let $\Pic^d(C)$ be the Picard variety of isomorphism classes $[L]$ of 
line bundles $L$ of degree $d$ on $C$. 
We denote by
$$
\map{\barlambda}{\Div^d (C)}{\Pic^d(C)}
$$
the natural morphism,
and consider the induced homomorphism
$$
\map{\barlambda_*}{\pione(\Div^d (C), D_0)}{\pione(\Pic^d(C), \barlambda(D_0))=H_1(C, \Z)}.
$$
\begin{proposition}\label{prop:barlambda}
\setrmkakko
Suppose that $d\ge g$.
\rmkakko
We have $\Sing(\barlambda)=\barlambda\inv(\barlambda(\Sing(\barlambda)))$.
\rmkakko
If $d\ge 2g-1$ then $\Sing (\barlambda)=\emptyset$.
If $d\le 2g-2$ then $\dim \Sing (\barlambda)\le g-1$
and $\dim \barlambda(\Sing (\barlambda))\le 2g-2-d$.
\end{proposition}
\begin{proof}
Note that $\barlambda$ is surjective because $d\ge g$.
For $D\in \Div^d (C)$, we have
$$
\barlambda\inv (\barlambda(D))=|\OOO_C(D)|\;\;\cong\;\; \P^{d-g+s(D)},
$$
where $s(D):=h^0(C, K_C(-D))$.
Hence $D\in \Sing(\barlambda)$
if and only if $s(D)>0$,
and therefore the assertion (1) follows, and  moreover, we have
$$
\dim \barlambda(\Sing(\barlambda))\le \dim \Sing(\barlambda)-(d-g+1).
$$
On the other hand,
we have $s(D)>0$ if and only if
$D$ is a sub-divisor of a member of the $(g-1)$-dimensional
linear system $|K_C|$.
Since $\deg K_C=2g-2$,
we obtain the proof.
\end{proof}
\begin{remark}
Suppose $d\ge g$.
Then $\Sing(\barlambda)$ is the locus of \emph{special divisors} 
of degree $d$ on $C$,
and $\barlambda(\Sing(\barlambda))$ is the locus of 
\emph{special line bundles} of degree $d$ on $C$.
\end{remark}
\begin{proposition}\label{prop:barlambdastar}
Suppose that $d\ge g$.
Then $\barlambda_*$ is an isomorphism.
\end{proposition}
\begin{proof}
The general fiber of  $\barlambda$ is isomorphic to $\P^{d-g}$.
By Proposition~\ref{prop:barlambda},
the assumption $d\ge g$ implies that
$\barlambda(\Sing(\barlambda))\subset \Pic^d(C)$
is of codimension $\ge 2$.
Hence  Proposition~\ref{prop:barlambdastar}
follows from Nori's lemma~(Proposition~\ref{prop:nori}).
\end{proof}
\begin{proposition}\label{prop:veryample}
{\rm (1)}
Suppose that $d\ge g+2$.
Then there exists a Zariski closed subset $\Xi_1\subset\Pic^d(C)$ of codimension 
 $\ge 2$ such that the complete linear system $|L|$ is base-point free
for any $[L]\in \Pic^d(C)\setminus\Xi_1$. 

{\rm (2)}
Suppose that $d\ge g+4$.
Then there exists a Zariski closed subset $\Xi_2\subset \Pic^d(C)$
of codimension $\ge 2$ such that
$|L|$ is very ample for any $[L]\in \Pic^d(C)\setminus  \Xi_2$.
\end{proposition}
\begin{proof}
Suppose that $d\ge g+2$, and let $L$ be a line bundle of degree $d$.
If  $|L|$  has a base point $p$,
then $L(-p)$ is a special line bundle, and hence
$[L]\in \Pic^{d}(C)$ is contained in the image of the morphism
\begin{equation}\label{eq:timesC1}
\barlambda\sprime(\Sing(\barlambda\sprime))\times C\;\;\to\;\; \Pic^d (C)
\end{equation}
given by $([M], p)\mapsto [M(p)]$,
where $\barlambda\sprime: \Div^{d-1}(C)\to\Pic^{d-1}(C)$ is the natural morphism.
Since  $\dim \barlambda\sprime(\Sing(\barlambda\sprime))\le 2g-d-1$ by Proposition~\ref{prop:barlambda},
the image of~\eqref{eq:timesC1} is of codimension $\ge 2$.
\par
\smallskip
Suppose that $d\ge g+4$.
If a base-point free line bundle $L$ of degree $d$ is not very ample,
then there exist points $p$, $q$ of $C$
such that $h^0(L(-p-q))=h^0(L(-p))$ holds,
and hence $L(-p-q)$ is a special line bundle of degree $d-2$.
We complete the proof by the same argument as above.
\end{proof}
We denote by
$$
\map{\lambda}{\rDiv^d (C)}{\Pic^d(C)}
$$
the restriction of $\barlambda$ to $\rDiv^d(C)$,
and consider the homomorphism 
$$
\map{\lambda_*}{B(C, d):=\pione(\rDiv^d (C), D_0)}{H_1(C, \Z)=\pione(\Pic^d(C))}
$$
induced by $\lambda$.
From Proposition~\ref{prop:barlambdastar},
we obtain the following:
\begin{corollary}\label{cor:SB2}
Suppose that $d\ge g$.
Then the simple braid group $\SB(C, D_0)$ 
defined in Definition~\ref{def:SB} is equal to the kernel of  the homomorphism
$\lambda_*$.
\end{corollary}
Let $\sigma: \IbI\to  (\rDiv^d(C), D_0)$
be a loop.
Then there exist  paths $\sigma_i: I\to C$ for $i=1, \dots, d$ 
such that $\sigma_i(0)=p_i$ and
such that $\sigma(t)=\sigma_1(t)+\dots+\sigma_d(t)$
for all $t\in I$.
The homology class $\lambda_*([\sigma])\in H_1(C, \Z)$
is represented by the $1$-cycle 
obtained as the conjunction of the paths $\sigma_1, \dots, \sigma_d$.
\par
\medskip
Let $\Gamma^d (C) \subset \Div^d (C)$ be the
big diagonal in $\Div^d(C)=C^d/\SSSS_d$,
where $\symgroup_d$ is the symmetric group
acting on the Cartesian product $C^d$ of $d$ copies of $C$
by permutation of the components.
We have
$$
\rDiv^d (C)=\Div^d (C) \setminus \Gamma^d (C).
$$
For $[L]\in \Pic^d (C)$, 
we put
$$
\Gamma (L):=\Gamma^d(C)\cap \barlambda\inv ([L])
\quand
|L|\spred:=\lambda\inv ([L])=|L|\setminus \Gamma (L),
$$
where 
$\barlambda\inv ([L])$ is identified with $|L|$.
\begin{remark}\label{rem:dualhyp}
Suppose that $L$  is very ample,
and let $C_L\subset \P^{d-g+s(L)}$ 
denote the image of the embedding of $C$ by  $|L|$.
Then,
under the identification $|L|\cong (\P^{d-g+s(L)})\dual$,  $\Gamma (L)$ is equal to 
the dual hypersurface $C_L\dual$ of $C_L$,
and hence it is of degree
$$
d\dual:=2(d+g-1).
$$
\end{remark}
\begin{proposition}\label{prop:Lgeneral}
Suppose that $d\ge g+4$.
If $[L]\in \Pic^d(C)$ is general,
then the inclusion $|L|\spred\inj \rDiv^d(C)$ induces an isomorphism
$$
\pione (|L|\spred, D_0)\;\;\cong\;\; \SB(C, D_0),
$$
where $D_0$ is a point of $|L|\spred$.
\end{proposition}
\begin{proof}
We put 
$\Xi:=\barlambda(\Sing(\barlambda))\cup\Xi_2$, 
where $\Xi_2$ is the Zariski closed subset in 
Proposition~\ref{prop:veryample}.
Then $\Xi$ is a Zariski closed subset of codimension $\ge 2$ in $\Pic^d(C)$
and $\barlambda\inv(\Xi)$ is of codimension $\ge 2$ in $\Div^d(C)$
by Proposition~\ref{prop:barlambda}.
Moreover $\barlambda\inv(\Xi)$  contains $\Sing(\barlambda)$,
and $L\sprime$ is very ample if $[L\sprime]\notin \Xi$.
We consider the restriction
$$
\map{f}{X:=\rDiv^d(C)\setminus \lambda\inv (\Xi)}{Y:=\Pic^d(C)\setminus \Xi}
$$
of $\lambda$ to $X=\rDiv^d(C)\setminus \lambda\inv (\Xi)$.
We have
\begin{eqnarray*}
&& \pione(Y, [L])\;\;=\;\;\pione(\Pic^d(C), [L])\;\;=\;\;H_1 (C, \Z), \\
&& \pione(X, D_0)\;\;=\;\;
\pione(\rDiv^d(C), D_0)\;\;=\;\;\B(C, D_0), \\
&&  \pi_2(Y)\;\;=\;\;\pi_2(\Pic^d(C))\;\;=\;\;0.
\end{eqnarray*}
By the last equality, the morphism $f$ satisfies~\cond{Z}.
Since $f$ is smooth  with every fiber being non-empty  Zariski open subsets of
$\P^{d-g}$, 
 the conditions~\cond{C1} and~\cond{C2}
 are also satisfied.
Therefore we can apply  Theorem~\ref{thm:ZvK}.
Using Proposition~\ref{prop:proj} and Remark~\ref{rem:dualhyp}, %
the lifted monodromy action of $\pione (X\spc, D_0)$ on $\pione (|L|^{\red}, D_0)$ is trivial.
Combining this result with Corollary~\ref{cor:RRReq},  
we see that  $\pione (|L|^{\red}, D_0)$  is equal to
the kernel of the homomorphism $\B(C, D_0)\to H_1(C, \Z)$
induced by $f$,
which is $\SB(C, D_0)$ by Corollary~\ref{cor:SB2}. 
\end{proof}
Now we prove our third  main result.
\begin{proof}[Proof of Theorem~\ref{thm:SB}]
We denote by $L$ the line bundle on $C\subset \P^M$
corresponding to the hyperplane section,
and let $C_L\subset \P^N$ be the image of the embedding 
of $C$ by $|L|$.
Then $C\subset \P^M$ is the image of a projection
$C_L\to \P^M$
with the center being disjoint from $C_L\subset \P^N$.
Let $\rho: C\to\Pt$ be a general projection.
By this sequence of the linear projections
$\P^N\ratmap  \P^M\ratmap\P^2$,
we have the canonical embeddings of linear subspaces
$$
(\P^2)\dual\inj (\P^M)\dual\inj (\P^N)\dual.
$$
Let $\rho(C)\dual\subset (\P^2)\dual$,
$C\dual\subset (\P^M)\dual$ and $(C_L)\dual \subset (\P^N)\dual$
be the dual hypersurfaces 
of $\rho(C)\subset \P^2$, $C\subset \P^M$ and $C_L\subset \P^N$,
respectively.
Then we have
$$
\rho(C)\dual = (\P^2)\dual \cap C\dual =(\P^2)\dual \cap (C_L)\dual,
\quad 
C\dual =(\P^M)\dual \cap (C_L)\dual.
$$
We will consider the homomorphisms
$$
\pione ((\P^2)\dual \setminus \rho(C)\dual)
\;\to \;
\pione((\P^M)\dual\setminus C\dual)
\;\to \;
\pione((\P^N)\dual\setminus (C_L)\dual)
$$
induced by the inclusions.
Since $C\subset \P^M$ is Pl\"ucker general by the assumption,
the degree $d\dual$ of $\rho(C)\dual$, the number 
 $\delta\dual$ of ordinary nodes on $\rho(C)\dual$ and 
 the number  $\kappa\dual$ of ordinary cusps on $\rho(C)\dual$
 are given by the Pl\"ucker formula; 
$$
d\dual=2 d+2 g-2,
\quad
\delta\dual=2 d^2+4 d g+2 g^2-10 d-14 g+12,
\quad
\kappa\dual=3 d+6 g-6.
$$
(See~\cite[Chap.~7]{MR2107253}, for example.)
In particular,
the section $\rho(C)\dual$ of $(C_L)\dual$ by $(\P^2)\dual\subset (\P^N)\dual$
 is equisingular to the \emph{general} plane section of $(C_L)\dual$.
 By the classical Zariski hyperplane section 
theorem~(\cite{MR932724},~\cite{MR820315},~\cite{MR1503330}),
we see that the inclusion induces an isomorphism
$$
\pione ((\P^2)\dual \setminus \rho(C)\dual)\;\cong\; \pione((\P^N)\dual\setminus (C_L)\dual).
$$
On the other hand,
the scheme-theoretic intersection of 
$(C_L)\dual$ and  $(\P^2)\dual$ in $(\P^N)\dual$
is reduced, and hence
the scheme-theoretic intersection of 
$C\dual$ and  $(\P^2)\dual$ in $(\P^M)\dual$
is also reduced,
and thus the inclusion induces a surjective homomorphism
$$
\pione ((\P^2)\dual \setminus \rho(C)\dual)\;\surj\;\pione((\P^M)\dual\setminus C\dual).
$$
Therefore we conclude that the inclusions 
induce isomorphisms
$$
\pione ((\P^2)\dual \setminus \rho(C)\dual)
\;\cong \;
\pione((\P^M)\dual\setminus C\dual)
\;\cong \;
\pione((\P^N)\dual\setminus (C_L)\dual).
$$
Remark that $(\P^M)\dual\setminus C\dual$ is equal to $U_0(C, \P^M)$,
and $(\P^N)\dual\setminus (C_L)\dual$ is equal to $|L|^{\red}$.
Therefore it is enough to show that
$\pione (|L|^{\red})$ or $\pione ((\P^2)\dual \setminus \rho(C)\dual)$ is isomorphic to 
the simple braid group $\SB_d^g$.
Note that, since $[L]$ is not necessarily a general point of $\Pic^d(C)$,
we cannot apply Proposition~\ref{prop:Lgeneral}.
We overcome this difficulty using  Harris' theorem~\cite{MR837522}.
\par
\medskip
Note that $\rho (C)$ is a plane curve of degree $d$ with $\delta:=(d-1)(d-2)/2-g$ ordinary nodes
and no other singularities.
Let $\P_*(H^0(\Pt, \OOO(d)))$ be the space of all plane curves of degree $d$, and 
let $\SSS_{d, \delta}\subset\P_*(H^0(\Pt, \OOO(d)))$ be the locus
of reduced plane curves $\Gamma\subset\Pt$ of degree $d$ such that
$\Sing \Gamma$ consists of only $\delta$ ordinary nodes.
In~\cite{MR837522},
Harris gave an affirmative answer to the Severi problem,
in virtue of  which 
we know that $\SSS_{d, \delta}$ is irreducible.
We then denote by $\SSS_{d, \delta}\spc \subset \SSS_{d, \delta}$
the locus of  $\Gamma\in \SSS_{d, \delta}$
such that the dual curve $\Gamma\dual$ has only ordinary nodes and ordinary cusps as its singularities.
Then $\SSS_{d, \delta}\spc$ is a Zariski open subset of $\SSS_{d, \delta}$
containing $\rho(C)$.
\par
\medskip
Let $C\sprime$ be an arbitrary compact Riemann surface of genus $g$,
and let $[L\sprime]$ be a \emph{general} point of $\Pic^d(C\sprime)$.
Since $d\ge g+4$,
we see from Proposition~\ref{prop:veryample} that
$|L\sprime|$ is very ample of dimension $d-g$.
We denote by $C\sprime_{L\sprime}\subset \P\sp{d-g}$ the image of the embedding 
$C\sprime\inj  \P\sp{d-g}$ by $|L\sprime|$,
and consider the general projection $\rho\sprime: C\sprime_{L\sprime}\to \Pt$.
Then  $\rho\sprime(C\sprime_{L\sprime})$ is a point of $\SSS_{d, \delta}$.
Since $\SSS_{d, \delta}$ is irreducible, 
 we can connect the two points
$\rho(C)\in \SSS_{d, \delta}$ and $\rho\sprime(C\sprime_{L\sprime})\in\SSS_{d, \delta}$
by an irreducible closed curve $T\subset \SSS_{d, \delta}$.
We put
$T^0:=T\cap \SSS_{d, \delta}\spc$,
which is a Zariski open subset of $T$ containing $\rho(C)$.
When $\Gamma$ moves on $\SSS_{d, \delta}\spc$
the dual curves $\Gamma\dual$ form an equisingular family of plane curves.
Therefore we have 
\begin{equation}\label{eq:T0}
\pione ((\Pt)\dual\setminus \rho(C)\dual)\;\;\cong\;\;
\pione ((\Pt)\dual\setminus \Gamma\dual)\quad\textrm{for any $\Gamma\in T\sp 0$}.
\end{equation}
On the other hand, by Propositions~\ref{prop:veryample}~and~\ref{prop:Lgeneral},
there exists a Zariski open dense subset $T\sp 1\subset T$
containing $\rho\sprime(C\sprime_{L\sprime})$
such that the complete linear system
$|\OOO_\Gamma (1)|$ of a hyperplane section of $\Gamma\subset\Pt$
is very ample  on the normalization $\Gamma\sp\sim$ of $\Gamma$ for any $\Gamma\in T\sp 1$, 
that $\dim |\OOO_\Gamma (1)|=d-g$ for any $\Gamma\in T\sp 1$, and that
\begin{equation}\label{eq:T1}
\pione ((\Pt)\dual\setminus \Gamma\dual)\;\;\cong\;\;
\pione(|\OOO_\Gamma (1)|\spred)
\;\;\cong\;\;
\SB^d_g 
\quad\textrm{for any $\Gamma\in T\sp 1$}.
\end{equation}
Here we have used the classical Zariski hyperplane section theorem again.
Since $T^0\cap T^1\ne\emptyset$,
we complete  the proof of Theorem~\ref{thm:SB}
by combining the isomorphisms~\eqref{eq:T0},~\eqref{eq:T1}.
\end{proof}
\section{The conjecture of Auroux,  Donaldson, Katzarkov and Yotov}
\label{sec:ADKY}
Let $X\subset \PN$ be a smooth non-degenerate projective surface of degree $d$, and 
let $B\subset \Pt$ be the branch curve of a general projection
$X\to \Pt$.
The fundamental group $\pione (\Pt\setminus B)$
has been studied intensively by 
Moishezon, Teicher and Robb
(\cite{MR644819}, \cite{MR903386}, \cite{MR1203688}, \cite{MR1360512}, 
\cite{MR1689261}, \cite{MR1492521}, \cite{MR1468277}, \dots\dots).
In many examples,
it has turned out that  $\pione (\Pt\setminus B)$ is rather ``small".
In~\cite[Conjectures 1.3 and 1.6]{MR2081427},
Auroux,  Donaldson, Katzarkov and Yotov
formulated  the following conjecture
(not only for algebraic surfaces but also for symplectic $4$-manifolds),
and confirmed it for some new examples.
\par
\medskip
Note that there exist  natural homomorphisms
$$
\pione (\Pt\setminus  B)\;\to\; \SSSS_d\quand\pione (\Pt\setminus B)\;\to\;   H_1 (\Pt\setminus B)\cong \Z/\deg(B)\Z.
$$
For a smooth projective  surface $X$ and a line bundle $L$ on $X$, we denote by
$$
\map{\lambda_{(X,L)}}{H^2(X, \Z)}{\Z^2}
$$
the homomorphism given by $\lambda_{(X, L)}(\alpha):=(\alpha \cup c_1(L), \alpha \cup c_1(K_S+3L))$,
where $\cup $ denotes the cup-product.
\begin{conjecture}\label{conj:ADKY}%
Let $L$ be an ample line bundle  of a smooth projective  surface $S$,
and 
let $X_m\subset \P^{N(m)}$ be the 
 image of the embedding of $S$  by the complete linear system $|L\sp{\otimes m}|$.
 We denote by 
$B_m\subset \Pt$  the branch curve of a general projection 
$X_m\to \Pt$.
Let $G^0_m$ be the kernel of the natural homomorphism
$$
\pione (\Pt\setminus B_m)\;\to\;\SSSS_d\times\Z/\deg(B_m)\Z.
$$
Suppose that $S$ is simply-connected and that $m$ is large enough.
Then the abelianization of $G^0_m$ is isomorphic to
$(\Z^2/\Im(\lambda_{(X, mL)}))^{d-1}$, and the commutator subgroup $[G^0_m, G^0_m]$
is a quotient of $(\Z/2\Z)^2$.
\end{conjecture}
For a smooth non-degenerate projective surface $X\subset\PN$,
the fundamental  groups $\pione (U_0(X, \PN))$ and
 $\pione (\Pt\setminus B)$ are related as follows.
Note that the target space $\Pt$ of the general projection $X\to \Pt$ is 
identified with the closed subvariety
$$
\set{L\in \Grass^{2}(\PN)}{\textrm{$L$ contains the center of the projection}}
$$
of $\Grass^{2}(\PN)$, 
and $\Pt\setminus B$ is identified with the pull-back of $U_0(X, \PN)$
by this embedding $\Pt\inj \Grass^{2}(\PN)$.
\begin{proposition}\label{prop:surj}
The inclusion
$\Pt\setminus B\inj U_0(X, \PN)$ induces 
 a surjective  homomorphism
$\pione (\Pt\setminus B)\surj \pione (U_0(X, \PN))$.
\end{proposition}
\begin{proof}
Consider the incidence variety
$$
\renewcommand{\arraystretch}{1.4}
\begin{array}{ccc}
\set{(L, M)\in \Grass^2(\PN)\times \Grass^3(\PN)}{L\supset M}&\maprightsp{\pr_1} & \Grass^2(\PN)\\
\lower 3pt \llap{${}^{\pr_2}$}\mapdown &&\\
\Grass^3(\PN), 
\end{array}
$$
where $\pr_1$ and $\pr_2$ are the natural projections,
and put
$$
\UUU:=\pr_1\inv (U_0(X, \PN)).
$$
Since $\pr_1$ is smooth with every fiber being isomorphic to $\P^{N-2}$,
we see that $\UUU$ is smooth, irreducible, and that $\pr_1|_{\UUU}$ induces an isomorphism
$\pione (\UUU)\cong \pione (U_0(X, \PN))$.
For  $M\in \Grass^3(\PN)$,
the target space $\Pi_M$ of the projection
$$
\map{\rho_M}{X}{\Pi_M}
$$
with the center $M$ is the fiber of $\pr_2$ over $M$,
and we have
$$
\Pi_M\setminus B_M\cong (\pr_2|_{\UUU})\inv (M)=\pr_2\inv (M)\cap \UUU, 
$$
where $B_M\subset \Pi_M$ is the branch curve of $\rho_M$.
Hence it is enough to show that the inclusion of the general fiber 
of $\pr_2|_{\UUU}$ over $M$ induces a surjective homomorphism
\begin{equation}\label{eq:surjUUU}
\pione ((\pr_2|_{\UUU})\inv (M))\;\surj\;\pione(\UUU).
\end{equation}
Since $\pr_2$ is smooth, so is  $\pr_2|_{\UUU}$.
Moreover the locus of all $M\in \Grass^3(\PN)$ such that $(\pr_2|_{\UUU})\inv (M)=\emptyset$
is contained in a Zariski closed subset of codimension $\ge 2$ in $\Grass^3(\PN)$. 
Hence Nori's lemma~(Proposition~\ref{prop:nori})
implies the surjectivity~\eqref{eq:surjUUU}.
\end{proof}
Thus we see that
the group $\pione (U_0(X, \PN))$ is ``smaller" than $\pione (\Pt\setminus B)$.
In view of Corollary~\ref{cor:SB} and Conjecture~\ref{conj:ADKY},
we expect that the image $\varGamma_{\Lambda}$ of the monodromy~\eqref{eq:monhom}
should be  ``large".
\par\medskip
The group $\varGamma_{\Lambda}$ is generated by the Dehn twists
associated with 
the ordinary nodes of the singular members of the pencil 
$\{Y_t\}_{t\in \Lambda}$.
Hence the group $\varGamma_{\Lambda}$ and its action on $\SB(Y_0, Z_\Lambda)$
can be visualized by drawing on $Y_0$ the reduced divisor  $Z_{\Lambda}$ 
and the vanishing cycles for the singular members of the pencil.
\par\medskip

As for the largeness of $\varGamma_{\Lambda}$,
we have the following  result of Smith~\cite[Theorem 1.3 and Corollary 4.3]{MR1838364}.
\begin{theorem}[Smith]
The vanishing cycles of the Lefschetz fibration
$\YYY\to\Lambda$ fill up the fiber $Y_0$;
that is, their complement is a bunch of discs.
Moreover distinct points of $Z_\Lambda$ are on distinct discs.
\end{theorem}
The second assertion follows from the argument 
in the proof of~\cite[Theorem~5.1]{MR1838364},
and the fact that the homology classes of the sections of $\YYY\to\Lambda$
corresponding to the points of $Z_\Lambda$
are distinct.
\begin{remark}\label{rem:LSXm}
In the calculation of  $\pione (U_0(X_m, \P^{N(m)}))$ by means of  Corollary~\ref{cor:SB}, 
the assumption $d\ge g+4$ 
is satisfied 
when $m$ is large enough.
Indeed,  the degree $d$ of $X_m$ is given by $d=m^2 L^2$,  while 
the genus $g$ of the general hyperplane section $Y_0$ of $X_m$ is given by
$g=(m^2 L^2+mL\cdot K_X)/2 +1$.
\end{remark}
\bibliographystyle{plain}
\def\cprime{$'$} \def\cprime{$'$} \def\cprime{$'$} \def\cprime{$'$}


\begin{thebibliography}{10}

\bibitem{MR2081427}
D.~Auroux, S.~K. Donaldson, L.~Katzarkov, and M.~Yotov.
\newblock Fundamental groups of complements of plane curves and symplectic
  invariants.
\newblock {\em Topology}, 43(6):1285--1318, 2004.

\bibitem{MR0375281}
J.~S. Birman.
\newblock {\em Braids, links, and mapping class groups}.
\newblock Princeton University Press, Princeton, N.J., 1974.
\newblock Annals of Mathematics Studies, No. 82.

\bibitem{MR1682991}
J.~A. Carlson and D.~Toledo.
\newblock Discriminant complements and kernels of monodromy representations.
\newblock {\em Duke Math. J.}, 97(3):621--648, 1999.


\bibitem{MR0366922}
D.~Cheniot.
\newblock Une d\'emonstration du th\'eor\`eme de {Z}ariski sur les sections
  hyperplanes d'une hypersurface projective et du th\'eor\`eme de {V}an
  {K}ampen sur le groupe fondamental du compl\'ementaire d'une courbe
  projective plane.
\newblock {\em Compositio Math.}, 27:141--158, 1973.



\bibitem{MR644816}
I.~Dolgachev and A.~Libgober.
\newblock On the fundamental group of the complement to a discriminant variety.
\newblock In {\em Algebraic geometry (Chicago, Ill., 1980)}, volume 862 of {\em
  Lecture Notes in Math.}, pages 1--25. Springer, Berlin, 1981.



\bibitem{MR932724}
M.~Goresky and R.~MacPherson.
\newblock {\em Stratified {M}orse theory}, volume~14 of {\em Ergebnisse der
  Mathematik und ihrer Grenzgebiete (3)}.
\newblock Springer-Verlag, Berlin, 1988.

\bibitem{MR820315}
H.~A. Hamm and L{\^e}~D{\~u}ng Tr{\'a}ng.
\newblock Lefschetz theorems on quasiprojective varieties.
\newblock {\em Bull. Soc. Math. France}, 113(2):123--142, 1985.

\bibitem{MR837522}
J.~Harris.
\newblock On the {S}everi problem.
\newblock {\em Invent. Math.}, 84(3):445--461, 1986.

\bibitem{KulikovMPIpreprint}
Vik.~S. Kulikov and I.~Shimada.
\newblock On the fundamental groups of complements to the dual hypersurfaces of
  projective curves.
\newblock preprint, MPI 96-32.

\bibitem{MR592569}
K.~Lamotke.
\newblock The topology of complex projective varieties after {S}. {L}efschetz.
\newblock {\em Topology}, 20(1):15--51, 1981.



\bibitem{MR644819}
B. Moishezon.
\newblock Stable branch curves and braid monodromies.
\newblock In {\em Algebraic geometry (Chicago, Ill., 1980)}, volume 862 of {\em
  Lecture Notes in Math.}, pages 107--192. Springer, Berlin, 1981.

\bibitem{MR903386}
B. Moishezon and M.Teicher.
\newblock Simply-connected algebraic surfaces of positive index.
\newblock {\em Invent. Math.}, 89(3):601--643, 1987.

\bibitem{MR1203688}
B.~Moishezon.
\newblock On cuspidal branch curves.
\newblock {\em J. Algebraic Geom.}, 2(2):309--384, 1993.

\bibitem{MR1360512}
B. Moishezon and M.Teicher.
\newblock Fundamental groups of complements of branch curves as solvable
  groups.
\newblock In {\em Proceedings of the Hirzebruch 65 Conference on Algebraic
  Geometry (Ramat Gan, 1993)}, 
  pages 329--345, Ramat Gan, 1996. Bar-Ilan Univ.

\bibitem{MR732347}
M.~V. Nori.
\newblock Zariski's conjecture and related problems.
\newblock {\em Ann. Sci. \'Ecole Norm. Sup. (4)}, 16(2):305--344, 1983.

\bibitem{MR1468277}
A.~Robb.
\newblock On branch curves of algebraic surfaces.
\newblock In {\em Singularities and complex geometry (Beijing, 1994)}, volume~5
  of {\em AMS/IP Stud. Adv. Math.}, pages 193--221. Amer. Math. Soc.,
  Providence, RI, 1997.

\bibitem{MR1282219}
I.~Shimada.
\newblock Remarks on fundamental groups of complements of divisors on algebraic
  varieties.
\newblock {\em Kodai Math. J.}, 17(2):311--319, 1994.

\bibitem{MR1341806}
I.~Shimada.
\newblock Fundamental groups of open algebraic varieties.
\newblock {\em Topology}, 34(3):509--531, 1995.

\bibitem{MR1354002}
I.~Shimada.
\newblock A generalization of {L}efschetz-{Z}ariski theorem on fundamental
  groups of algebraic varieties.
\newblock {\em Internat. J. Math.}, 6(6):921--932, 1995.

\bibitem{MR1421396}
I.~Shimada.
\newblock A note on {Z}ariski pairs.
\newblock {\em Compositio Math.}, 104(2):125--133, 1996.

\bibitem{MR1428061}
I.~Shimada.
\newblock Fundamental groups of complements to singular plane curves.
\newblock {\em Amer. J. Math.}, 119(1):127--157, 1997.

\bibitem{MR1474860}
I.~Shimada.
\newblock On the commutativity of fundamental groups of complements to plane
  curves.
\newblock {\em Math. Proc. Cambridge Philos. Soc.}, 123(1):49--52, 1998.

\bibitem{MR2011641}
I.~Shimada.
\newblock Equisingular families of plane curves with many connected components.
\newblock {\em Vietnam J. Math.}, 31(2):193--205, 2003.


\bibitem{MR1988200}
I.~Shimada.
\newblock Fundamental groups of algebraic fiber spaces.
\newblock {\em Comment. Math. Helv.}, 78(2):335--362, 2003.

\bibitem{MR1952329}
I.~Shimada.
\newblock On the {Z}ariski-van {K}ampen theorem.
\newblock {\em Canad. J. Math.}, 55(1):133--156, 2003.

\bibitem{MR1952330}
I.~Shimada.
\newblock Zariski hyperplane section theorem for {G}rassmannian varieties.
\newblock {\em Canad. J. Math.}, 55(1):157--180, 2003.

\bibitem{MR1838364}
I.~Smith.
\newblock Geometric monodromy and the hyperbolic disc.
\newblock {\em Q. J. Math.}, 52(2):217--228, 2001.

\bibitem{MR1492521}
M.~Teicher.
\newblock Braid groups, algebraic surfaces and fundamental groups of
  complements of branch curves.
\newblock In {\em Algebraic geometry---Santa Cruz 1995}, volume~62 of {\em
  Proc. Sympos. Pure Math.}, pages 127--150. Amer. Math. Soc., Providence, RI,
  1997.

\bibitem{MR1689261}
M.~Teicher.
\newblock The fundamental group of a {${\bf CP}\sp 2$} complement of a branch
  curve as an extension of a solvable group by a symmetric group.
\newblock {\em Math. Ann.}, 314(1):19--38, 1999.

\bibitem{vanKampen}
E.~R. van Kampen.
\newblock On the fundamental group of an algebraic curve.
\newblock {\em Amer.~J.~Math.}, 55:255--260, 1933.

\bibitem{MR2107253}
C.~T.~C. Wall.
\newblock {\em Singular points of plane curves}, volume~63 of {\em London
  Mathematical Society Student Texts}.
\newblock Cambridge University Press, Cambridge, 2004.

\bibitem{MR1503330}
O.~Zariski.
\newblock A theorem on the {P}oincar\'e group of an algebraic hypersurface.
\newblock {\em Ann. of Math. (2)}, 38(1):131--141, 1937.

\end{thebibliography}
%
\end{document}